\documentclass{article}
\usepackage{fullpage}

\newtheorem{theorem}{Theorem}[section]

\newtheorem{proposition}[theorem]{Proposition}

\newtheorem{lemma}[theorem]{Lemma}

\begin{document}

\title{Classification of Normal Operators \\
in Spaces with Indefinite Scalar Product of Rank 2}
\author{{O.V.Holtz \kern0.07\textwidth V.A.Strauss} \\
Department of Applied Mathematics \\
Chelyabinsk State Technical University \\
454080 Chelyabinsk, Russia}
\date{}
\maketitle

\begin{abstract}
\indent A finite-dimensional complex space with indefinite scalar
product $[\cdot \: , \cdot ]$
having $v_{-}=2$ negative squares and $v_{+} \geq 2$ positive ones is considered.
The paper presents a classification of operators that are normal with respect to
this product. It relates to the paper \cite{1}, where the similar classification
was obtained by Gohberg and Reichstein for the case $v=min\{v_{-},v_{+}\}=1$.
\end{abstract}

\section{Introduction}
\indent Consider a complex linear space $C^{n}$ with an indefinite scalar
product $[ \cdot \: , \cdot ]$. By definition, the latter is
a nondegenerate sesquilinear Hermitian form. If the ordinary scalar product
$(\cdot \:, \cdot)$
is fixed, then there exists a nondegenerate Hermitian operator $H$ such that
$[x,y]=(Hx,y)$ $\; \forall x,y\in C^{n}$.
If $A$ is a linear operator ($A: C^{n} \rightarrow C^{n}$), then the
{\em{$H$-adjoint\/}} of $A$ (denoted by $A^{[*]}$) is defined by the identity
$[A^{[*]}x,y]=[x,Ay]$ (hence $A^{[*]}=H^{-1}A^{*}H$).
An operator $N$ is called {\em{$H$-normal\/}} if $NN^{[*]}=N^{[*]}N$, an operator
$U$ is called {\em{$H$-unitary\/}} if $UU^{[*]}=I$, where $I$ is the identity transformation.

Let $V$ be a nontrivial subspace of $C^{n}$. $V$ is called {\em{neutral\/}} if
$[x,y]=0$ for all $x,y \in V$. In this case we may write $[V,V]=0$.
$V$ is called {\em{nondegenerate\/}}
if from $x \in V$ and $\forall y \in V \; [x,y]=0$ it follows
that $x=0$. The subspace $V^{[\perp]}$ is defined as the set of all
vectors $x \in C^{n}$: $[x,y]=0$ $\; \forall y \in V$. If $V$ is
nondegenerate, then $V^{[\perp]}$ is also nondegenerate and
$V \dot{+} V^{[\perp]}=C^{n}$.

A linear operator $A$ acting in $C^{n}$ is called {\em{decomposable\/}} if
there exists a nondegenerate subspace $V \subset C^{n}$ such that both
$V$ and $V^{[\perp]}$ are invariant for $A$. Then $A$ is the {\em{orthogonal
sum of $A_{1}=A|_{V}$ and $A_{2}=A|_{V^{[\perp]}}$}}. Since the conditions
$AV^{[\perp]} \subseteq V^{[\perp]}$ and $A^{[*]}V \subseteq V$ are
equivalent, an operator $A$ is decomposable if there exists a nondegenerate
subspace $V$ which is invariant both for $A$ and $A^{[*]}$.

Pairs of matrices $\{ A_{1}, H_{1}\}$ and $\{A_{2}, H_{2}\}$, where $H_{1}$
and $H_{2}$ are Hermitian, are called {\em{unitarily similar\/}} if
$A_{2}=T^{-1}A_{1}T$, $H_{2}=T^{*}H_{1}T$ for some invertible $T$;
in case when $H_{1}=H_{2}$ they are {\em{$H_{1}$-unitarily similar}}.

Throughout what follows by a rank of a space we mean $v=\min\{v_{-},v_{+}\}$,
where $v_{-}\:$($v_{+}$) is the number of negative (positive)
squares of the quadratic form $[x,x]$, or (it is the same) the number
of negative (positive) eigenvalues of the operator $H$.
Note that without loss of generality it can be assumed that $v_{-} \leq v_{+}$
(otherwise $H$ can be replaced by $-H$; the latter (invertible and Hermitian
operator) has opposite eigenvalues).

Our aim is to obtain a complete classification for $H$-normal operators
acting in the space $C^{n}$ of rank $2$, i.e., to find a set of
canonical forms such that any $H$-normal operator could be reduced to one
and only one of these forms. This means that for any invertible Hermitian
matrix $H$ with $v=2$ and for any $H$-normal matrix $N$ we must
point out one and only one of the canonical pairs of matrices
$\{ \tilde{N}, \tilde{H} \}$ such that the pair $\{N,H\}$ is unitarily similar
to $\{ \tilde{N},\tilde{H} \}$.

Since any $H$-normal operator $N: \: C^{n} \rightarrow C^{n}$ is an orthogonal
sum of $H$-normal operators each of which has one or two distinct eigenvalues
(Lemma~1 from~\cite{1}), it is sufficient to solve our problem only
for indecomposable operators having one or two distinct eigenvalues.

Thus, in this paper we consider only indecomposable operators having one or two
distinct eigenvalues and assume that $2=v_{-} \leq v_{+}$.

Finally let us introduce some notation. Denote the identity matrix of order
$r \times r$ by $I_{r}$, the $r \times r$ matrix with 1's on the secondary
diagonal and zeros elsewhere by $D_{r}$, and a block diagonal matrix with
$A$, $B$, $\ldots$, $C$ diagonal blocks by $A \oplus B \oplus \ldots \oplus C$:
$$
I_{r}=\left( \begin{array}{cccc}
		1 & & & 0    \\
		& \cdot & &   \\
		& & \cdot &   \\
		0 & & & 1  \end{array}
	\right), \;\;\;
D_{r}=\left( \begin{array}{cccc}
		0 & & & 1    \\
		& & \cdot &   \\
		& \cdot & &   \\
		1 & & & 0  \end{array}
	\right), $$
$$A \oplus B \oplus \ldots \oplus C=\left( \begin{array}{cccc}
		A & & & 0    \\
		& B & &  \\
		& & \cdot &   \\
		0 & & & C  \end{array}
	\right). $$

\begin{center}
{\bf{Acknowledgements}}
\end{center}
We would like to express our gratitude to Prof. Heinz Langer
who drew our attention to this problem and to
Prof. Andr\'{e} Ran for his attention to
our work and very helpful comments of this paper.

\section{Some Properties of Indecomposable $H$-normal Operators}
The results of this section hold for any finite-dimensional space with
indefinite scalar product.

\begin{proposition}
   Let an indecomposable $H$-normal operator $N$ acting in $C^{n}$ ($n > 1$)
   have the only eigenvalue $\lambda$; then there exists a decomposition of
   $C^{n}$ into a direct sum of subspaces
   \begin{equation}
   S_{0}=\{ x \in C^{n}: \; (N - \lambda I)x=(N^{[*]}-\overline{\lambda} I)x=0\}, \label{S0}
   \end{equation}
   $S$, $S_{1}$ such that
\begin{equation}
N=\left( \begin{array}{ccc}
        N'=\lambda I & * & * \\
        0 & N_{1} & * \\
        0 & 0 & N''=\lambda I
        \end{array}
        \right), \;
   H=\left( \begin{array}{ccc}
        0 & 0 & I \\
        0 & H_{1} & 0 \\
        I & 0 & 0
      \end{array}
   \right), \label{pred1}
 \end{equation}
where $N':$ $S_{0} \rightarrow S_{0}$, $N_{1}:$ $S \rightarrow S$,
$N'':$ $S_{1} \rightarrow S_{1}$, the internal operator $N_{1}$ is
$H_{1}$-normal, and the pair $\{N_{1},H_{1}\}$ is determined up to the
unitary similarity.
\end{proposition}
{\bf{Proof:}}
Since $N$ and $N^{[*]}$ commute, the subspace $S_{0}$ defined by (\ref{S0})
is nontrivial. For $N$ to be indecomposable $S_{0}$ must be neutral. Indeed,
otherwise $\exists v \in S_{0}: \: Nv= \lambda v,\; N^{[*]}= \overline{\lambda} v,\; [v,v] \neq 0$,
therefore, $V=span\{v\}$ is a nondegenerate subspace that is invariant both
for $N$ and $N^{[*]}$, hence, $N$ is decomposable. Thus, $S_{0}$ is neutral.
Let us take advantage of the following well-known result:
{\em{for any neutral subspace $V_{1} \subset C^{n}$ there exists a
subspace $V_{2}$ ($V_{1} \cap V_{2}=\{ 0 \}$) such that}}
\begin{equation}
H|_{(V_{1}\dot{+}V_{2})}=\left( \begin{array}{cc}
		0 & I  \\
		I & 0
		\end{array}
	\right). \label{S1}
\end{equation}
Therefore, for $S_{0}$ there exists a neutral subspace $S_{1}$ such that
$H|_{(S_{0} \dot{+} S_{1})}$ has form (\ref{S1}). Since the subspace
$(S_{0}\dot{+}S_{1})$ is nondegenerate, the subspace
$S=(S_{0} \dot{+} S_{1})^{[\perp]}$ is also nondegenerate and
$C^{n}=S_{0}\dot{+}S\dot{+}S_{1}$. As $\forall v \in C^{n}$
$(N- \lambda I)v \in S_{0}^{[\perp]}$ and
$(N^{[*]}- \overline{\lambda} I)v \in S_{0}^{[\perp]}$,
the matrices $N$ and $H$ has form (\ref{pred1}) with respect to the
decomposition $C^{n}=S_{0}\dot{+}S\dot{+}S_{1}$. Since $N$ is $H$-normal,
the internal operator $N_{1}$ is $H_{1}$-normal.

It is seen that only the subspace $S_{0}$ is fixed; $S$ and $S_{1}$ may
change. However, the pair $\{N_{1},H_{1}\}$ is unique in a sense, namely,
it is determined up to the unitary similarity. Indeed, any transformation
$T$ such that $TS_{0} \subseteq S_{0}$ has the form
$$ T=\left( \begin{array}{ccc}
        T_{1} & T_{2} & T_{3} \\
        0 & T_{4} & T_{5} \\
        0 & T_{6} & T_{7} \end{array}
   \right). $$
Since
$$ \tilde{H}=\left( \begin{array}{ccc}
        0 & 0 & I \\
        0 & \widetilde{H_{1}} & 0 \\
        I & 0 & 0
      \end{array}
   \right), $$
from condition $\tilde{H}=T^{*}HT$ it follows that $T_{6}=0$,
$\widetilde{H_{1}}=T_{4}^{*}H_{1}T_{4}$. As $\tilde{N}=T^{-1}NT$,
$\widetilde{N_{1}}=T_{4}^{-1}N_{1}T_{4}$ so that the pair $\{N_{1},H_{1}\}$
is unitarily similar to $\{ \widetilde{N_{1}},\widetilde{H_{1}}\}$, Q.E.D.

{\bf{Remark:}} the decomposition $C^{n}=S_{0}\dot{+}S\dot{+}S_{1}$ was
constructed in \cite{1}, section~6 so that the first part of this statement
is borrowed from \cite{1}.

{\bf{Corollary:}} to go over from one decomposition
$C^{n}=S_{0}\dot{+}S\dot{+}S_{1}$ to another by means of a
transformation $T$ it is necessary that $T$ would be block triangular
with respect to both decompositions.

\begin{theorem}
If an $H$-normal operator $N$ acting in a space $C^{n}$ of rank $k \geq 1$
is indecomposable, then either (A) or (B) holds:
\begin{description}
\item{(A)} $N$ has two eigenvalues and $n=2k$;
\item{(B)} $N$ has one eigenvalue and $2k \leq n \leq 4k$.
\end{description}
\end{theorem}
{\bf{Proof:}} First show that $n \geq 2k$. Indeed,
$n=v_{-}+v_{+} \geq 2\: min\{v_{-},v_{+}\}=2k$. Now prove (A). Let $N$ have
two distinct eigenvalues. Then, according to Lemma~1 form~\cite{1}, $C^{n}$
is a direct sum of two neutral subspaces of the same dimension $m$ which are
invariant for $N$ and $N^{[*]}$. Since in a space with indefinite scalar product
no neutral space can be of dimension more than rank of a space, $m \leq k$ and
$n \leq 2k$. But it is established before that $n \geq 2k$. Hence, $n=2k$ and
the proof of (A) is completed.

Now prove (B), i.e., show that if $N$ has one eigenvalue, then $n \leq 4k$.
For $k=1$ the proof is given in Theorem~1, \cite{1}. Suppose inductively that for all
$i \leq k$ the size of indecomposable operators having one eigenvalue is
not more than $4i \times 4i$. Let $v_{-}=k+1$, $v_{+} \geq v_{-}$, $N$
have the only eigenvalue $\lambda$. According to Proposition~1, one can assume
that the matrices $N$ and $H$ has form (\ref{pred1}). Let
$N_{1}=N_{1}^{(1)} \oplus \ldots \oplus N_{1}^{(p)}$ be a decomposition of
the internal operator $N_{1}$ into an orthogonal sum of indecomposable
operators, $H_{1}=H_{1}^{(1)} \oplus \ldots \oplus H_{1}^{(p)}$,
$S=S^{(1)} \oplus \ldots \oplus S^{(p)}$ be the corresponding decompositions
of $H_{1}$ and $S$. Let $v_{-}^{(i)}$ be the number of negative eigenvalues
of $H_{1}^{(i)}$ ($i=1, \ldots p$). If $dim \: S_{0}=s$, then
$\sum_{i=1}^{p}v_{-}^{(i)}=k+1-s$. Let
$$ H_{1}'=\sum_{v_{-}^{(i)}>0} H_{1}^{(i)}, \;\;
 H_{1}''=\sum_{v_{-}^{(i)}=0} H_{1}^{(i)}. $$
Then $H_{1}=H_{1}' \oplus H_{1}''$, $N_{1}=N_{1}' \oplus N_{1}''$,
where $N_{1}'$, $N_{1}''$ are the corresponding sums of operators $N_{1}^{(i)}$.
Since for any $i=1, \ldots p$ rank of the subspace $S_{1}^{(i)}$ is not more
than $v_{-}^{(i)}$,
$v_{-}^{(i)} \leq k$ (because $k+1-s \leq k$), and the size of an indecomposable
operator in a space of rank $0$ is equal to $1$, by the inductive hipothesis
$dim S^{(i)} \leq 4v_{-}^{(i)}$, hence $dim S' \leq 4(k+1-s)$.
Since $H_{1}''$ has only positive eigenvalues, $N_{1}''$ is a usual normal
operator having one eigenvalue $\lambda$, therefore, $N_{1}''=\lambda I$ so that
$$ N=\left( \begin{array}{cccc}
        \lambda I & * & M_{1} & * \\
        0 & N_{1}' & 0 & * \\
        0 & 0 & \lambda I & * \\
        0 & 0 & 0 & \lambda I \end{array}
   \right), \;\;
 N^{[*]}=\left( \begin{array}{cccc}
        \overline{\lambda} I & * & M_{2} & * \\
        0 & N_{1}'^{[*]} & 0 & * \\
        0 & 0 & \overline{\lambda} I & * \\
        0 & 0 & 0 & \overline{\lambda} I \end{array}
   \right).  $$
If $dim \: S''=r > 2s$, then the system
\begin{eqnarray*}
M_{1}X=0 \\
M_{2}X=0
\end{eqnarray*}
has a nontrivial solution $X=(x_{1}, \ldots ,x_{r})^{T}$ (where $Y^{T}$ is $Y$
transposed). Therefore, there exists a nonzero vector $v=\sum_{i=1}^{r}x_{i}w_{i}$
($w_{i}$ are the basis vectors of $S''$) that satisfies the condition
$(N -\lambda I)v=(N^{[*]}- \overline{\lambda} I)v=0$, i.e., $v \in S_{0}$.
But $S_{0} \cap S=\{0\}$. This contradiction proves that $dim\: S'' \leq 2s$.
Thus, $n=2\:dim\:S_{0}+dim\:S'+dim\:S'' \leq 2s+4(k+1-s)+2s=4(k+1)$, Q.E.D.

Since an indecomposable operator cannot have more than two eigenvalues
(Lemma~1, \cite{1}), either (A) or (B) is true so that the proof of the
theorem is completed.

\section{The Classification of Indecomposable $H$-nor\-mal Operators}
The principal aim of this paper is to prove the following result:

\begin{theorem}
If an indecomposable $H$-normal operator $N$ ($N: \: C^{n} \rightarrow C^{n}$)
acts in a space with indefinite scalar product with $v_{-}=2$ negative squares
and $v_{+} \geq 2$ positive ones, then $4 \leq n \leq 8$ and the pair $\{N,H\}$
is unitarily similar to one and only one of canonical pairs
\mbox{\{(\ref{lemma1.1}),(\ref{lemma1.2})\}} -
\mbox{\{(\ref{lemma11.1}),(\ref{lemma11.2})\}}. The choice of the particular
canonical form is determined as follows.

If $N$ has two distinct eigenvalues $\lambda_{1}$, $\lambda_{2}$, then
$\{N,H\}$ is unitarily similar to \mbox{\{(\ref{lemma1.1}),(\ref{lemma1.2})\}}:
$$ N=\left( \begin{array}{cccc}
		 \lambda_{1} & 1    &        0 &            0    \\
		 0 &           \lambda_{1} & 0 &            0    \\
		 0 &           0 &           \lambda_{2} &  0    \\
		 0 &           0 &           x & \lambda_{2}
		\end{array}
		\right),
\; x \in C,
$$
\begin{equation}
 for \; x \neq 0 \;
 \left[ \begin{array}{cc}
       {\cal I}m\{ \lambda_{1}- \lambda_{2} \}>0 & if \;
      {\cal I}m \{ \lambda_{1}- \lambda_{2} \} \neq 0, \\
      {\cal R}e \{ \lambda_{1} - \lambda_{2}\}>0 & otherwise,
      \end{array}  \right. \label{lemma1.1}
\end{equation}
\begin{equation}
   H=\left( \begin{array}{cc}
		 0 & I_{2} \\
		 I_{2} & 0 \end{array}
		\right). \label{lemma1.2}
\end{equation}
If $N$ has one eigenvalue $\lambda$, $dim \: S_{0}=1$, the internal operator
$N_{1}$ is indecomposable, and $n=4$, then $\{N,H\}$ is unitarily similar to
\mbox{\{(\ref{lemma2.1}),(\ref{lemma2.2})\}}:
\begin{equation}
 N=\left( \begin{array}{cccc}
        \lambda & 1 & ir_{1} & ir_{2}z \\
        0 & \lambda & z & 0         \\
        0 & 0 & \lambda & z^{2}     \\
        0 & 0 & 0 & \lambda
        \end{array}
\right), \; |z|=1, \; r_{1},r_{2} \in \Re, \label{lemma2.1}
\end{equation}
\begin{equation}
 H=D_{4}.  \label{lemma2.2}
\end{equation}
If $N$ has one eigenvalue $\lambda$, $dim \: S_{0}=1$, $N_{1}$ is
indecomposable, and $n=5$, then $\{N,H\}$ is unitarily similar to
one and only one of pairs \mbox{\{(\ref{lemma3.1}),(\ref{lemma3.4})\}},
\mbox{\{(\ref{lemma3.2}),(\ref{lemma3.4})\}},
\mbox{\{(\ref{lemma3.3}),(\ref{lemma3.4})\}}:
\begin{equation}
 N=\left( \begin{array}{ccccc}
	\lambda & 1 & 0 & 0 & ir_{3} \\
	0 & \lambda & 1 & ir_{1} & -2r_{1}^{2}+ir_{2}  \\
	0 & 0 & \lambda & 1 & 2ir_{1}     \\
	0 & 0 & 0 & \lambda  & 1 \\
	0 & 0 & 0 & 0 & \lambda
	\end{array}
    \right), \;r_{1},r_{2},r_{3} \in \Re, \label{lemma3.1}
\end{equation}
\begin{equation}
 N=\left( \begin{array}{ccccc}
	\lambda & 1 & 0 & 0 & ir_{3} \\
	0 & \lambda & z & r_{1} & -2z^{2}r_{1}^{2}Im^{2}z+ir_{2}z^{2}  \\
	0 & 0 & \lambda & z & -2ir_{1}z^{2}Imz   \\
	0 & 0 & 0 & \lambda  & z^{2} \\
	0 & 0 & 0 & 0 & \lambda
	\end{array}
        \right), \; \begin{array}{c}
                          |z|=1, \; z \neq i, \\
                           0< arg\:z < \pi, \\
                           r_{1},r_{2},r_{3} \in \Re, \end{array} \label{lemma3.2}
\end{equation}
\begin{equation}
 N=\left( \begin{array}{ccccc}
	\lambda & 1 & 0 & 0 & r_{3} \\
	0 & \lambda & i & r_{1} & 2r_{1}^{2}+ir_{2}  \\
	0 & 0 & \lambda & i & 2ir_{1}     \\
	0 & 0 & 0 & \lambda  & -1 \\
	0 & 0 & 0 & 0 & \lambda
	\end{array}
\right), \; r_{1},r_{2},r_{3} \in \Re,  \label{lemma3.3}
\end{equation}
\begin{equation}
 H=D_{5}.  \label{lemma3.4}
\end{equation}
If $N$ has one eigenvalue $\lambda$, $dim \: S_{0}=1$, $N_{1}$ is
decomposable, and $n=4$, then $\{N,H\}$ is unitarily similar to 
one and only one of pairs \mbox{\{(\ref{lemma4.1}),(\ref{lemma4.4})\}},
\mbox{\{(\ref{lemma4.2}),(\ref{lemma4.4})\}},
\mbox{\{(\ref{lemma4.3}),(\ref{lemma4.4})\}}:
\begin{equation}
 N=\left( \begin{array}{cccc}
	\lambda & 1 & 0 & 0  \\
	0 & \lambda & 0 & z  \\
	0 & 0 & \lambda & 0  \\
	0 & 0 & 0 & \lambda
	\end{array}
    \right), \; |z|=1, \label{lemma4.1}
\end{equation}
\begin{equation}
 N=\left( \begin{array}{cccc}
	\lambda & 1 & 1 & 0 \\
	0 & \lambda & 0 & z \\
	0 & 0 & \lambda & (1+ir)z  \\
	0 & 0 & 0 & \lambda
	\end{array}
\right), \; |z|=1, \; r \in \Re > 0,  \label{lemma4.2}
\end{equation}
\begin{equation}
 N=\left( \begin{array}{ccccc}
	\lambda & 1 & -1 & 0 \\
	0 & \lambda & 0 & z \\
	0 & 0 & \lambda & -(1+ir)z     \\
	0 & 0 & 0 & \lambda
	\end{array}
\right), \; |z|=1, \; r \in \Re > 0, \label{lemma4.3}
\end{equation}
\begin{equation}
H=D_{4}. \label{lemma4.4}
\end{equation}
If $N$ has one eigenvalue $\lambda$, $dim \: S_{0}=1$, $N_{1}$ is
decomposable, and $n=5$, then $\{N,H\}$ is unitarily similar to 
\mbox{\{(\ref{lemma5.1}),(\ref{lemma5.2})\}}:
\begin{equation}
 N=\left( \begin{array}{ccccc}
	\lambda & 1 & 0 & \frac{1}{2}r_{1}^{2}+ir_{2} & 0 \\
	0 & \lambda & 0 & z & 0 \\
	0 & 0 & \lambda & 0 & r_{1}  \\
	0 & 0 & 0 & \lambda & z^{2}  \\
	0 & 0 & 0 & 0 & \lambda	\end{array}
\right), \; |z|=1, \; r_{1},r_{2} \in \Re, \; r_{1}>0, \label{lemma5.1}
\end{equation}
\begin{equation}
 H=D_{5}.  \label{lemma5.2}
\end{equation}
If $N$ has one eigenvalue $\lambda$, $dim \: S_{0}=1$, $N_{1}$ is
decomposable, and $n=6$, then $\{N,H\}$ is unitarily similar to 
one and only one of pairs \mbox{\{(\ref{lemma6.1}),(\ref{lemma6.3})\}},
\mbox{\{(\ref{lemma6.2}),(\ref{lemma6.3})\}}:
\begin{equation}
 N=\left( \begin{array}{cccccc}
	\lambda & 1   &   2ir_{1} & 0      & 0 & 0 \\
	0 &   \lambda &   1       & ir_{1} & 0 & 2r_{1}^{2}-r_{2}^{2}/2+ir_{3} \\
	0 &       0   & \lambda   & 1      & 0 & 0  \\
	0 &       0   & 0     & \lambda    & 0 & 1  \\
        0 &       0   & 0     &     0      & \lambda & r_{2}  \\
	0 &       0   & 0     &     0      & 0  & \lambda \end{array}
\right),\; r_{1},r_{2}\in \Re, \; r_{2}>0, \label{lemma6.1}
\end{equation}
$$ N=\left( \begin{array}{cccccc}
	\lambda & 1   &  -2ir_{1}{\cal I}mz & 0      & 0 & 0 \\
	0 &   \lambda &   z       & r_{1} & 0 & (2r_{1}^{2}{\cal I}m^{2}z-r_{2}^{2}/2+ir_{3})z^{2} \\
	0 &       0   & \lambda   & z      & 0 & 0  \\
	0 &       0   & 0     & \lambda    & 0 & z^{2}  \\
        0 &       0   & 0     &     0      & \lambda & r_{2}  \\
	0 &       0   & 0     &     0      & 0  & \lambda \end{array}
\right), $$
\begin{equation}
|z|=1, \; 0<arg\:z<\pi,\;r_{1},r_{2},r_{3}\in \Re, \; r_{2}>0, \label{lemma6.2}
\end{equation}
\begin{equation}
 H=\left( \begin{array}{cccc}
		0 & 0     & 0 & I_{1} \\
		0 & D_{3} & 0 & 0 \\
		0 & 0     & I_{1} & 0 \\
                I_{1} & 0     & 0 & 0
		\end{array}
	\right).  \label{lemma6.3}
\end{equation}
If $N$ has one eigenvalue $\lambda$, $dim \: S_{0}=2$, and $n=4$, 
then $\{N,H\}$ is unitarily similar to
one and only one of pairs \mbox{\{(\ref{lemma7.1}),(\ref{lemma7.3})\}},
\mbox{\{(\ref{lemma7.2}),(\ref{lemma7.3})\}}:
\begin{equation}
 N=\left( \begin{array}{cccc}
	 \lambda & 0 & z & re^{-i\pi/3}z \\
	 0 & \lambda & 0 & e^{i\pi/3}z  \\
	 0 & 0 & \lambda & 0  \\
	 0 & 0 & 0 & \lambda	\end{array}
\right), \; \begin{array}{c}
|z|=1, \; r \in \Re \geq \sqrt{3}, \\
0\leq arg\:z < \pi \;\; if \; r>\sqrt{3}, \end{array} \label{lemma7.1}
\end{equation}
\begin{equation}
 N=\left( \begin{array}{cccc}
	 \lambda & 0 & 0 & 0 \\
	 0 & \lambda & 1 & 0 \\
	 0 & 0 & \lambda & 0 \\
	 0 & 0 & 0 & \lambda	\end{array}
\right), \label{lemma7.2}
\end{equation}
\begin{equation}
 H=\left( \begin{array}{cc}
		0 & I_{2} \\
		I_{2} & 0 \end{array}
	\right). \label{lemma7.3}
\end{equation}
If $N$ has one eigenvalue $\lambda$, $dim \: S_{0}=2$, and $n=5$,
then $\{N,H\}$ is unitarily similar to
one and only one of pairs \mbox{\{(\ref{lemma8.1}),(\ref{lemma8.3})\}},
\mbox{\{(\ref{lemma8.2}),(\ref{lemma8.3})\}}:
\begin{equation}
 N=\left( \begin{array}{ccccc}
        \lambda & 0 & 1 & 0 & 0 \\
        0 & \lambda & 0 & 1 & 0 \\
        0 & 0 & \lambda & z & 0 \\
        0 & 0 & 0 & \lambda & 0 \\
        0 & 0 & 0 & 0 & \lambda \end{array}
\right), \; |z|=1,  \label{lemma8.1}
\end{equation}
\begin{equation}
 N=\left( \begin{array}{ccccc}
        \lambda & 0 & 1 & 0 & 0 \\
        0 & \lambda & 0 & r & z \\
        0 & 0 & \lambda & z^{2} & 0 \\
        0 & 0 & 0 & \lambda & 0 \\
        0 & 0 & 0 & 0 & \lambda \end{array}
\right), \; |z|=1, \; r \in \Re>0,  \label{lemma8.2}
\end{equation}
\begin{equation}
 H=\left( \begin{array}{ccc}
                0 & 0 & I_{2} \\
                0 & I_{1} & 0 \\
                I_{2} & 0 & 0 \end{array}
        \right).  \label{lemma8.3}
\end{equation}
If $N$ has one eigenvalue $\lambda$, $dim \: S_{0}=2$, and $n=6$,
then $\{N,H\}$ is unitarily similar to
\mbox{\{(\ref{lemma9.1}),(\ref{lemma9.2})\}}:
\begin{equation}
 N=\left( \begin{array}{cccccc}
        \lambda & 0 & 1 & 0 & ir_{1} & 0 \\
	0 & \lambda & 0 & 1 & r_{2} & ir_{1} \\
	0 & 0 & \lambda & 0 & z & 0 \\
	0 & 0 & 0 & \lambda & 0 & z \\
	0 & 0 & 0 & 0 & \lambda & 0 \\
	0 & 0 & 0 & 0 & 0 & \lambda	\end{array}
\right), \; \begin{array}{c} |z|=1,\;z \neq -1,
                      \\ r_{1},r_{2} \in \Re,\;r_{2}>0,
             \end{array}  \label{lemma9.1}
\end{equation}
\begin{equation}
 H=\left( \begin{array}{ccc}
                0 & 0 & I_{2} \\
		0 & I_{2} & 0 \\
		I_{2} & 0 & 0 \end{array}
	\right).  \label{lemma9.2}
\end{equation}
If $N$ has one eigenvalue $\lambda$, $dim \: S_{0}=2$, and $n=7$, 
then $\{N,H\}$ is unitarily similar to 
\mbox{\{(\ref{lemma10.1}),(\ref{lemma10.2})\}}:
$$N=\left( \begin{array}{ccccccc}
        \lambda & 0 & 1 & 0 & 0 & 0 & 0 \\
        0 & \lambda & 0 & 1 & 0 & 0 & 0 \\
        0 & 0 & \lambda & 0 & 0 & -z_{1} \overline{z_{2}}cos \alpha & sin \alpha cos \beta \\
        0 & 0 & 0 & \lambda & 0 & z_{1} sin \alpha & z_{2}cos \alpha cos \beta  \\
        0 & 0 & 0 & 0 & \lambda & 0 & sin \beta \\
        0 & 0 & 0 & 0 & 0 & \lambda & 0 \\
        0 & 0 & 0 & 0 & 0 & 0 & \lambda \end{array}
        \right), $$
\begin{eqnarray}
& |z_{1}|=|z_{2}|=1, \;\; 0<\alpha, \beta\leq \pi/2, & \nonumber \\
& z_{1}=1 \;\; if \; \beta=\pi/2, \;\; z_{2}=1 \;\; if \; \alpha=\pi/2, & \label{lemma10.1}
\end{eqnarray}
\begin{equation}
H=\left( \begin{array}{ccc}
        0 & 0 & I_{2} \\
        0 & I_{3} & 0 \\
        I_{2} & 0 & 0 \end{array}
        \right).  \label{lemma10.2}
\end{equation}
If $N$ has one eigenvalue $\lambda$, $dim \: S_{0}=2$, and $n=8$,
then $\{N,H\}$ is unitarily similar to
\mbox{\{(\ref{lemma11.1}),(\ref{lemma11.2})\}}:
$$N- \lambda I=\left( \begin{array}{cccccccc}
        0 & 0 & 1 & 0 & 0 & 0 & 0 & 0 \\
        0 & 0 & 0 & 1 & 0 & 0 & 0 & 0 \\
        0 & 0 & 0 & 0 & 0 & 0 & -z_{1}\overline{z_{2}}sin \alpha \/cos \beta & cos \alpha \/ cos \gamma \\
        0 & 0 & 0 & 0 & 0 & 0 & z_{1}cos \alpha\/ cos \beta & z_{2}sin \alpha \/ cos \gamma \\
        0 & 0 & 0 & 0 & 0 & 0 & \sin \beta & 0 \\
        0 & 0 & 0 & 0 & 0 & 0 & 0 & \sin \gamma \\
        0 & 0 & 0 & 0 & 0 & 0 & 0 & 0 \\
        0 & 0 & 0 & 0 & 0 & 0 & 0 & 0 \end{array}
        \right),$$
\begin{eqnarray}
& |z_{1}|=|z_{2}|=1, \; 0 \leq \alpha< \pi/2, \; 0< \beta< \gamma \leq \pi/2, & \nonumber \\
& z_{1}=1 \;\; if \; \gamma=\pi/2, \;\; z_{2}=1 \;\; if \; \alpha=0 & \label{lemma11.1}
\end{eqnarray}
\begin{equation}
H=\left( \begin{array}{ccc}
        0 & 0 & I_{2} \\
        0 & I_{4} & 0 \\
        I_{2} & 0 & 0 \end{array}
        \right).  \label{lemma11.2}
\end{equation}
\end{theorem}

The following sections contain the proof of this theorem.
\section{Two Distinct Eigenvalues of $N$}
\indent Suppose an indecomposable $H$-normal operator $N$ has $2$ distinct
eigenvalues. Then (Lemma~1, \cite{1}) $C^{n}={\cal Q}_{1}\dot{+}{\cal Q}_{2}$,
$dim\,{\cal Q}_{1}=dim\,{\cal Q}_{2}=m$,
$[{\cal Q}_{1},{\cal Q}_{1}]=0$, $[{\cal Q}_{2},{\cal Q}_{2}]=0$,
$N{\cal Q}_{1} \subseteq {\cal Q}_{1}$,
$N{\cal Q}_{2} \subseteq {\cal Q}_{2}$,
$N_{1}=N|_{{\cal Q}_{1}}$ $(N_{2}=N|_{{\cal Q}_{2}})$ has only one eigenvalue
$\lambda_{1}$ ($\lambda_{2}$).
According to Theorem~1, $m=2$ and $n=4$.
Note that the subspaces ${\cal Q}_{1}$ and ${\cal Q}_{2}$
are determined up to interchanging. Since $N$ is indecomposable, at least one
of the operators $N_{1}$, $N_{2}$ is not scalar. Consequently, one can assume
$N_{1} \neq \lambda_{1} I$. If both $N_{1}$ and $N_{2}$ are not scalar, then
we can fix ${\cal I}m\{ \lambda_{1}- \lambda_{2} \}>0$ if
${\cal I}m \{ \lambda_{1}- \lambda_{2}\} \neq 0$ and
${\cal R}e \{ \lambda_{1} - \lambda_{2}\}>0$ if
${\cal I}m \{ \lambda_{1}- \lambda_{2} \}=0$
(let us remember that $\lambda_{1} \neq \lambda_{2}$). Now ${\cal Q}_{1}$ and
${\cal Q}_{2}$ are determined uniquely.

As $H$ is nondegenerate, for any basis in ${\cal Q}_{1}$ there exists a
basis in ${\cal Q}_{2}$ such that
$$H=\left( \begin{array}{cc}
		0 & I  \\
		I & 0
		\end{array}
	\right). $$
Let us fix a basis in ${\cal Q}_{1}$ such that
\begin{equation}
N_{1}=\left( \begin{array}{cc}
		\lambda_{1} & 1  \\
		0 & \lambda_{1}
		\end{array}
	\right). \label{JF}
\end{equation}
$N$ is $H$-normal if and only if
\begin{equation}
	N_{1}{N_{2}}^{*}={N_{2}}^{*}N_{1}.           \label{1}
\end{equation}
From (\ref{1}) it follows that $N_{2}^{*}=\alpha N_{1}+ \beta I$.
As $N_{2}=\overline{\alpha} N_{1}^{*}+\overline{\beta} I$ has the only
eigenvalue $\lambda_{2}$, we conclude
$N_{2}=\lambda_{2}I+x(N_{1}^{*}-\overline{\lambda_{1}}I)$ ($x \in C$).
Thus, we have reduced $N$ to the form
\begin{equation}
N=\left( \begin{array}{cc}
                \lambda_{1} & 1 \\
                0 & \lambda_{1} \end{array} \right) \oplus
   \left( \begin{array}{cc}
		\lambda_{2} & 0  \\
		x & \lambda_{2}	\end{array}
	\right), \; x \in C. \label{TECF}
\end{equation}
Show that forms (\ref{TECF}) with different values of $x$ are not $H$-unitarily
Jsimilar. To this end suppose that some matrix $T$ satisfies the conditions
\begin{eqnarray}
NT=T\tilde{N}, \label{11}  \\
TT^{[*]}=I, \label{12}
\end{eqnarray}
where $N=N_{1}\oplus N_{2}$, $\tilde{N}=N_{1} \oplus \widetilde{N_{2}}$,
$N_{1}$ has form (\ref{JF}),
$$N_{2}=\left( \begin{array}{cc}
		\lambda_{2} & 0  \\
		x & \lambda_{2}
		\end{array}
	\right), \;\;
 \widetilde{N_{2}}=\left( \begin{array}{cc}
		\lambda_{2} & 0  \\
		\tilde{x} & \lambda_{2}
		\end{array}
	\right). $$
From (\ref{11}) it follows that $T$ is block diagonal with respect to the
decomposition $C^{n}={\cal Q}_{1}\dot{+}{\cal Q}_{2}$: $T=T_{1} \oplus T_{2}$,
$T_{1}$ satisfying the condition $N_{1}=T_{1}^{-1}N_{1}T_{1}$.
Taking into account (\ref{12}), we get $T_{2}=T_{1}^{*-1}$, therefore,
$\widetilde{N_{2}}=T_{2}^{-1}N_{2}T_{2}=N_{2}$, i.e., $\tilde{x}=x$.

It can easily be checked that (\ref{TECF}) is indecomposable so that we have
proved the following lemma:
\begin{lemma}
   If an indecomposable $H$-normal operator acts in a space $C^{n}$
   of rank 2 and has 2 distinct
   eigenvalues $\lambda_{1}$, $\lambda_{2}$,
   then $n=4$ and the pair $\{N,H\}$ is unitarily similar
   to canonical pair \{(\ref{lemma1.1}),(\ref{lemma1.2})\}:
$$ N=\left( \begin{array}{cccc}
		 \lambda_{1} & 1    &        0 &            0    \\
		 0 &           \lambda_{1} & 0 &            0    \\
		 0 &           0 &           \lambda_{2} &  0    \\
		 0 &           0 &           x & \lambda_{2}
		\end{array}
		\right),
\; x \in C, $$
$$ for \; x \neq 0 \;
 \left[ \begin{array}{cc}
       {\cal I}m\{ \lambda_{1}- {\lambda_{2}} \}>0 & if \;
      {\cal I}m \{ \lambda_{1}- {\lambda_{2}}\} \neq 0, \\
      {\cal R}e \{ \lambda_{1} - \lambda_{2}\}>0 & otherwise,
      \end{array}  \right.
$$
$$ H=\left( \begin{array}{cc}
		 0 & I_{2} \\
		 I_{2} & 0 \end{array}
		\right),
$$
where the number $x$ forms a complete and minimal invariant of
the pair $\{N,H\}$ under the unitary similarity (in short, we say that $x$
is an $H$-unitary invariant). In other words, every pair $\{N,H\}$
satisfying the hypothesis of the lemma is unitary similar to pair
\{(\ref{lemma1.1}),(\ref{lemma1.2})\} and pairs
\{(\ref{lemma1.1}),(\ref{lemma1.2})\} with different values of
$x$ are not $H$-unitarily similar to each other.
\end{lemma}

\section{One Eigenvalue of $N$}
\indent Throughout what follows we will assume that $N$ has
only one eigenvalue $\lambda$ so that $N$ and $H$ have form (\ref{pred1}).
Since the neutral subspace $S_{0}$ cannot be more than two-dimensional,
there appear two cases to be considered: $dim\:S_{0}=1$ and $dim \: S_{0}=2$.
Now let us prove the following proposition which holds for all spaces with
indefinite scalar product:

\begin{proposition}
An $H$-normal operator such that $dim \:S_{0}=1$ is indecomposable.
\end{proposition}
{\bf{Proof:}} Assume the converse.
Suppose some nondegenerate subspace $V$ is invariant both for $N$ and for
$N^{[*]}$. Let us denote $V_{1}=V$, $V_{2}=V^{[\perp]}$, $N_{1}=N|_{V_{1}}$,
$N_{2}=N|_{V_{2}}$, $H_{1}=H|_{V_{1}}$, $H_{2}=H|_{V_{2}}$. The
following conditions must hold:
$N_{1}N_{1}^{[*]}=N_{1}^{[*]}N_{1}$,
$N_{2}N_{2}^{[*]}=N_{2}^{[*]}N_{2}$. Here
$N_{i}^{[*]}$ is the $H_{i}$-adjoint of $N_{i}$ ($i=1,2$). Let us define
$$
S_{0}^{i}=\{ x \in V_{i}: \: (N_{i}- \lambda I)x=(N_{i}^{[*]}- \overline{\lambda} I)x=0 \}, \;\;i=1,2.
$$
Since the operators $N_{1}$ and $N_{1}^{[*]}$ ($N_{2}$ and $N_{2}^{[*]}$)
commute, $dim\:S_{0}^{i} \geq 1$ ($i=1,2$), therefore,
$dim \{ S_{0}=S_{0}^{1}+S_{0}^{2} \} \geq 2$. This contradicts
the condition $dim\:S_{0}=1$. Thus, $N$ is indecomposable.

If $dim S_{0}=1$, then rank of $S$ is equal to $1$, therefore, to classify
the internal operator $N_{1}$ we may apply Theorem~1 from~\cite{1}.
Since the indecomposability (or decomposability) of $N_{1}$ is a property
which does not change under the unitary similarity of the pair $\{N_{1},H_{1}\}$,
we must consider both the case when $N_{1}$ is indecomposable and that when
$N_{1}$ is decomposable.

\subsection{$dim \: S_{0}=1$ and $N_{1}$ is Indecomposable}
If $N_{1}$ is indecomposable, then, according to Theorem~1,
$ 2 \leq dim \: S \leq 4$ (recall that rank of $S$ is equal to $1$).
Therefore, $4 \leq n \leq 6$. Let us
consider the alternatives $n=4,5,6$ one after another.

\subsubsection{$n=4$}
According to Theorem~1 of~\cite{1}, one can assume that
$N_{1}$ and $H_{1}$ are reduced to the form
$$ N_{1}=\left( \begin{array}{cc}
                        \lambda & z \\
                        0 & \lambda
                        \end{array}
                  \right), \; | z | = 1, \;\;
H_{1}=D_{2}.$$
Hence
$$N- \lambda I=\left( \begin{array}{cccc}
                 0 & a & b & c \\
                 0 & 0 & z & d \\
                 0 & 0 & 0 & e \\
                 0 & 0 & 0 & 0
                 \end{array}
        \right), \;\;
 H=D_{4}. $$
Throughout what follows only $H$-unitary transformations are used
unless otherwise stipulated. This means that for each case we fix some form
of the matrix $H$ and find out to what form it is possible to reduce the
matrix $N$ without the change of $H$.

The condition of the $H$-normality of $N$ is equivalent to the system
\begin{eqnarray}
  a \overline{z} & = & \overline{e}z \label{sys1} \\
  {\cal R}e\{a \overline{b}\} & = & {\cal R}e\{d \overline{e}\}. \label{sys2}
\end{eqnarray}
If $a=0$, then $e=0$, therefore, the vector $v_{2}$ from $S$ ($v_{i}$ are the
basis vectors) belongs to $S_{0}$, which is impossible.
Thus, $a \neq 0$. Replace the vector $v_{1}$
by $av_{1}$ and $v_{4}$ by $v_{4}/\overline{a}$. This transformation
reduces $N - \lambda I$ to the form
 $$N- \lambda I=\left( \begin{array}{cccc}
                0 & 1 & b' & c' \\
                0 & 0 & z & d' \\
                0 & 0 & 0 & z^{2} \\
                0 & 0 & 0 & 0
                \end{array}
        \right).  $$
Further, apply the transformation
 $$T=\left( \begin{array}{cccc}
        1 & z\overline{d}\:' & 0 & 0 \\
        0 & 1 & 0 & 0 \\
        0 & 0 & 1 & -\overline{z}d'\\
        0 & 0 & 0 & 1
        \end{array}
\right) $$
to the matrix $N - \lambda I$. We obtain:
 $$ N- \lambda I=\left( \begin{array}{cccc}
                0 & 1 & b'' & c'' \\
                0 & 0 & z & 0 \\
                0 & 0 & 0 & z^{2} \\
                0 & 0 & 0 & 0
                \end{array}
        \right). $$
It follows from (\ref{sys2}) that $b''=ir_{1}$ ($r_{1} \in \Re$).
Taking the transformation
$$ T=\left( \begin{array}{cccc}
    1 & 0 & \frac{1}{2}\overline{z}{\cal R}e\{c''\overline{z}\} & 0 \\
    0 & 1 & 0 & -\frac{1}{2}z{\cal R}e\{c''\overline{z}\} \\
    0 & 0 & 1 & 0  \\
    0 & 0 & 0 & 1
    \end{array}
  \right), $$
we reduce $N - \lambda I$ to the final form
\begin{equation}
 N- \lambda I=\left( \begin{array}{cccc}
                0 & 1 & ir_{1} & ir_{2}z \\
                0 & 0 & z & 0         \\
                0 & 0 & 0 & z^{2}     \\
                0 & 0 & 0 & 0
                \end{array}
        \right), \; |z|=1, \; r_{1},r_{2} \in \Re, \label{N2}
\end{equation}
where $r_{2}={\cal I}m\{c''\overline{z}\}$.

Let us prove that the numbers $z$, $r_{1}$, $r_{2}$ are $H$-unitary
invariants. Indeed, let $T$ be an $H$-unitary transformation
of the matrix $N$ to the form $\tilde{N}$, where
$$ \tilde{N}- \lambda I=\left( \begin{array}{cccc}
                        0 & 1 & i\widetilde{r_{1}} & i\widetilde{r_{2}}\tilde{z} \\
                        0 & 0 & \tilde{z} & 0         \\
                        0 & 0 & 0 & \tilde{z}^{2}     \\
                        0 & 0 & 0 & 0
                \end{array}
\right), \; |\tilde{z}|=1, \; \widetilde{r_{1}}, \widetilde{r_{2}} \in \Re. $$
This means that $T$ satisfies conditions (\ref{11}) and (\ref{12}).
From Corollary of Proposition~1 it follows that $T$ is block triangular with
respect to the decomposition $C^{n}=S_{0}+S+S_{1}$. According to Theorem~1
from \cite{1}, $z$ is an $H_{1}$-unitary invariant of $N_{1}$. $T_{4}=T|_{S}$ is
a $H_{1}$-unitary transformation of $N_{1}$ to the form $\widetilde{N_{1}}$,
therefore, $z$ is also an $H$-unitary invariant of $N$, i.e., $\tilde{z}=z$.
Applying condition
(\ref{11}), we see that $T$ is uppertriangular and its diagonal terms
are equal to each other. From (\ref{12}) it follows that $|t_{11}|=1$.
Therefore, without loss of generality one can assume that $t_{11}=1$
(we replace our matrix $T$ by the matrix $T'=\overline{t_{11}}T$; the latter
has the same properties (\ref{11}), (\ref{12}) ).

Thus,
$$ T=\left( \begin{array}{cccc}
                1 & t_{12} & t_{13} & t_{14} \\
                0 & 1 & t_{23} & t_{24} \\
                0 & 0 & 1 & t_{34}  \\
                0 & 0 & 0 & 1
                \end{array} \label{T}
        \right). $$
For $T$ to be $H$-unitary it is neccessary and sufficient to have
\begin{eqnarray}
   \overline{t_{34}}+t_{12} & = & 0 \label{n_sys1} \\
   \overline{t_{24}}+\overline{t_{23}}t_{12}+t_{13} & = & 0 \label{n_sys2} \\
   {\cal R}e t_{14}+{\cal R}e \{ t_{12} \overline{t_{13}} \} & = & 0 \label{n_sys3} \\
   {\cal R}e t_{23} & = & 0, \label{n_sys4}
\end{eqnarray}
for $T$ to reduce $N$ to the form $\tilde{N}$
it is neccessary and sufficient to have
\begin{eqnarray}
t_{23}+ir_{1} & = & i \widetilde{r_{1}}+zt_{12}  \label{m_sys11} \\
t_{24}+ir_{1}t_{34}+ir_{2}z & = & i\widetilde{r_{2}}z+z^{2}t_{13} \label{m_sys22} \\
zt_{34} & = & z^{2}t_{23}. \label{m_sys33}
\end{eqnarray}
Express $t_{34}$ in terms of $t_{23}$ from (\ref{m_sys33}) and $t_{12}$
in terms of $t_{23}$ from (\ref{m_sys11}):
$t_{34}=zt_{23},\;t_{12}=\overline{z}(ir_{1}-i\widetilde{r_{1}})+\overline{z}t_{23}$.
Substituting these expressions in (\ref{n_sys1}), we get:
$2 {\cal R}e t_{23}=i(\widetilde{r_{1}}-r_{1})$. Since ${\cal R}e t_{23}=0$
(condition (\ref{n_sys4})), $\widetilde{r_{1}}=r_{1}$. Further, let us express $t_{24}$ in terms of
$t_{13}$ and $t_{23}$ (condition (\ref{m_sys22})):
$t_{24}=(i \widetilde{r_{2}}-ir_{2})z+z^{2}t_{13}-ir_{1}zt_{23}$.
Then condition (\ref{n_sys2}) can be written in the form
$$(ir_{2}-i\widetilde{r_{2}})+\overline{zt_{13}}+zt_{13}+ir_{1}\overline{t_{23}}+ | t_{23} |^{2}=0.$$
As ${\cal R}et_{23}=0$, $ir_{1}\overline{t_{23}} \in \Re $, consequently,
$\overline{zt_{13}}+zt_{13}+ir_{1}\overline{t_{23}}+ |t_{23}|^{2} \in \Re $.
But $i(r_{2}-\widetilde{r_{2}}) \in \Im $. Therefore, $\widetilde{r_{2}}=r_{2}$.
Thus, the numbers $z$, $r_{1}$, $r_{2}$ are $H$-unitary invariants.

Due to Proposition~2 matrix (\ref{N2}) is indecomposable so that we have proved
 the following lemma:
\begin{lemma}
If an indecomposable $H$-normal operator $N$ ($N: \: C^{4} \rightarrow C^{4}$)
has the only eigenvalue $\lambda$, $dim \:S_{0}=1$,
the internal operator $N_{1}$ is indecomposable, then the pair $\{N,H\}$ is
unitarily similar to canonical pair \{(\ref{lemma2.1}),(\ref{lemma2.2})\}:
$$N=\left( \begin{array}{cccc}
        \lambda & 1 & ir_{1} & ir_{2}z \\
        0 & \lambda & z & 0         \\
        0 & 0 & \lambda & z^{2}     \\
        0 & 0 & 0 & \lambda
        \end{array}
\right), \; |z|=1, \; r_{1},r_{2} \in \Re, $$
$$ H=D_{4}, $$
where $z$, $r_{1}$, $r_{2}$ are $H$-unitary invariants.
\end{lemma}

\subsubsection{$n=5$}
 According to Theorem~1 of~\cite{1}, it can be assumed that
the pair $\{N_{1},H_{1}\}$ has either form (\ref{th1.3}) or (\ref{th1.4}):
\begin{equation}
N_{1}=\left( \begin{array}{ccc}
	         	\lambda & z & r \\
          		0 & \lambda & z \\
                        0 & 0 & \lambda
         		\end{array}
		  \right), \; | z | = 1, \; 0<arg\:z< \pi,\; r \in \Re, \;\;H_{1}=D_{3},  \label{th1.3}
\end{equation}
\begin{equation}
N_{1}=\left( \begin{array}{ccc}
	         	\lambda & 1 & ir \\
          		0 & \lambda & 1  \\
                        0 & 0 & \lambda
         		\end{array}
	                \right), \;\; r \in \Re, \;\; H_{1}=D_{3}. \label{th1.4}
\end{equation}
For a while we will consider both the cases together, assuming that
$$N_{1}=\left( \begin{array}{ccc}
	         	\lambda & z' & x \\
          		0 & \lambda & z' \\
                        0 & 0 & \lambda
         		\end{array}
		        \right), \; | z'| = 1, \; 0 \leq arg\:z'< \pi, \; x \in C.$$
Then
$$N - \lambda I=\left( \begin{array}{ccccc}
	    0 & a & b & c & d   \\
            0 & 0 & z' & x & e \\
            0 & 0 & 0 & z' & f  \\
            0 & 0 & 0 & 0 & g   \\
            0 & 0 & 0 & 0 & 0 	\end{array}
		  \right).$$
The condition of the $H$-normality is equivalent to the system
\begin{eqnarray}
  a \overline{z'} & = & \overline{g}z' \label{sys51} \\
 a\overline{x}+b\overline{z'} & = & \overline{g}x+\overline{f}z' \label{sys52}  \\
2{\cal R}e\{a \overline{c}\}+ | b |^{2} & = & 2{\cal R}e\{e \overline{g}\}+ | f |^{2}. \label{sys53}
\end{eqnarray}
As above (see the case when $n=4$), one can check that $a \neq 0$, hence
$a$ can be assumed equal to $1$, so $g=z'^{2}$.
Having in mind these equalities, take the ($H$-unitary) transformation
$$T=\left( \begin{array}{ccccc}
1 & \overline{z'}b & \overline{z'}(c-x\overline{z'}b) & 0 & -\frac{1}{2} | c-x \overline{z'}b|^{2}   \\
0 & 1 & 0 & 0 & 0 \\
0 & 0 & 1 & 0 & -z'(\overline{c}-\overline{x}z' \overline{b})  \\
0 & 0 & 0 & 1 & -z'\overline{b}   \\
0 & 0 & 0 & 0 & 1   	\end{array}
 \right).$$
It reduces $N -\lambda I$ to the form
$$
N - \lambda I=\left( \begin{array}{ccccc}
	    0 & 1 & 0 & 0 & d'   \\
            0 & 0 & z' & x & e' \\
            0 & 0 & 0 & z' & f'  \\
            0 & 0 & 0 & 0 & z'^{2}   \\
            0 & 0 & 0 & 0 & 0
         	\end{array}
		  \right).
$$
Now apply either the transformation
$$T =\left( \begin{array}{ccccc}
	    1 & 0 & 0 & {\cal R}e d'/({\cal R}e z'^{2}+1)  & 0   \\
            0 & 1 & 0 & 0 & -{\cal R}e d'/({\cal R}e z'^{2}+1) \\
            0 & 0 & 1 & 0 & 0  \\
            0 & 0 & 0 & 1 & 0   \\
            0 & 0 & 0 & 0 & 1
         	\end{array}
		  \right) \; (z' \neq i) $$
or
$$T =\left( \begin{array}{ccccc}
	       	1 & 0 & 0 & -\frac{1}{2}i{\cal I}m d' & 0   \\
         	0 & 1 & 0 & 0 & -\frac{1}{2}i{\cal I}m d' \\
            0 & 0 & 1 & 0 & 0  \\
            0 & 0 & 0 & 1 & 0   \\
            0 & 0 & 0 & 0 & 1
         	\end{array}
		  \right) \; (z'=i)$$
to the matrix $N - \lambda I$. We get
$$N- \lambda I =\left( \begin{array}{ccccc}
    0 & 1 & 0 & 0 & i({\cal I}m d'+{\cal I}m\{d'\overline{z'}^{2}\})/({\cal R}e z'^{2}+1)  \\
    0 & 0 & z' & x & e' \\
    0 & 0 & 0 & z' & f'  \\
    0 & 0 & 0 & 0 & z'^{2}   \\
    0 & 0 & 0 & 0 & 0 	\end{array}
  \right) \; (z' \neq i)$$
or
$$N- \lambda I =\left( \begin{array}{ccccc}
	       	0 & 1 & 0 & 0 & {\cal R}e d'  \\
         	0 & 0 & i & x & e' \\
                0 & 0 & 0 & i & f'  \\
                0 & 0 & 0 & 0 & -1   \\
                0 & 0 & 0 & 0 & 0
         	\end{array}  \right) \; (z'=i). $$

Now we shall distinguish cases (\ref{th1.3}) and (\ref{th1.4}).

(a) ${z'=1, \;\; x=ir_{1} \; (r_{1} \in \Re).\:}$  Conditions (\ref{sys52}), (\ref{sys53}) of the
$H$-normality of $N$ yield:
$f'=2ir_{1}$, $e'=-2r_{1}^{2}+ir_{2}$. Denote $({\cal I}md'-{\cal I}m\{d'\overline{z'}^{2}\})/({\cal R}ez'^{2}+1)$
by $r_{3}$.  We have
\begin{equation}
N- \lambda I =\left( \begin{array}{ccccc}
            0 & 1 & 0 & 0 & ir_{3}   \\
            0 & 0 & 1 & ir_{1} & -2r_{1}^{2}+ir_{2} \\
            0 & 0 & 0 & 1 & 2ir_{1}  \\
            0 & 0 & 0 & 0 & 1   \\
            0 & 0 & 0 & 0 & 0
         	\end{array}
		  \right), \;r_{1},r_{2},r_{3} \in \Re.  \label{apr*}
\end{equation}
There remains to check the $H$-unitary invariance of the numbers $r_{1},\;r_{2},\;r_{3}$.
To prove this, let us suppose that some $H$-unitary matrix $T$ reduces (\ref{apr*})
to the form
$$\tilde{N}-\lambda I=\left( \begin{array}{ccccc}
				0 & 1 & 0 & 0 & i\widetilde{r_{3}} \\
				0 & 0 & 1 & i\widetilde{r_{1}} & -2\widetilde{r_{1}}^{2}+i\widetilde{r_{2}} \\
				0 & 0 & 0 & 1 & 2i\widetilde{r_{1}} \\
				0 & 0 & 0 & 0 & 1 \\
				0 & 0 & 0 & 0 & 0	\end{array}
	\right), \; \widetilde{r_{1}}, \widetilde{r_{2}}, \widetilde{r_{3}} \in \Re. $$
 From condition (\ref{11}) $NT=T\tilde{N}$ it follows that $T$ is uppertriangular
 with diagonal terms which are equal to each other.
 According to Theorem~1 from \cite{1}, $r_{1}$ is an $H_{1}$-unitary invariant
 for $N_{1}$. We already know that in this case $r_{1}$ must be an $H$-unitary
 invariant (see the previos case $n=4$), i.e., $\widetilde{r_{1}}=r_{1}$.
 For $T$ to be $H$-unitary, i.e., to satisfy (\ref{12}),
 $| t_{11}|$ must be equal to $1$. Therefore, as in case $n=4$, one can
 assume that $t_{11}=1$. Thus, $T$ has the form
 \begin{equation}
 T=\left( \begin{array}{ccccc}
				1 & t_{12} & t_{13} & t_{14} & t_{15} \\
				0 & 1 & t_{23} & t_{24} & t_{25} \\
				0 & 0 & 1 & t_{34} & t_{35} \\
				0 & 0 & 0 & 1 & t_{45} \\
				0 & 0 & 0 & 0 & 1	\end{array}
	\right). \label{apr**}
 \end{equation}
Condition (\ref{11}) amounts to system (\ref{e_sys1}) - (\ref{e_sys6}),
(\ref{12}) to system (\ref{h_sys1}) - (\ref{h_sys6}):
\begin{eqnarray}
& t_{23}=t_{12} & \label{e_sys1}\\
& t_{24}=ir_{1}t_{12}+t_{13} & \label{e_sys2} \\
& t_{25}+ir_{3}=i\widetilde{r_{3}}+(-2r_{1}^{2}+i\widetilde{r_{2}})t_{12}+2ir_{1}t_{13}+t_{14} & \label{e_sys3} \\
& t_{34}=t_{23} & \label{e_sys4} \\
& t_{35}+ir_{1}t_{45}+ir_{2}=i\widetilde{r_{2}}+2ir_{1}t_{23}+t_{24} & \label{e_sys5} \\
& t_{45}=t_{34}, & \label{e_sys6} \\
& \overline{t_{45}}+t_{12}=0 & \label{h_sys1}\\
& \overline{t_{35}}+\overline{t_{34}}t_{12}+t_{13}=0 & \label{h_sys2} \\
& \overline{t_{25}}+\overline{t_{24}}t_{12}+\overline{t_{23}}t_{13}+t_{14}=0 & \label{h_sys3} \\
& 2{\cal R}et_{15}+2{\cal R}e\{t_{12}\overline{t_{14}}\}+| t_{13}|^{2}=0 & \label{h_sys4} \\
& \overline{t_{34}}+t_{23}=0 & \label{h_sys5} \\
& 2{\cal R}et_{24}+| t_{23} |^{2}=0. & \label{h_sys6}
\end{eqnarray}
Express $t_{35}$ in terms of $t_{23}$, $t_{24}$, $t_{45}$ from
(\ref{e_sys5}) and substitute this expression in (\ref{h_sys2}),
taking into account that $t_{12}=t_{23}=t_{34}=t_{45}$ and expressing
$t_{24}$ in terms of $t_{12}$ and $t_{13}$ from condition (\ref{e_sys2}).
We obtain $ir_{2}-i\widetilde{{r}_{2}}=2ir_{1}\overline{t_{12}}+2{\cal R}et_{13}+| t_{12} | ^{2}$.
Since ${\cal R}et_{12}=0$ (equation (\ref{h_sys1})), we have $2ir_{1}\overline{t_{12}} \in \Re $,
hence, the right hand side of the condition obtained is real and the left one
is imaginary. Therefore, $\widetilde{{r}_{2}}=r_{2}$.

Since $t_{13}=t_{24}-ir_{1}t_{12}$ (condition (\ref{e_sys2})), $t_{25}$
can be expressed in terms of $t_{12}$, $t_{24}$ and $t_{14}$ in the following
way (see condition (\ref{e_sys3})):
$t_{25}=i(\widetilde{r_{3}}-r_{3})+ir_{2}t_{12}+2ir_{1}t_{24}+t_{14}$.
By substituting this expression in (\ref{h_sys3}), we get
$ir_{3}-i\widetilde{r_{3}}=ir_{2}\overline{t_{12}}+ir_{1}(2\overline{t_{24}}+|t_{12}|^{2})+2{\cal R}e\{t_{12}\overline{t_{24}}\}+2{\cal R}et_{14}$.
Because of condition (\ref{h_sys6})
$ir_{1}(2\overline{t_{24}}+|t_{12}|^{2})$ is real as well as the rest terms
of the right hand side, hence, $\widetilde{r_{3}}=r_{3}$. We have proved
the $H$-unitary invariance of $r_{1},\;r_{2},\;r_{3}$.

(b) ${z'=z,\;| z |=1,\;0< arg\;z< \pi, \;x=r_{1} \in \Re.\: }$ Applying
conditions (\ref{sys52}), (\ref{sys53}) of the $H$-normality of $N$, we get
$$N-\lambda I=\left( \begin{array}{ccccc}
				0 & 1 & 0 & 0 & ir_{3} \\
				0 & 0 & z & r_{1} & -2z^{2}r_{1}^{2}Im^{2}z+ir_{2}z^{2} \\
				0 & 0 & 0 & z & -2ir_{1}z^{2}Im\:z \\
				0 & 0 & 0 & 0 & z^{2} \\
				0 & 0 & 0 & 0 & 0	\end{array}
	\right)\;r_{1},r_{2},r_{3} \in \Re \;\; (z \neq i) $$
or
$$ N-\lambda I=\left( \begin{array}{ccccc}
				0 & 1 & 0 & 0 & r_{3} \\
				0 & 0 & i & r_{1} & 2r_{1}^{2}+ir_{2} \\
				0 & 0 & 0 & i & 2ir_{1} \\
				0 & 0 & 0 & 0 & -1 \\
				0 & 0 & 0 & 0 & 0	\end{array}
	\right)\;r_{1},r_{2},r_{3} \in \Re \; \;(z=i). $$
We shall join these cases, assuming that
$$
N-\lambda I=\left( \begin{array}{ccccc}
		   0 & 1 & 0 & 0 & ix \\
		   0 & 0 & z & r_{1} & -2z^{2}r_{1}^{2}Im^{2}z+ir_{2}z^{2} \\
		   0 & 0 & 0 & z & -2ir_{1}z^{2}Im\:z \\
		   0 & 0 & 0 & 0 & z^{2} \\
		   0 & 0 & 0 & 0 & 0	\end{array}
	\right), $$
where
$$x=\left[ \begin{array}{cc}
	r_{3}   \in \Re, & z \neq i \\
	-ir_{3} \in \Im \; (r_{3} \in \Re), & z=i. \end{array}
	\right. $$
Let us prove the $H$-unitary invariance of the numbers $z$, $r_{1}$, $r_{2}$, $r_{3}$ (or $x$).
Suppose some matrix $T$ realizes the $H$-unitary
transformation of $N$ to the form $\tilde{N}$, where
$$\tilde{N}-\lambda I=\left( \begin{array}{ccccc}
0 & 1 & 0 & 0 & i\tilde{x} \\
0 & 0 &\tilde{z} & \widetilde{r_{1}} & -2\tilde{z}^{2}\widetilde{r_{1}}^{2}Im^{2}\tilde{z}+i\widetilde{r_{2}}\tilde{z}^{2} \\
0 & 0 & 0 & \tilde{z} & -2i\widetilde{r_{1}}\tilde{z}^{2}{\cal I}m\tilde{z} \\
0 & 0 & 0 & 0 & \tilde{z}^{2} \\
0 & 0 & 0 & 0 & 0	\end{array}
\right).$$
By Theorem~1 of~\cite{1}, $z$ and $r_{1}$ are $H_{1}$-unitary invariants,
hence, they are $H$-unitary invariants, i.e., $\tilde{z}=z$, $\widetilde{r_{1}}=r_{1}$.
Further, from (\ref{11}) it follows that $T$ is uppertriangular with diagonal
terms which are equal to each other. Applying (\ref{12}), we get that $T$
has form (\ref{apr**}). Now condition (\ref{12}) is equivalent to system
(\ref{h_sys1}) - (\ref{h_sys6}),
condition (\ref{11}) to system (\ref{e_sys11}) - (\ref{e_sys66}):
\begin{eqnarray}
&t_{23} = zt_{12}& \label{e_sys11}\\
&t_{24} = r_{1}t_{12}+zt_{13}& \label{e_sys22} \\
&t_{25}+ix = i\tilde{x}+(-2z^{2}r_{1}^{2}Im^{2}z+i\widetilde{r_{2}}z^{2})t_{12}-2ir_{1}z^{2}Imz\:t_{13}+z^{2}t_{14}& \label{e_sys33} \\
&t_{34} = t_{23}& \label{e_sys44} \\
&zt_{35}+r_{1}t_{45}+ir_{2}z^{2} = i\widetilde{r_{2}}z^{2}-2ir_{1}z^{2}Imz\:t_{23}+z^{2}t_{24}& \label{e_sys55} \\
&zt_{45} = z^{2}t_{34}.& \label{e_sys66}
\end{eqnarray}
Express $t_{35}$ in terms of $t_{23}$, $t_{24}$, $t_{45}$ and, taking into account
the equalities
$t_{12}=\overline{z}t_{23}$ (\ref{e_sys11}),
$t_{13}=\overline{z}(t_{24}-r_{1}t_{12})$ (\ref{e_sys22}),
$t_{34}=t_{23}$ (\ref{e_sys44}), $t_{45}=zt_{23}$ (\ref{e_sys66}),
substitute the obtained expression in (\ref{h_sys2}). After multipling both sides by
$\overline{z}$, we have:
$(ir_{2}-i\widetilde{r_{2}})=-2ir_{1}Imz\:t_{23}+t_{24}+\overline{t_{24}}+|t_{23}|^{2}-r_{1}(\overline{z}t_{23}+z\overline{t_{23}})$.
Since ${\cal R}et_{23}=0$ (\ref{h_sys5}), the right hand side of this equality
is real. Consequently, $\widetilde{r_{2}}=r_{2}$.

Now let us express $t_{25}$ in terms of $t_{23}$, $t_{24}$, $t_{14}$ from
(\ref{e_sys33}) :
$t_{25}=i(\tilde{x}-x)-2r_{1}^{2}zIm^{2}z\:t_{23}+ir_{2}zt_{23}-2ir_{1}zImz\:t_{24}+2ir_{1}^{2}Imz\:t_{23}+z^{2}t_{14}$.
Rewrite condition (\ref{h_sys3}) in the form
$t_{25}+t_{24}\overline{t_{12}}+t_{23}\overline{t_{13}}+\overline{t_{14}}=0$,
multiply its both sides by $\overline{z}$ and  substitute the expression for
$t_{25}$ in it. We obtain:
$i(x-\tilde{x})\overline{z}=-2r_{1}^{2}Im^{2}z\:t_{23}+ir_{2}t_{23}-2ir_{1}Imz\:t_{24}+2ir_{1}^{2}\overline{z}Imz\:t_{23}+zt_{14}+\overline{z}\overline{t_{14}}+t_{23}\overline{t_{24}}+t_{24}\overline{t_{23}}-zr_{1}|t_{23}|^{2}$.
Since $-2r_{1}^{2}Im^{2}z+2ir_{1}^{2}\overline{z}Imz=ir_{1}^{2}Imz\:Rez$ and
$-2ir_{1}Imz\:t_{24}-r_{1}z|t_{23}|^{2}=r_{1}(2Rez\:Ret_{24}+2Imz\:Imt_{24})$,
the right hand side is real. Therefore, $Im[i\overline{z}(x-\tilde{x})]=0$.
If $z \neq i$, then this condition means $(r_{3}-\widetilde{r_{3}})Rez=0$,
hence $\widetilde{{r}_{3}}=r_{3}$ because $Re \: z \neq 0$. If $z=i$, then
 $Im[i(\widetilde{{r}_{3}}-r_{3})]=0$, hence also we get
$\widetilde{{r}_{3}}=r_{3}$. This concludes the proof of the $H$-unitary
invariance of $z,\;r_{1}\;r_{2},\;r_{3}$.

Due to Proposition~2 all obtained forms are indecomposable. They are not
$H$-unitarily similar because their internal matrices $N_{1}$ are not
$H_{1}$-unitarily similar due to Theorem~1 of \cite{1}. Thus, we have proved
the following lemma:
\begin{lemma}
If an indecomposable $H$-normal operator $N$ ($N: \: C^{5} \rightarrow C^{5}$)
has the only eigenvalue $\lambda$, $dim \:S_{0}=1$, the internal operator
$N_{1}$ is indecomposable, then the pair
$\{ N,H \}$ is unitarily similar to one and only one of canonical pairs
\{(\ref{lemma3.1}),(\ref{lemma3.4})\},
\{(\ref{lemma3.2}),(\ref{lemma3.4})\},
\{(\ref{lemma3.3}),(\ref{lemma3.4})\}:
$$ N=\left( \begin{array}{ccccc}
	\lambda & 1 & 0 & 0 & ir_{3} \\
	0 & \lambda & 1 & ir_{1} & -2r_{1}^{2}+ir_{2}  \\
	0 & 0 & \lambda & 1 & 2ir_{1}     \\
	0 & 0 & 0 & \lambda  & 1 \\
	0 & 0 & 0 & 0 & \lambda
	\end{array}
    \right), \;r_{1},r_{2},r_{3} \in \Re, $$
$$ N=\left( \begin{array}{ccccc}
	\lambda & 1 & 0 & 0 & ir_{3} \\
	0 & \lambda & z & r_{1} & -2z^{2}r_{1}^{2}Im^{2}z+ir_{2}z^{2}  \\
	0 & 0 & \lambda & z & -2ir_{1}z^{2}Imz   \\
	0 & 0 & 0 & \lambda  & z^{2} \\
	0 & 0 & 0 & 0 & \lambda
	\end{array}
        \right), \; \begin{array}{c}
                          |z|=1, \; z \neq i, \\
                           0< arg\:z < \pi, \\
                           r_{1},r_{2},r_{3} \in \Re, \end{array} $$
$$ N=\left( \begin{array}{ccccc}
	\lambda & 1 & 0 & 0 & r_{3} \\
	0 & \lambda & i & r_{1} & 2r_{1}^{2}+ir_{2}  \\
	0 & 0 & \lambda & i & 2ir_{1}     \\
	0 & 0 & 0 & \lambda  & -1 \\
	0 & 0 & 0 & 0 & \lambda
	\end{array}
\right), \; r_{1},r_{2},r_{3} \in \Re,  $$
$$ H=D_{5}, $$
where  $z, \; r_{1}, \; r_{2}, \; r_{3}$ are $H$-unitary invariants.
\end{lemma}

\subsubsection{$n=6$}
 In this case, according to Theorem~1 from~\cite{1}, the matrices
$N_{1}$ and $H_{1}$ can be written in the form
$$N_{1}=\left( \begin{array}{cccc}
				\lambda & \cos \alpha & \sin \alpha & 0 \\
				0 & \lambda & 0 & 1 \\
				0 & 0 & \lambda & 0 \\
				0 & 0 & 0 & \lambda	\end{array}
	\right), \; 0< \alpha \leq \pi /2,   \;\;
H_{1}=\left( \begin{array}{ccc}
				0 & 0 & 1 \\
				0 & I_{2} & 0 \\
				1 & 0 & 0 	\end{array}
	\right) \]
 so that
 $$ N- \lambda I=\left( \begin{array}{cccccc}
            0 & a & b & c & d & e   \\
				0 & 0 & \cos \alpha & \sin \alpha & 0 & f \\
				0 & 0 & 0 & 0 & 1 & g \\
				0 & 0 & 0 & 0 & 0 & h \\
                                0 & 0 & 0 & 0 & 0 & p \\
				0 & 0 & 0 & 0 & 0 & 0 \end{array}
	\right), \;
H=\left( \begin{array}{cccccc}
        0 & 0 & 0 & 0 & 0 & 1   \\
	0 & 0 & 0 & 0 & 1 & 0 \\
	0 & 0 & 1 & 0 & 0 & 0 \\
	0 & 0 & 0 & 1 & 0 & 0 \\
        0 & 1 & 0 & 0 & 0 & 0 \\
        1 & 0 & 0 & 0 & 0 & 0 \end{array}
	\right). $$
 The condition of the $H$-normality of $N$ is equivalent to the following system:
\begin{eqnarray}
  a & = & \overline{p} \cos \alpha \label {hi1} \\
  0 & = & \overline{p} \sin \alpha \label{hi2} \\
  b \cos \alpha + c \sin \alpha & = & \overline{g} \nonumber \\
  2 {\cal R}e \{a \overline{d} \}+|b|^{2}+|c|^{2} & = & 2 {\cal R}e\{f \overline{p} \}+|g|^{2}+|h|^{2}. \nonumber
\end{eqnarray}
From (\ref{hi2}) and the condition $0< \alpha \leq \pi /2 $ it follows that $p=0$.
Then from (\ref{hi1}) we obtain that also $a=0$. Hence, the  vector $v_{2} \in S$
belongs to $S_{0}$, which is impossible. This contradiction proves that for
indecomposable operator $N: \: C^{6} \rightarrow C^{6}$ $dim S_{0} \neq 1$.

Recall that if $n > 6$, then the operator $N_{1}$ is always decomposable
(Theorem~1 of~\cite{1}). Thus, we have obtained the classification for all
indecomposable operators $N$ having also indecomposable internal operator
$N_{1}$.

\subsection{$dim\:S_{0}=1$ and $N_{1}$ is Decomposable}
If the operator $N_{1}$ is decomposable, then it can be represented as an
orthogonal sum of indecomposable operators $N_{1}^{(1)}$, $\ldots$,
$N_{1}^{(p)}$: $N_{1}=N_{1}^{(1)} \oplus \ldots \oplus N_{1}^{(p)}$,
$H_{1}=H_{1}^{(1)} \oplus \ldots \oplus H_{1}^{(p)}$. Without loss of
generality it can be assumed that $H_{1}^{(1)}$ has one negative eigenvalue.
Denote $H_{1}^{(1)}$ by $H_{2}$, $N_{1}^{(1)}$ by $N_{2}$,
$H_{1}^{(2)} \oplus \ldots \oplus H_{1}^{(p)}$ by $H_{3}$,
$N_{1}^{(2)} \oplus \ldots \oplus N_{1}^{(p)}$ by $N_{3}$.
Since $H_{3}$ has only positive eigenvalues, one can assume that $H_{3}=I$.
$N_{3}$ is a usual normal operator having the only eigenvalue $\lambda$,
hence, $N_{3}=\lambda I$.

Show that the size of $N_{3}$ is equal to $1 \times 1$.
Indeed, let $dim V_{2}=k$, $dim V_{3}=l>1$ ($V_{2}$ and $V_{3}$ are the
subspaces of $S$ corresponding to $N_{2}$ and $N_{3}$, respectively),
$V_{2}=span \{ w_{1}^{(2)}, w_{2}^{(2)}, \ldots, w_{k}^{(2)}\}$,
$V_{3}=span \{w_{1}^{(3)}$,$w_{2}^{(3)}$,\linebreak $\ldots$,$w_{l}^{(3)}\}$.
Then, by the above,
$$
N=\left( \begin{array}{cccc}
	\lambda & M_{1} & M_{2} & * \\
	0 & N_{2} & 0 & * \\
	0 & 0 & \lambda I & * \\
	0 & 0 & 0 & \lambda	\end{array}
	\right), \;\;\;
N^{[*]}=\left( \begin{array}{cccc}
	\overline{\lambda} & M_{3} & M_{4} & * \\
	0 & N^{[*]}_{2} & 0 & * \\
	0 & 0 & \overline{\lambda} I & * \\
	0 & 0 & 0 & \overline{\lambda} \end{array}
	\right),
$$
where $M_{1}=(a_{1}, a_{2}, \ldots ,a_{k})$, $M_{2}=(b_{1}, b_{2}, \ldots ,b_{l})$,
$M_{3}=(c_{1}, c_{2}, \ldots ,c_{k})$, $M_{4}=(d_{1}, d_{2}, \ldots ,d_{l})$.
Because of the $H_{2}$-normality of $N_{2}$ $dim S_{0}^{(2)} \geq 1$
($S_{0}^{(2)}=\{ x \in V_{2}: \: (N_{2}- \lambda I)x=(N_{2}^{[*]}- \overline{\lambda} I)x=0 \}$),
hence, without loss of generality it can be  assumed that $w_{1}^{(2)} \in S_{0}^{(2)}$.
Since $l>1$, $\exists \{ \alpha _{i} \}_{1}^{n+1}$
($\sum_{1}^{n+1} |\alpha_{i}| \neq 0 $):
\begin{eqnarray}
 \sum_{1}^{n} \alpha_{i} b_{i} + \alpha_{n+1} a_{1} & = & 0 \\
 \sum_{1}^{n} \alpha_{i} d_{i} + \alpha_{n+1} c_{1} & = & 0.
\end{eqnarray}
Therefore, $\exists v=\sum_{1}^{n} \alpha _{i} w_{i}^{(3)}+ \alpha_{n+1} w_{1}^{(2)} \neq 0$:
$(N-\lambda I)v=(N^{[*]}- \overline{\lambda} I)v=0$, i.e., some nonzero vector
from $S$ belongs to $S_{0}$. This is impossible so that $dim \: V_{3}=1$.

As $N_{2}$ is indecomposable and rank of $V_{2}$ is less than or equal to $1$,
$dim \: V_{2} \leq 4$ in accordance with Theorem~1.
Thus, $1 \leq dim V_{2} \leq 4$, $dim V_{3}=1$ so that
$4 \leq n \leq 7$. Consider the cases $n=4,\;5,\;6,\;7$ one after another.

\subsubsection{$n=4$}
 Then $dim V_{2}=1$, $dim V_{3}=1$,
$$N- \lambda I=\left( \begin{array}{cccc}
	0 & a & b & c \\
	0 & 0 & 0 & d \\
	0 & 0 & 0 & e \\
	0 & 0 & 0 & 0	\end{array}
	\right), \;\;
H=\left( \begin{array}{cccc}
	0 & 0 & 0 & 1 \\
	0 &-1 & 0 & 0 \\
	0 & 0 & 1 & 0 \\
	1 & 0 & 0 & 0	\end{array}
	\right).
$$
Since $H_{1}=-1 \oplus 1$ is congruent to $D_{2}$, we will assume that
$H_{1}=D_{2}$ so that $H=D_{4}$. Having fixed $H=D_{4}$, we will apply, as is
customary, only $H$-unitary transformations.

The condition of the $H$-normality of $N$ is now equivalent to the following:
\begin{equation}
{\cal R}e\{ a \overline{b}\}={\cal R}e \{d \overline{e}\}. \label{apr***}
\end{equation}
Since the assumption $a=b=0$ contradicts the condition $S \cap S_{0}=\{0\}$
(because then either $v_{2}$ or $v_{3}$ belongs to $S_{0}$),
one can assume that $a \neq 0$ and, therefore, $a=1$ (see the paragraph
after (\ref{sys2})). Keeping in mind that $a=1$, reduce $N -\lambda I$ to the
form
 $$N- \lambda I=\left( \begin{array}{cccc}
	0 & 1 & b'=sgn {\cal R}e b & c' \\
	0 & 0 & 0 & d' \\
	0 & 0 & 0 & e' \\
	0 & 0 & 0 & 0	\end{array}
	\right), $$
having applied either the transformation
$$ T=\left( \begin{array}{cccc}
\sqrt{|{\cal R}e b|} & 0 & 0 & 0 \\
0 & \sqrt{|{\cal R}e b|} & -i{\cal I}m b/\sqrt{|{\cal R}e b|} & 0 \\
0 & 0 & 1/\sqrt{|{\cal R}e b|} & 0 \\
0 & 0 & 0 & 1/\sqrt{|{\cal R}e b|} \end{array}
\right)\;({\cal R}e b \neq 0) $$
or
$$ T=\left( \begin{array}{cccc}
1 & 0 & 0 & 0 \\
0 & 1 & -b & 0 \\
0 & 0 & 1 & 0 \\
0 & 0 & 0 & 1 \end{array}
\right) \; ({\cal R}e\:b = 0).$$

Now consider the three cases (${\cal R}e\:b'=0,1$ or $-1$) separately.

 (a) ${b'=0.\:}$ Since ${\cal R}e \{ d' \overline{e'} \}=0$ (condition
(\ref{apr***}) of the $H$-normality of $N$) and $d' \neq 0$ (otherwise
$v_{3} \in S_{0}$), the representation $d'=\varrho_{1}z$, $e'=i \varrho_{2}z$
($|z|=1$, $\varrho_{1},\varrho_{2} \in \Re$, $\varrho_{1}>0$) is valid.
Therefore, taking
$$T=\left( \begin{array}{cccc}
	\sqrt{\varrho_{1}} & 0 & 0 & 0 \\
	0 & \sqrt{\varrho_{1}} & 0 & 0 \\
	0 & i\varrho_{2}/\sqrt{\varrho_{1}} & 1/\sqrt{\varrho_{1}} & 0 \\
	0 & 0 & 0 & 1/\sqrt{\varrho_{1}} \end{array}
	\right),
$$
we reduce $N-\lambda I$ to the form
 $$N- \lambda I=\left( \begin{array}{cccc}
	0 & 1 & 0 & c'' \\
	0 & 0 & 0 & z \\
	0 & 0 & 0 & 0 \\
	0 & 0 & 0 & 0 \end{array}
	\right).
$$ 
One can assume that $c''=0$. To achieve this it is sufficient to apply the
transformation
$$T=\left( \begin{array}{cccc}
	1 & 0 & \overline{c''} & 0 \\
	0 & 1 & 0 & -c'' \\
	0 & 0 & 1 & 0 \\
	0 & 0 & 0 & 1 \end{array}
	\right). $$
There remains to prove that $z$ is an $H$-unitary invariant.
Indeed, any matrix $T$ satisfying condition (\ref{11})
$(N-\lambda I)T=T(\tilde{N}-\lambda I)$ for the matrices
$$ N- \lambda I=\left( \begin{array}{cccc}
	0 & 1 & 0 & 0 \\
	0 & 0 & 0 & z \\
	0 & 0 & 0 & 0 \\
	0 & 0 & 0 & 0	\end{array}
	\right), \; \;
\tilde{N}- \lambda I=\left( \begin{array}{cccc}
	0 & 1 & 0 & 0 \\
	0 & 0 & 0 & \tilde{z} \\
	0 & 0 & 0 & 0 \\
	0 & 0 & 0 & 0	\end{array}
	\right), \; \; |z|=|\tilde{z}|=1 $$
and condition (\ref{12}) $TT^{[*]}=I$  has the form
 $$T=t_{11}\left( \begin{array}{cccc}
	1 & * & * & * \\
	0 & 1 & 0 & * \\
	0 & 0 & 1 & * \\
	0 & 0 & 0 & 1	\end{array}
	\right), \;\; |t_{11}|=1.$$
This follows the desired equality $z=\tilde{z}$.

(b) ${b'=1.\:}$ As ${\cal R}e \{d'\overline{e'}\}=1$
(condition (\ref{apr***})), $d'=\varrho z$, $e'=(1/ \varrho+ir) z$
($|z|=1$, $\varrho, r \in \Re$, $\varrho > 0$). Consider the transformation
\begin{equation}
T=I_{1} \oplus
    \left( \begin{array}{cc}
	-it/(1-it) & 1/(1-it) \\
	1/(1-it) & -it/(1-it) \end{array}
	\right) \oplus I_{1},\;\;\; t \in \Re, \label{apr****}
\end{equation}
where $t$ is a root of the equation  $1+t^{2}=1/ \varrho^{2}+(t \varrho+r)^{2}$.
Its discriminant ${\cal D}/4=1/{\varrho^{2}}+\varrho^{2}+r^{2}-2$ is
nonnegative so that $t$ is in fact real. Subjecting to (\ref{apr****}),
the matrix $N - \lambda I$ becomes the following:
$$N - \lambda I=\left( \begin{array}{cccc}
	0 & 1 & 1 & c'' \\
	0 & 0 & 0 & z' \\
	0 & 0 & 0 & (1+ir')z' \\
	0 & 0 & 0 & 0	\end{array}
	\right), \; |z'|=1, \; r' \in \Re.$$
Note that if $r'=0$, then there exists a nonzero vector 
$v=\alpha v_{2}+ \beta v_{3} \in S_{0}$, which is impossible. 
Applying (\ref{apr****}) with $t=-\frac{1}{2}r'$, we can replace $r'$
by $-r'$. Thus, we can assume $r'>0$. Finally, to get $c''=0$ it is
sufficient to take
$$T=\left( \begin{array}{cccc}
	1 & t_{12} & t_{13} & -{\cal R}e \{ t_{12} \overline{t_{13}}\} \\
	0 & 1 & 0 & -\overline{t_{13}} \\
	0 & 0 & 1 & -\overline{t_{12}} \\
	0 & 0 & 0 & 1	\end{array}
	\right), $$
where $t_{12}=e^{{-i\varphi}/2}(rc_{1}''-2c_{2}'')/(2r)$,
$t_{13}=e^{{-i\varphi}/2}c_{2}''/r$ (we mean that $z'=e^{i\varphi}$,
$c_{1}''={\cal R}e \{c''e^{-i\varphi/2}\}$,
$c_{2}''={\cal I}m \{c''e^{-i\varphi/2}\}$).

Thus, we have reduced the matrix $N- \lambda I$ to the form
 $$  N- \lambda I=\left( \begin{array}{cccc}
	0 & 1 & 1 & 0 \\
	0 & 0 & 0 & z  \\
	0 & 0 & 0 & (1+ir)z \\
	0 & 0 & 0 & 0	\end{array}
	\right), \; | z |=1, \; r \in \Re>0. $$
Now there remains to show that the numbers $z$ and $r$ are $H$-unitary
invariants.

First note that for a block triangular matrix
\begin{equation}
T=\left( \begin{array}{ccc}
	T_{1} & T_{2} & T_{3} \\
	0 & T_{4} & T_{5} \\
	0 & 0 & T_{6} \end{array}
	\right) \label{t-t}
\end{equation}
to reduce $N - \lambda I$ to the form $\tilde{N} - \lambda I$,
where
$$N- \lambda I=\left( \begin{array}{ccc}
	0 & N_{1} & N_{2} \\
	0 & N_{3} & N_{4} \\
	0 & 0 & 0  \end{array}
	\right), \; \;
\tilde{N}- \lambda I=\left( \begin{array}{ccc}
	0 & \widetilde{N_{1}} & \widetilde{N_{2}} \\
	0 & \widetilde{N_{3}} & \widetilde{N_{4}} \\
	0 & 0 & 0  \end{array}
	\right), $$
it is necessary and sufficient to have
\begin{eqnarray}
  N_{1} T_{4} & = & T_{1} \widetilde{N_{1}}+T_{2} \widetilde{N_{3}}  \label{eqvv1} \\
  N_{1} T_{5} + N_{2} T_{6} & = & T_{1} \widetilde{N_{2}} + T_{2} \widetilde{N_{4}}  \label{eqvv2} \\
  N_{3} T_{4} & = & T_{4} \widetilde{N_{3}} \label{eqvv3} \\
  N_{3} T_{5} + N_{4} T_{6} & = & T_{4} \widetilde{N_{4}}. \label{eqvv4}
\end{eqnarray}
If
$$ H=\left( \begin{array}{ccc}
	0 & 0 & I \\
	0 & H_{1} & 0 \\
	I & 0 & 0 \end{array}
	\right),  $$
then  for (\ref{t-t}) to be $H$-unitary it is necessary and sufficient to have
\begin{eqnarray}
  T_{1} T_{6}^{*} & = & I \label{h_uu1} \\
  T_{4} H_{1} T_{2}^{*} + T_{5} T_{1}^{*} & = & 0 \label{h_uu2} \\
  T_{1} T_{3}^{*} + T_{2} H_{1} T_{2}^{*} + T_{3} T_{1}^{*} & = & 0 \label{h_uu3} \\
  T_{4} H_{1} T_{4}^{*} H_{1} & = & I. \label{h_uu4}
\end{eqnarray}

Since any $H$-unitary transformation $T$ such that
 $$  \left( \begin{array}{cccc}
	0 & 1 & 1 & 0 \\
	0 & 0 & 0 & z  \\
	0 & 0 & 0 & (1+ir)z \\
	0 & 0 & 0 & 0	\end{array}
	\right) T=T \left( \begin{array}{cccc}
	                  0 & 1 & 1 & 0 \\
	                  0 & 0 & 0 & \tilde{z}  \\
	                  0 & 0 & 0 & (1+i\tilde{r})\tilde{z} \\
	                  0 & 0 & 0 & 0	\end{array}
       \right),
$$
$| z |=| \tilde{z} |=1,\; r, \; \tilde{r} \in \Re > 0$,
has to be block triangular (by Corollary of Proposition~1),
systems (\ref{eqvv1}) - (\ref{eqvv4}), (\ref{h_uu1}) - (\ref{h_uu4})
are applicable. Combining (\ref{eqvv1}) and (\ref{h_uu4}), we get
$|t_{11}|=1$, hence (condition (\ref{h_uu1})) $t_{44}=t_{11}$. Now
from (\ref{eqvv1}) and (\ref{eqvv4}) it follows that
$(2+ir)z=(2+i\tilde{r})\tilde{z}$, hence $\tilde{z}=z$, $\tilde{r}=r$, Q.E.D.

(c)  ${b'=-1. \:}$ The matrix $N-\lambda I$ can be carried into the form
 $$ N-\lambda I= \left( \begin{array}{cccc}
	0 & 1 & -1 & 0 \\
	0 & 0 & 0 & z  \\
	0 & 0 & 0 & -(1+ir)z \\
	0 & 0 & 0 & 0	\end{array}
     \right),  \; | z |=1, \; r \in \Re>0, $$
where $z$ and $r$ are $H$-unitary invariants. The proof is analogous to
the case (b) above.

Thus, we have obtained the canonical form for each case considered.
By using conditions (\ref{eqvv1}) - (\ref{h_uu4}) one can easily check
that these forms are not $H$-unitarily similar to each other. They are
indecomposable due to Proposition~2. Thus, we have proved the following lemma:
\begin{lemma}
If an indecomposable $H$-normal operator $N$ ($N: \: C^{4} \rightarrow C^{4}$)
has the only eigenvalue $\lambda$, $dim\:S_{0}=1$,
the internal operator $N_{1}$ is decomposable, then the pair $\{ N,H \}$ is
unitarily similar to one and only one of canonical pairs
\{(\ref{lemma4.1}),(\ref{lemma4.4})\},
\{(\ref{lemma4.2}),(\ref{lemma4.4})\},
\{(\ref{lemma4.3}),(\ref{lemma4.4})\}:
$$ N=\left( \begin{array}{cccc}
	\lambda & 1 & 0 & 0  \\
	0 & \lambda & 0 & z  \\
	0 & 0 & \lambda & 0  \\
	0 & 0 & 0 & \lambda
	\end{array}
    \right), \; |z|=1, $$
$$ N=\left( \begin{array}{cccc}
	\lambda & 1 & 1 & 0 \\
	0 & \lambda & 0 & z \\
	0 & 0 & \lambda & (1+ir)z  \\
	0 & 0 & 0 & \lambda
	\end{array}
\right), \; |z|=1, \; r \in \Re > 0, $$
$$ N=\left( \begin{array}{ccccc}
	\lambda & 1 & -1 & 0 \\
	0 & \lambda & 0 & z \\
	0 & 0 & \lambda & -(1+ir)z     \\
	0 & 0 & 0 & \lambda
	\end{array}
\right), \; |z|=1, \; r \in \Re > 0, $$
$$ H=D_{4}, $$
where  $z$, $r$ are $H$-unitary invariants.
\end{lemma}
\subsubsection{$n=5$}
 Then $dim V_{2}=2$, $dim V_{3}=1$ and, according
to Theorem~1 from \cite{1}, after interchanging the $3$-rd and $4$-th
rows and colomns, we get:
$$N- \lambda I=\left( \begin{array}{ccccc}
	0 & a & b & c & d \\
	0 & 0 & 0 & z & e \\
	0 & 0 & 0 & 0 & f \\
	0 & 0 & 0 & 0 & g \\
	0 & 0 & 0 & 0 & 0 \end{array}
	\right), \; | z |=1, \;
H=\left( \begin{array}{ccccc}
	0 & 0 & 0 & 0 & 1 \\
	0 & 0 & 0 & 1 & 0 \\
	0 & 0 & 1 & 0 & 0 \\
	0 & 1 & 0 & 0 & 0 \\
	1 & 0 & 0 & 0 & 0 \end{array}
	\right). \]

The condition of the $H$-normality of $N$ is equivalent to the system
\begin{eqnarray}
  a \overline{z} & = & \overline{g} z \label{ooo1} \\
  2{\cal R}e \{a \overline{c}\}+| b |^{2} & = & 2{\cal R}e\{e \overline{g} \}+| f |^{2}. \label{ooo2}
\end{eqnarray}
It is readily seen that $a \neq 0$, consequently, it can be assumed that $a=1$
and $g=z^{2}$ (see the paragraph after (\ref{sys2})).
Further, take the ($H$-unitary) transformation
$$T=\left( \begin{array}{ccccc}
	1 & 0 & 0 & 0 & 0 \\
	0 & 1 & -b & -\frac{1}{2} | b |^{2} & 0 \\
	0 & 0 & 1 & \overline{b} & 0 \\
	0 & 0 & 0 & 1 & 0 \\
	0 & 0 & 0 & 0 & 1 \end{array}
	\right) $$
and reduce $N-\lambda I$ to the form
$$ N- \lambda I=\left( \begin{array}{ccccc}
	0 & 1 & 0 & c' & d' \\
	0 & 0 & 0 & z & e' \\
	0 & 0 & 0 & 0 & f' \\
	0 & 0 & 0 & 0 & z^{2} \\
	0 & 0 & 0 & 0 & 0 \end{array}
	\right). $$
Applying now the transformation
$$ T=I_{1} \oplus \left( \begin{array}{ccc}
	1 & 0 & i{\cal I}m\{e'\overline{z}^{2}\} \\
	0 & e^{i\:arg\:f'} & 0 \\
	0 & 0 & 1 \end{array}
	\right) \oplus I_{1}, $$
we get
$$ N-\lambda I=\left( \begin{array}{ccccc}
	0 & 1 & 0 & c''& d'' \\
        0 & 0 & 0 & z & r_{1}z^{2} \\
        0 & 0 & 0 & 0 & r_{2} \\
        0 & 0 & 0 & 0 & z^{2} \\
        0 & 0 & 0 & 0 & 0  \end{array}
	\right), \; r_{1},r_{2} \in \Re, \; r_{2} \geq 0. $$
We can assume that $r_{2}>0$ because otherwise $v_{3} \in S_{0}$, which
is impossible. From condition (\ref{ooo2}) of the $H$-normality of $N$ it
follows that
$c''=r_{1}+\frac{1}{2}r_{2}^{2}+ir_{3}$ ($r_{3} \in \Re$).
Keeping in mind these conditions, apply the transformation
$$T=\left( \begin{array}{ccccc}
	1 & t_{12} & t_{13} & 0 & -\frac{1}{2}|t_{13}|^{2} \\
	0 & 1 & 0 & 0 & 0 \\
	0 & 0 & 1 & 0 & -\overline{t_{13}} \\
	0 & 0 & 0 & 1 & -\overline{t_{12}} \\
	0 & 0 & 0 & 0 & 1 \end{array}
	\right), $$
where $t_{12}=r_{1}\overline{z}$,
$t_{13}=(d''-r_{1}z(r_{1}+\frac{1}{2}r_{2}^{2}+ir_{3}))/r_{2}$,
to the matrix $N - \lambda I$. Then $c'''=\frac{1}{2}r_{2}^{2}+ir_{3}$,
$d'''=0$, the rest terms of $N - \lambda I$ do not change.
Renaming $r_{2}$ and $r_{3}$, write out the final form of $N -\lambda I$:
$$N- \lambda I=\left( \begin{array}{ccccc}
	0 & 1 & 0 & \frac{1}{2}r_{1}^{2}+ir_{2} & 0 \\
	0 & 0 & 0 & z & 0 \\
	0 & 0 & 0 & 0 & r_{1} \\
	0 & 0 & 0 & 0 & z^{2} \\
	0 & 0 & 0 & 0 & 0 \end{array}
	\right), \; r_{1}, r_{2} \in \Re,\; r_{1}>0,\;|z|=1. $$

To prove the $H$-unitary invariance of $z$, $r_{1}$, $r_{2}$ assume that
$$ \tilde{N}- \lambda I=\left( \begin{array}{ccccc}
	0 & 1 & 0 & \frac{1}{2}\widetilde{r_{1}}^{2}+i\widetilde{r_{2}} & 0 \\
	0 & 0 & 0 & \tilde{z} & 0 \\
	0 & 0 & 0 & 0 & \tilde{r}_{1} \\
	0 & 0 & 0 & 0 & \tilde{z}^{2} \\
	0 & 0 & 0 & 0 & 0 \end{array}
	\right), \; \widetilde{r_{1}},\widetilde{r_{2}} \in \Re, \; \widetilde{r_{1}}>0, \; | \tilde{z} |=1,
$$
and there exists a matrix $T$ such that $NT=T\tilde{N}$ (condition (\ref{11}))
and $TT^{[*]}=I$ (condition (\ref{12})). Recall that $T$ has block form
(\ref{t-t}) so that conditions (\ref{eqvv1}) - (\ref{h_uu4}) hold.
From (\ref{eqvv3}) it follows that $t_{23}=0$ and $zt_{44}=\tilde{z}t_{22}$.
Since $t_{22}\overline{t_{44}}=1$ (\ref{h_uu4}), $z|t_{44}|^{2}=\tilde{z}$,
i.e., $\tilde{z}=z$, $|t_{44}|=1$. Therefore, one can assume that
$$T_{4}=\left( \begin{array}{ccc}
	1 & 0 & it \\
	0 & t_{33} & 0 \\
	0 & 0 & 1  \end{array}
	\right), \; |t_{33}|=1, \; t \in \Re $$
because it is allowed to divide $T$ by its term $t_{22}=t_{44}$ of modulus
$1$. Now from (\ref{eqvv4}) it follows that $t_{45}=itz$,
$\widetilde{r_{1}}t_{33}=r_{1}$. As $r_{1}$, $\widetilde{r_{1}}>0$,
$t_{33}=1$ and $\widetilde{r_{1}}=r_{1}$. Since
$t_{12}=-\overline{t_{45}}$ (condition (\ref{h_uu2})) and
$t_{24}+(\frac{1}{2}r_{1}^{2}+ir_{2})t_{44}=(\frac{1}{2}\widetilde{r_{1}}^{2}+i\widetilde{r_{2}})t_{11}+\tilde{z}t_{12}$
(condition (\ref{eqvv1})), $\widetilde{r_{2}}=r_{2}$. This completes the
proof of the $H$-unitary invariance of $z$, $r_{1}$, $r_{2}$.

Due to Proposition~2 the obtained form is
indecomposable. Thus, we have proved the following lemma:
\begin{lemma}
If an indecomposable $H$-normal operator $N$ ($N: \: C^{5} \rightarrow C^{5}$)
has the only eigenvalue $\lambda$, $dim \:S_{0}=1$,
the internal operator $N_{1}$ is decomposable, then the pair $\{N,H\}$
is unitarily similar to canonical pair 
\{(\ref{lemma5.1}),(\ref{lemma5.2})\}:
$$ N=\left( \begin{array}{ccccc}
	\lambda & 1 & 0 & \frac{1}{2}r_{1}^{2}+ir_{2} & 0 \\
	0 & \lambda & 0 & z & 0 \\
	0 & 0 & \lambda & 0 & r_{1}  \\
	0 & 0 & 0 & \lambda & z^{2}  \\
	0 & 0 & 0 & 0 & \lambda	\end{array}
\right), \; |z|=1, \; r_{1},r_{2} \in \Re, \; r_{1}>0, $$
$$ H=D_{5}, $$
where  $r_{1}$, $r_{2}$, $z$ are $H$-unitary invariants.
\end{lemma}

\subsubsection{$n=6$}
In this case $dim V_{2}=3$, $dim V_{3}=1$. The matrices
$N- \lambda I$ and $H$, according to Theorem~1 from \cite{1},
have the form:
\begin{equation}
N- \lambda I=\left( \begin{array}{cccccc}
        0 & a & b & c & d & e \\
	0 & 0 & z & r & 0 & f \\
	0 & 0 & 0 & z & 0 & g \\
	0 & 0 & 0 & 0 & 0 & h \\
	0 & 0 & 0 & 0 & 0 & p \\
	0 & 0 & 0 & 0 & 0 & 0 \end{array}
	\right), \; | z |=1, \; r \in \Re \label{case1}
\end{equation}
or
\begin{equation}
N- \lambda I=\left( \begin{array}{cccccc}
        0 & a & b & c & d & e  \\
	0 & 0 & 1 & ir & 0 & f \\
	0 & 0 & 0 & 1 & 0 & g \\
	0 & 0 & 0 & 0 & 0 & h \\
	0 & 0 & 0 & 0 & 0 & p \\
	0 & 0 & 0 & 0 & 0 & 0 \end{array}
	\right), \; r \in \Re, \label{case2}
\end{equation}
$$H=\left( \begin{array}{cccc}
        0 & 0 & 0 & 1     \\
	0 & D_{3} & 0 & 0 \\
	0 & 0 & 1 & 0     \\
	1 & 0 & 0 & 0     \end{array}
	\right).\]
For a while we will consider these two cases together, assuming that
$$ N- \lambda I=\left( \begin{array}{cccccc}
        0 & a & b & c & d & e \\
	0 & 0 & z & x & 0 & f \\
	0 & 0 & 0 & z & 0 & g \\
	0 & 0 & 0 & 0 & 0 & h \\
	0 & 0 & 0 & 0 & 0 & p \\
	0 & 0 & 0 & 0 & 0 & 0 \end{array}
	\right), \; | z |=1, \; x \in C. \]
Then the condition of the $H$-normality of $N$ is equivalent to the system
\begin{eqnarray}
  a \overline{z} & = & z \overline{h} \label{system1} \\
  a \overline{x}+b \overline{z} & = & x \overline{h}+ z \overline{g} \label{system2} \\
  2 {\cal R}e \{ a\overline{c} \} + |b|^{2}+|d|^{2} & = & 2 {\cal R}e \{ f \overline{h} \}+|g|^{2}+|p|^{2}. \label{system3}
\end{eqnarray}
As is customary, we can assume that $a=1$, $h=z^{2}$. Let us use the
($H$-unitary) transformation
$$T=\left( \begin{array}{cccccc}
	1 & 0 & 0 & 0 & 0 & 0 \\
	0 & 1 & 0 & -\frac{1}{2}|d|^{2} & -d & 0 \\
	0 & 0 & 1 & 0 & 0 & 0 \\
	0 & 0 & 0 & 1 & 0 & 0 \\
	0 & 0 & 0 & \overline{d} & 1 & 0 \\
	0 & 0 & 0 & 0 & 0 & 1 \end{array}
	\right). \]
It reduces $N-\lambda I$ to the form
$$
N- \lambda I=\left( \begin{array}{cccccc}
	0 & 1 & b' & c' & 0 & e' \\
	0 & 0 & z & x & 0 & f' \\
	0 & 0 & 0 & z & 0 & g' \\
	0 & 0 & 0 & 0 & 0 & z^{2} \\
	0 & 0 & 0 & 0 & 0 & p' \\
	0 & 0 & 0 & 0 & 0 & 0 \end{array}
	\right).
\]
Further, take the transformation
$$T=\left( \begin{array}{cccccc}
	1 & z \overline{g'} & \overline{z}c'-x\overline{g'} & 0 & 0 & -\frac{1}{2}|\overline{z}c'-x \overline{g'}|^{2} \\
	0 & 1 & 0 & 0 & 0 & 0 \\
	0 & 0 & 1 & 0 & 0 & -z \overline{c'}+\overline{x}g' \\
	0 & 0 & 0 & 1 & 0 & -\overline{z}g' \\
	0 & 0 & 0 & 0 & 1 & 0 \\
	0 & 0 & 0 & 0 & 0 & 1 \end{array}
	\right) \]
and carry the matrix $N-\lambda I$ into the form
$$ N- \lambda I=\left( \begin{array}{cccccc}
	0 & 1 & b'' & 0 & 0 & e'' \\
	0 & 0 & z & x & 0 & f'' \\
	0 & 0 & 0 & z & 0 & 0 \\
	0 & 0 & 0 & 0 & 0 & z^{2} \\
	0 & 0 & 0 & 0 & 0 & p'' \\
	0 & 0 & 0 & 0 & 0 & 0 \end{array}
	\right). \]
Now note that $p'' \neq 0$ because otherwise $v_{5} \in S_{0}$.
Since the rotation of the vector $v_{5}$ about any angle
does not change the matrix $H$, we can assume that $p''=r_{2} \in \Re >0$
(we put $\widetilde{v_{5}}=e^{i\:arg\:p''}v_{5}$).
The transformation
$$ T=\left( \begin{array}{cccccc}
	1 & 0 & 0 & 0 & e''/r_{2} & -\frac{1}{2}|e''/r_{2}|^{2} \\
	0 & 1 & 0 & 0 & 0 & 0 \\
	0 & 0 & 1 & 0 & 0 & 0 \\
	0 & 0 & 0 & 1 & 0 & 0 \\
	0 & 0 & 0 & 0 & 1 & -\overline{e}''/r_{2} \\
	0 & 0 & 0 & 0 & 0 & 1 \end{array}
	\right) \]
reduces the matrix $N-\lambda I$ to the final form:
$$ N- \lambda I=\left( \begin{array}{cccccc}
	0 & 1 & b''' & 0 & 0 & 0 \\
	0 & 0 & z & x & 0 & f''' \\
	0 & 0 & 0 & z & 0 & 0 \\
	0 & 0 & 0 & 0 & 0 & z^{2} \\
	0 & 0 & 0 & 0 & 0 & r_{2} \\
	0 & 0 & 0 & 0 & 0 & 0 \end{array}
	\right). \]

Now we will distinguish the cases (\ref{case1}) and (\ref{case2}).

(a) ${z=1, x \in \Im. \:}$ According to conditions (\ref{system2}) and
(\ref{system3}) of the $H$-norm\-al\-ity of $N$,
$$ N- \lambda I=\left( \begin{array}{cccccc}
	0 & 1 & 2ir_{1} & 0 & 0 & 0 \\
	0 & 0 & 1 & ir_{1} & 0 & 2r_{1}^{2}-r_{2}^{2}/2+ir_{3} \\
	0 & 0 & 0 & 1 & 0 & 0 \\
	0 & 0 & 0 & 0 & 0 & 1 \\
	0 & 0 & 0 & 0 & 0 & r_{2} \\
	0 & 0 & 0 & 0 & 0 & 0 \end{array}
	\right), \begin{array}{c} r_{1},\;r_{2},\;r_{3}\in \Re, \\
                                  r_{2}>0. \end{array} \]
Let us show that $r_{1}$, $r_{2}$, $r_{3}$ are $H$-unitary invariants.
Indeed, suppose some matrix $T$ satisfies  conditions (\ref{12}) $TT^{[*]}=I$
and (\ref{11}) $(N - \lambda I)T=T(\tilde{N} - \lambda I)$, where
$$ \tilde{N}- \lambda I=\left( \begin{array}{cccccc}
	0 & 1 & 2i\widetilde{r_{1}} & 0 & 0 & 0 \\
	0 & 0 & 1 & i\widetilde{r_{1}} & 0 & 2\widetilde{r_{1}}^{2}-\widetilde{r_{2}}^{2}/2+i\widetilde{r_{3}} \\
	0 & 0 & 0 & 1 & 0 & 0 \\
	0 & 0 & 0 & 0 & 0 & 1 \\
	0 & 0 & 0 & 0 & 0 & \widetilde{r_{2}} \\
	0 & 0 & 0 & 0 & 0 & 0 \end{array}
	\right), \begin{array}{c} \widetilde{r_{1}},\;\widetilde{r_{2}},\;\widetilde{r_{3}}\in \Re, \\
                                  \widetilde{r_{2}}>0. \end{array} \]
From (\ref{11}) it follows that
$$ T=\left( \begin{array}{cccccc}
	t_{11} & t_{12} & t_{13} & t_{14} & t_{15} & t_{16} \\
	0 & t_{11} & t_{23} & t_{24} & 0 & t_{26} \\
	0 & 0 & t_{11} & t_{34} & 0 & t_{36} \\
	0 & 0 & 0 & t_{11} & 0 & t_{46} \\
	0 & 0 & 0 & t_{54} & t_{55} & t_{56} \\
	0 & 0 & 0 & 0 & 0 & t_{11} \end{array}
	\right). \]
Using (\ref{h_uu4}), we get: $t_{54}=0$, $|t_{11}|=1$. As above (see the
argument before Lemma~5), we can assume that $t_{11}=1$.
Then $t_{34}=-\overline{t_{23}}$ (condition (\ref{h_uu4})) and
$i(\widetilde{r_{1}}-r_{1})=t_{34}-t_{23}$ (condition (\ref{eqvv3})),
hence, $\widetilde{r_{1}}=r_{1}$ and ${\cal R}e t_{23}=0$. Further, from
(\ref{eqvv4}) it follows that $r_{2}=\widetilde{r_{2}}t_{55}$,
from (\ref{h_uu4}) that $|t_{55}|=1$. As $r_{2},\widetilde{r_{2}}>0$,
$\widetilde{r_{2}}=r_{2}$ and $t_{55}=1$. Thus,
$$ T=\left( \begin{array}{cccc}
	1 & it & t_{24} & 0 \\
	0 & 1 & it & 0 \\
	0 & 0 & 1 & 0 \\
	0 & 0 & 0 & 1 \end{array}
	\right), \;\; t \in \Re,\; 2{\cal R}e t_{24}+t^{2}=0. \]
Substituting $T_{4}$ in (\ref{eqvv1}), we get $t_{12}=it$,
$t_{13}=t_{24}-r_{1}t$; replacing $T_{5}$ by $-T_{4}H_{1}T_{2}^{*}$
in (\ref{eqvv4}), we have
$i\widetilde{r_{3}}=ir_{3}-2 {\cal R}e t_{24}-t^{2}$, hence
$\widetilde{r_{3}}=r_{3}$. This completes the proof of the $H$-unitary
invariance of $r_{1}$, $r_{2}$, $r_{3}$.

(b) ${0<arg\: z<\pi, x \in \Re.\:}$ Applying the condition of the $H$-normality,
we get
$$
N- \lambda I=\left( \begin{array}{cccccc}
	0 & 1 & -2ir_{1}{\cal I}m z & 0 & 0 & 0 \\
	0 & 0 & z & r_{1} & 0 & (2r_{1}^{2}{\cal I}m^{2}z-r_{2}^{2}/2+ir_{3})z^{2} \\
	0 & 0 & 0 & z & 0 & 0 \\
	0 & 0 & 0 & 0 & 0 & z^{2} \\
	0 & 0 & 0 & 0 & 0 & r_{2} \\
	0 & 0 & 0 & 0 & 0 & 0 \end{array}
	\right),
\]
where $r_{1},\:r_{2},\:r_{3} \in \Re$, $r_{2}>0$. That the numbers
$z$, $r_{1}$, $r_{2}$, $r_{3}$ are $H$-unitary invariants can be checked
as in (a) above. That the forms obtained are not $H$-unitary similar can
also be checked by the reader by using formulas (\ref{eqvv1}) - (\ref{h_uu4}).

Because of Proposition~2 the forms obtained are indecomposable so that we
have proved the following lemma:
\begin{lemma}
If an indecomposable $H$-normal operator $N$ ($N: \: C^{6} \rightarrow C^{6}$)
has the only eigenvalue $\lambda$, $dim \:S_{0}=1$,
the internal operator $N_{1}$ is decomposable, then the pair $\{N,H\}$
is unitarily similar to one and only one of canonical pairs
\{(\ref{lemma6.1}),(\ref{lemma6.3})\},
\{(\ref{lemma6.2}),(\ref{lemma6.3})\}:
$$ N=\left( \begin{array}{cccccc}
	\lambda & 1   &   2ir_{1} & 0      & 0 & 0 \\
	0 &   \lambda &   1       & ir_{1} & 0 & 2r_{1}^{2}-r_{2}^{2}/2+ir_{3} \\
	0 &       0   & \lambda   & 1      & 0 & 0  \\
	0 &       0   & 0     & \lambda    & 0 & 1  \\
        0 &       0   & 0     &     0      & \lambda & r_{2}  \\
	0 &       0   & 0     &     0      & 0  & \lambda \end{array}
\right),\; r_{1},r_{2}\in \Re, \; r_{2}>0, \]
$$ N=\left( \begin{array}{cccccc}
	\lambda & 1   &  -2ir_{1}{\cal I}mz & 0      & 0 & 0 \\
	0 &   \lambda &   z       & r_{1} & 0 & (2r_{1}^{2}{\cal I}m^{2}z-r_{2}^{2}/2+ir_{3})z^{2} \\
	0 &       0   & \lambda   & z      & 0 & 0  \\
	0 &       0   & 0     & \lambda    & 0 & z^{2}  \\
        0 &       0   & 0     &     0      & \lambda & r_{2}  \\
	0 &       0   & 0     &     0      & 0  & \lambda \end{array}
\right), \]
$$ |z|=1, \; 0<arg\:z<\pi,\;r_{1},r_{2},r_{3}\in \Re, \; r_{2}>0, \]
$$ H=\left( \begin{array}{cccc}
		0 & 0 & 0 & 1     \\
		0 & D_{3} & 0 & 0 \\
		0 & 0 & 1 & 0     \\
                1 & 0 & 0 & 0
		\end{array}
	\right), \]
where  $z$, $r_{1}$, $r_{2}$, $r_{3}$ are $H$-unitary invariants.
\end{lemma}

\subsubsection{$n=7$}
We will show that this alternative is impossible. Indeed, if $dim V_{2}=4$,
$dim V_{3}=1$, then, in accordance with Theorem~1 of~\cite{1},
$$N - \lambda I=\left(
\begin{array}{ccccccc}
 0 & a & b & c & d & e & f \\
 0 & 0 & cos \alpha & sin \alpha & 0 & 0 & g \\
 0 & 0 & 0 & 0 & 1 & 0 & h  \\
 0 & 0 & 0 & 0 & 0 & 0 & p  \\
 0 & 0 & 0 & 0 & 0 & 0 & q  \\
 0 & 0 & 0 & 0 & 0 & 0 & r  \\
 0 & 0 & 0 & 0 & 0 & 0 & 0
\end{array}
\right), \;0< \alpha \leq \pi/2,
\]
$$
H=\left(
\begin{array}{ccccccc}
 0 & 0 & 0 & 0 & 0 & 0 & 1 \\
 0 & 0 & 0 & 0 & 1 & 0 & 0 \\
 0 & 0 & 1 & 0 & 0 & 0 & 0  \\
 0 & 0 & 0 & 1 & 0 & 0 & 0  \\
 0 & 1 & 0 & 0 & 0 & 0 & 0  \\
 0 & 0 & 0 & 0 & 0 & 1 & 0  \\
 1 & 0 & 0 & 0 & 0 & 0 & 0
\end{array}
\right).
\]
Therefore,
the conditions of the $H$-normality of $N$ are as follows:
\begin{eqnarray*}
   a & = & \overline{q}cos \alpha \\
   0 & = & \overline{q}sin \alpha \\
  b\:cos \alpha +c\: sin \alpha & = & \overline{h} \\
  2{\cal R}e \{ a\overline{d}\}+|b|^{2}+|c|^{2}+|e|^{2} & = & 2{\cal R}e \{g\overline{q}\}+|h|^{2}+|p|^{2}+|r|^{2}.
\end{eqnarray*}
Since $sin \alpha \neq 0$, $q=0$, hence $a=0$. Thus,
$(N-\lambda I)v_{2}=(N^{[*]}-\overline{\lambda}I)v_{2}=0$ which
contradicts the condition $S_{0} \cap S=\{0\}$.

Thus, we have classified all indecomposable operators with one-dimen\-sional
subspace $S_{0}$. Now let us consider the case when $dim \: S_{0}=2$.

\subsection{$dim\:S_{0}=2$}
Let $S_{0}$ be 2-dimensional. Since the operator $H_{1}=H|_{S}$ has only
positive eigenvalues,
one can assume that $H_{1}=I$. $N_{1}$ is a usual normal operator having
the only eigenvalue $\lambda$, hence, $N_{1}=\lambda I$.
As a result, we have
\begin{equation}
N=\left( \begin{array}{ccc}
	\lambda I & N_{1} & N_{2} \\
	0 & \lambda I & N_{3} \\
	0 & 0 & \lambda I \end{array}
	\right), \label{lastref_n}
\end{equation}
\begin{equation}
H=\left( \begin{array}{ccc}
	0 & 0 & I \\
	0 & I & 0 \\
	I & 0 & 0 \end{array}
	\right). \label{lastref_h}
\end{equation}
Below we will not stipulate that the pair $\{N,H\}$ has form
\{(\ref{lastref_n}),(\ref{lastref_h})\}.

For $N$ to be $H$-normal it is necessary and sufficient to have
\begin{equation}
N_{1}N_{1}^{*}=N_{3}^{*}N_{3}. \label{h_normal}
\end{equation}

According to Theorem~1, for indecomposable operators $n \leq 8$. Let us
consider the cases $n=4,\;5,\;6,\;7,\;8$ one after another.

\subsubsection{$n=4$}
In this case $C^{4}=S_{0}\dot{+}S_{1}$,
$$ N - \lambda I=
\left( \begin{array}{cc}
	0 & N_{2} \\
	0 & 0 \end{array}
	\right)=
\left( \begin{array}{cccc}
	0 & 0 & a & b \\
	0 & 0 & c & d \\
	0 & 0 & 0 & 0 \\
	0 & 0 & 0 & 0 \end{array}
	\right).
$$
Condition (\ref{h_normal}) of the $H$-normality of $N$ does not restrict the
submatrix $N_{2}$ (its terms $a$, $b$, $c$, $d$). If $N_{2}=0$, the operator
$N$ is decomposable because the nondegenerate subspace $V=span \{v_{1},v_{3}\}$
is invariant for $N$ and $N^{[*]}$. Thus, $N_{2}$ can be either of rank $1$
or of rank $2$ ($rg \: N_{2}=1$ or $2$). 

(a) ${rg\: N_{2}=2.\:}$ Suppose an $H$-unitary transformation $T$
$$T=\left( \begin{array}{cc}
T_{1} & T_{2} \\
T_{3} & T_{4} \end{array}
\right)$$
reduces $N- \lambda I$ to the form $\tilde{N}-\lambda I$:
$$
N-\lambda I=\left( \begin{array}{cc}
0 & N_{2} \\
0 & 0 \end{array}
\right),\;\;\;
\tilde{N}-\lambda I=\left( \begin{array}{cc}
0 & \widetilde{N_{2}} \\
0 & 0 \end{array}
\right).
$$
Then conditions (\ref{eqn1}) - (\ref{eqn3}) must be satisfied:
\begin{eqnarray}
  N_{2}T_{3} & = & 0  \label{eqn1} \\
  N_{2}T_{4} & = & T_{1} \widetilde{N_{2}} \label{eqn2} \\
  0 & = & T_{3} \widetilde{N_{2}}. \label{eqn3}
\end{eqnarray}
Since $N_{2}$ is invertible, (\ref{eqn1}) holds only if
$T_{3}=0$. Hence, $T$ is $H$-unitary iff
\begin{eqnarray}
  T_{1}T_{4}^{*} & = & I  \label{hh1} \\
  T_{1}T_{2}^{*}+T_{2}T_{1}^{*} & = & 0. \label{hh2}
\end{eqnarray}
From system (\ref{hh1}) - (\ref{hh2}) it follows that without loss of
generality we can consider only block diagonal transformations
of the form $T=T_{1} \oplus T_{1}^{*-1}$ because $T_{2}$ does not figure
in equations (\ref{eqn1}) - (\ref{eqn3}).

Thus, the only condition (\ref{eqn2}) $N_{2}=T_{1}\widetilde{N_{2}}T_{1}^{*}$
must be satisfied. Applying Proposition~3 from Appendix, we obtain that
the submatrix $N_{2}$ can be reduced to one of the canonical forms
$$N_{2}=\left( \begin{array}{cc}
          z & \varrho e^{-i\pi/3}z \\
          0 & e^{i\pi/3}z \end{array}
\right), \;\;
N_{2}=\left( \begin{array}{cc}
          z_{1} & 0 \\
          0 & z_{2} \end{array}
\right), $$
where $z$, $z_{1}$, $z_{2}$, $\varrho$ ($|z|=1$, $\varrho \in \Re \geq \sqrt{3}$,
$0 \leq arg\:z < \pi$ if $\varrho> \sqrt{3}$,
$|z_{1}|=|z_{2}|=1$, $arg \:z_{1} \leq arg\:z_{2}$) are invariants.
For the latter form the operator $N$ is decomposable because the nondegenerate
subspace $V=span\{v_{1},v_{3}\}$ is invariant both for $N$ and $N^{[*]}$.
For the former we obtain the following canonical form:
$$ N - \lambda I=\left( \begin{array}{cccc}
	 0 & 0 & z & re^{-i\pi/3}z \\
	 0 & 0 & 0 & e^{i\pi/3}z  \\
	 0 & 0 & 0 & 0  \\
	 0 & 0 & 0 & 0  \end{array}
\right),\;\begin{array}{c}
|z|=1, \; r\in \Re \geq \sqrt{3}, \\
0\leq arg\:z < \pi \;\; if \; r>\sqrt{3}. \end{array} $$

(b) ${rg\: N_{2}=1.\:}$ Then
$$ N_{2}=\left( \begin{array}{cc}
	ka & kb \\
	la & lb \end{array}
	\right), \; |a|+|b| \neq 0, \; |k|+|l| \neq 0. $$
If $l\overline{a} = k \overline{b}$, then $v=bv_{3}-av_{4}\neq 0$
belongs both to $S_{0}$ and $S_{1}$, which is impossible
($S_{0} \cap S_{1}=\{0\}$). Thus, we can assume that
$l\overline{a} \neq k \overline{b}$. Taking the transformation
$T=T_{1} \oplus T_{1}^{*-1}$, where
$$ T_{1}=\left( \begin{array}{cc}
	\overline{a} & k \\
	\overline{b} & l \end{array}
	\right), $$
we obtain one more canonical form:
$$ N - \lambda I=\left( \begin{array}{cccc}
	0 & 0 & 0 & 0 \\
	0 & 0 & 1 & 0 \\
	0 & 0 & 0 & 0 \\
	0 & 0 & 0 & 0 \end{array}
	\right). $$
\begin{lemma}
If an indecomposable $H$-normal operator $N$ ($N: \: C^{4} \rightarrow C^{4}$)
has the only eigenvalue $\lambda$, $dim \:S_{0}=2$, then the pair $\{N,H\}$
is unitarily similar to one and only one of canonical pairs
\{(\ref{lemma7.1}),(\ref{lemma7.3})\},
\{(\ref{lemma7.2}),(\ref{lemma7.3})\}:
$$ N=\left( \begin{array}{cccc}
	 \lambda & 0 & z & re^{-i\pi/3}z \\
	 0 & \lambda & 0 & e^{i\pi/3}z  \\
	 0 & 0 & \lambda & 0  \\
	 0 & 0 & 0 & \lambda	\end{array}
\right),\;  \begin{array}{c} |z|=1, \; r \in \Re \geq \sqrt{3}, \\
               0\leq arg\:z <\pi \; if \; r>\sqrt{3}, \end{array} $$
$$ N=\left( \begin{array}{cccc}
	 \lambda & 0 & 0 & 0 \\
	 0 & \lambda & 1 & 0 \\
	 0 & 0 & \lambda & 0 \\
	 0 & 0 & 0 & \lambda \end{array} \right), $$
$$ H=\left( \begin{array}{cc}
		0 & I_{2} \\
		I_{2} & 0 \end{array}
	\right),  $$
where  $r,z$ are $H$-unitary invariants.
\end{lemma}
{\bf{Proof:}} The possibility to reduce $N$ to one of forms (\ref{lemma7.1}),
(\ref{lemma7.2}) is proved before the lemma. The argument in (a) above shows
that these forms are not similar, hence, they are not $H$-unitarily similar.
Thus, we must only prove the indecomposability of $N$.

Show that the first canonical form is indecomposable. Assume the converse.
Let some nondegenerate subspace $V$ be invariant for $N$ and $N^{[*]}$. Then
there exists a nonzero vector $w_{1} \in V: \; w_{1} \in S_{0}$. Therefore,
$\exists w_{2}=av_{3}+bv_{4}+v \in V$ ($v \in S_{0}$, $|a|+|b| \neq 0$).
\begin{eqnarray*}
(N-\lambda I)w_{2} & = & azv_{1}+b(r e^{-i\pi/3}zv_{1}+e^{i\pi/3}zv_{2}), \\
 (N^{[*]}-\overline{\lambda} I)w_{2} & = & a(\overline{z}v_{1}+re^{i\pi/3}\overline{z}v_{2})+be^{-i\pi/3}\overline{z}v_{2}.
\end{eqnarray*}
Since $min\:\{dim\: V, dim\: V^{[\perp]}\} \leq 2$, it can be assumed that
$dim V \leq 2$. As the vectors $w_{1}$ and $w_{2}$ are linearly independent,
we get $dim\:V=2$. Therefore, the vectors $(N-\lambda I)w_{2}$ and
$(N^{[*]}-\overline{\lambda}I)w_{2}$ must be linearly dependent, i.e.,
the following condition must be satisfied:
\begin{equation}
(a+bre^{-i\pi/3})(ar e^{i\pi/3}+be^{-i\pi/3})=abe^{i\pi/3}. \label{decompos1}
\end{equation}
Since (\ref{decompos1}) breaks if either $a$ or $b$ is equal to zero,
we can rewrite (\ref{decompos1}) as follows:
\begin{equation}
(\frac{a}{b})^{2}r e^{i\pi/3}+(\frac{a}{b})(e^{-i\pi/3}-e^{i\pi/3}+r^{2})+r e^{-2i\pi/3}=0. \label{decompos2}
\end{equation}
Discriminant of (\ref{decompos2}) is equal to $r^{4}-2r^{2}-3$.
Since $r^{2} \geq 3$, it is nonnegative. Therefore,
$$\frac{a}{b}=\frac{i\sqrt{3}-r^{2}\pm \sqrt{r^{4}-2r^{2}-3}}{r(1+i\sqrt{3})}.
$$
Consequently, $|\frac{a}{b}|^{2}=\frac{1}{2}(r^{2}-1 \mp \sqrt{r^{4}-2r^{2}-3})$,
therefore, $[w_{2},(N-\lambda I)w_{2}]=z|b|^{2}\overline{(|\frac{a}{b}|^{2}+\frac{a}{b}re^{i\pi/3}+e^{-i\pi/3})}=0$.
Thus, the subspace $V$ is degenerate, i.e., the operator $N$
is indecomposable.

For the second matrix $N$ we see that the vectors $(N-\lambda I)w_{2}$ and
$(N^{[*]}-\overline{\lambda}I)w_{2}$ ($w_{2}=av_{3}+bv_{4}+v$, $v \in S_{0}$)
can be linearly dependent only if $a=b=0$. Therefore, $N$ is also indecomposable.
This concludes the proof of the lemma.

\subsubsection{$n=5$}
The matrix $N -\lambda I$ has the form
$$ N - \lambda I=\left( \begin{array}{ccc}
        0 & N_{1} & N_{2} \\
        0 & 0 & N_{3} \\
        0 & 0 & 0 \end{array}
        \right)=
\left( \begin{array}{ccccc}
        0 & 0 & a & c & d \\
        0 & 0 & b & e & f \\
        0 & 0 & 0 & g & h \\
        0 & 0 & 0 & 0 & 0 \\
        0 & 0 & 0 & 0 & 0 \end{array}
        \right) \]
so that condition (\ref{h_normal}) of the $H$-normality of $N$ amounts to the system
\begin{eqnarray*}
  | a | & = & | g | \\
  a \overline{b} & = & \overline{g} h \\
  | b | & = & | h |.
\end{eqnarray*}
The latter means that $g=\overline{a}z$, $h=\overline{b}z$ ($|z|=1$).
Note that $a$ and $b$ are not equal to zero simultaneously because
otherwise $v_{3} \in S_{0}$, which is impossible.

Take the transformation $T=T_{1} \oplus I \oplus T_{1}^{*-1}$, where
$$ T_{1}=\left( \begin{array}{cc}
        a & t_{12} \\
        b & t_{22} \end{array}
        \right), \; at_{22} \neq bt_{12}, $$
and reduce $N- \lambda I$ to the form
$$ N- \lambda I=\left( \begin{array}{ccccc}
        0 & 0 & 1 & c' & d' \\
        0 & 0 & 0 & e' & f' \\
        0 & 0 & 0 & z & 0 \\
        0 & 0 & 0 & 0 & 0 \\
        0 & 0 & 0 & 0 & 0 \end{array}
        \right),  \;  |z|=1. $$
Now we fix the form of the submatrices $N_{1}$ and $N_{3}$ so that the
following transformations will change only the submatrix $N_{2}$. At first,
apply the transformation
\begin{equation}
T=\left( \begin{array}{ccc}
        I & T_{2} & -\frac{1}{2}T_{2}T_{2}^{*} \\
        0 & I & -T_{2}^{*} \\
        0 & 0 & I \end{array}
        \right), \label{BT}
\end{equation}
where $T_{2}^{*}=(0\;\;d')$, and reduce $N_{2}$ to the form
$$ N_{2}=\left( \begin{array}{cc}
        c'' & 0 \\
        e'' & f'' \end{array}
        \right). $$
Now let us consider two cases: $f''=0$ and $f'' \neq 0$.

(a) ${f''=0. \:}$ Then $e'' \neq 0$ because otherwise $v_{5} \in S_{0}$.
Subjecting $N - \lambda I$ to the transformation
$T=T_{1} \oplus I \oplus T_{1}^{*-1}$, where
$$ T_{1}=\left( \begin{array}{cc}
        1 & c'' \\
        0 & e'' \end{array}
        \right), $$
we get
$$ N_{2}=\left( \begin{array}{cc}
        0 & 0 \\
        1 & 0 \end{array}
        \right). $$

(b) ${f'' \neq 0. \:}$ Then one can assume that $|f''|=1$ (to this end it
is sufficient to put $\widetilde{v_{2}}=\sqrt{|f''|}v_{2}$,
$\widetilde{v_{5}}=v_{5}/\sqrt{|f''|}$). Thus, $f''=z_{1}$, $|z_{1}|=1$.

If $z_{1}^{2} \neq z$, then $N$ is decomposable. Indeed, applying
\begin{equation}
 T=\left( \begin{array}{ccc}
   T_{1} & -T_{1}T_{5}^{*} & -\frac{1}{2}T_{1}T_{5}^{*}T_{5} \\
   0 & I & T_{5} \\
   0 & 0 & T_{1}^{*-1} \end{array}
   \right), \label{apr*****}
\end{equation}
where
$$ T_{1}=\left( \begin{array}{cc}
        1 & z_{1}\overline{e''}/(1-\overline{z}z_{1}^{2}) \\
        0 & 1 \end{array}
        \right), \;\;
T_{5}=\left( \begin{array}{cc}
        0 & z_{1}^{2}\overline{e''}/(1-\overline{z}z_{1}^{2})  \end{array}
        \right), $$
we reduce $N_{2}$ to the diagonal form $N_{2}=c''' \oplus z_{1}$.
Now the nondegenerate subspace $V=span \{ v_{2},v_{5}\}$ is invariant for
$N$ and $N^{[*]}$, hence, $N$ is decomposable.

Let $z_{1}^{2}=z$. Note that if $e''=0$, then $N$ is decomposable
($V=span \{v_{2},v_{5}\}$ is nondegenerate, $NV \subseteq V$,
$N^{[*]}V \subseteq V$). Thus, $e'' \neq 0$. Taking transformation
(\ref{apr*****}) with
$$ T_{1}=\left( \begin{array}{cc}
        1 & iz_{1}c_{2}''/|e''| \\
        0 & e^{i\:arg\:e''} \end{array}
        \right), \;\;
T_{5}=\left( \begin{array}{cc}
 -z_{1}(c_{1}''+c_{2}''^{2}/|e''|^{2})/2 & iz_{1}^{2}c_{2}''/|e''|  \end{array}
   \right), $$
where $c_{1}''={\cal R}e\{c''\overline{z_{1}}\}$,
$c_{2}''={\cal I}m \{c''\overline{z_{1}}\}$, we reduce $N_{2}$ to the final
form
$$ N_{2}=\left( \begin{array}{cc}
        0 & 0 \\
        r & z_{1} \end{array}
        \right), \; r=|e''|>0. $$

\begin{lemma}
If an indecomposable $H$-normal operator $N$ ($N: \: C^{5} \rightarrow C^{5}$)
has the only eigenvalue $\lambda$, $dim \:S_{0}=2$, then the pair $\{N,H\}$
is unitarily similar to one and only one of canonical pairs
\{(\ref{lemma8.1}),(\ref{lemma8.3})\},
\{(\ref{lemma8.2}),(\ref{lemma8.3})\}:
$$ N=\left( \begin{array}{ccccc}
        \lambda & 0 & 1 & 0 & 0 \\
        0 & \lambda & 0 & 1 & 0 \\
        0 & 0 & \lambda & z & 0 \\
        0 & 0 & 0 & \lambda & 0 \\
        0 & 0 & 0 & 0 & \lambda \end{array}
\right), \; |z|=1, $$
$$ \ N=\left( \begin{array}{ccccc}
        \lambda & 0 & 1 & 0 & 0 \\
        0 & \lambda & 0 & r & z \\
        0 & 0 & \lambda & z^{2} & 0 \\
        0 & 0 & 0 & \lambda & 0 \\
        0 & 0 & 0 & 0 & \lambda \end{array}
\right), \; |z|=1, \; r \in \Re>0, $$
$$ H=\left( \begin{array}{ccc}
                0 & 0 & I_{2} \\
                0 & I_{1} & 0 \\
                I_{2} & 0 & 0 \end{array}
       \right), $$
where  $z$, $r$ are $H$-unitary invariants.
\end{lemma}
{\bf{Proof:}} The possibility to reduce $N$ to one of forms (\ref{lemma8.1}),
(\ref{lemma8.2}) is proved before the lemma. Hence, it is necessary to
show that these forms are indecomposable, are not $H$-unitarily similar
to each other and their terms $z$, $r$ are $H$-unitary invariants. These
statements may be proved as follows.

For the block triangular matrix
\begin{equation}
T=\left( \begin{array}{ccc}
   T_{1} & T_{2} & T_{3} \\
   0 & T_{4} & T_{5} \\
   0 & 0 & T_{6} \end{array}
   \right) \label{ttt}
\end{equation}
to satisfy condition (\ref{11})
$NT=T\tilde{N}$, where
$$ N - \lambda I=\left( \begin{array}{ccc}
   0 & N_{1} & N_{2} \\
   0 & 0 & N_{3} \\
   0 & 0 & 0 \end{array}
   \right), \;\;
 \tilde{N} - \lambda I=\left( \begin{array}{ccc}
   0 & \widetilde{N_{1}} & \widetilde{N_{2}} \\
   0 & 0 & \widetilde{N_{3}} \\
   0 & 0 & 0 \end{array}
   \right), $$
it is necessary and sufficient to have
\begin{eqnarray}
  N_{1} T_{4} & = & T_{1} \widetilde{N_{1}} \label{eqvvv1} \\
  N_{1} T_{5} + N_{2} T_{6} & = & T_{1} \widetilde{N_{2}} + T_{2} \widetilde{N_{3}} \label{eqvvv2} \\
  N_{3} T_{6} & = & T_{4} \widetilde{N_{3}} \label{eqvvv3}.
\end{eqnarray}
If $H$ has form (\ref{lastref_h}), then for (\ref{ttt}) to be $H$-unitary
it is necessary and sufficient to have
\begin{eqnarray}
  T_{1} T_{6}^{*} & = & I \label{unitar1}  \\
  T_{1} T_{5}^{*} + T_{2} T_{4}^{*} & = & 0 \label{unitar2} \\
  T_{1} T_{3}^{*} + T_{2} T_{2}^{*} + T_{3} T_{1}^{*} & = & 0 \label{unitar3} \\
  T_{4} T_{4}^{*} & = & I. \label{unitar4}
\end{eqnarray}
If an $H$-unitary transformation $T$ reduces matrix (\ref{lemma8.2}) (the second) to
form (\ref{lemma8.1}) (the first), then from
Corollary of Proposition~1 it follows that $T$ has block form (\ref{ttt})
and, according to (\ref{11}),
\begin{equation}
T_{1}=\left( \begin{array}{cc}
   t_{11} & t_{12} \\
   0 & t_{22} \end{array}
   \right).   \label{apr&}
\end{equation}
Apply condition (\ref{eqvvv2}), replacing $T_{6}$ by $T_{1}^{*-1}$
(\ref{unitar1}), $T_{2}$ by $-T_{1}T_{5}^{*}T_{4}$ (\ref{unitar2}).
Then we get: $z/\overline{t_{22}}=0$. This contradiction proves that
the canonical forms are not $H$-unitarily similar.

If
$$ \left( \begin{array}{ccccc}
        0 & 0 & 1 & 0 & 0 \\
        0 & 0 & 0 & r & z \\
        0 & 0 & 0 & z^{2} & 0 \\
        0 & 0 & 0 & 0 & 0 \\
        0 & 0 & 0 & 0 & \lambda \end{array}
\right)T=
T \left( \begin{array}{ccccc}
        0 & 0 & 1 & 0 & 0 \\
        0 & 0 & 0 & \tilde{r} & \tilde{z} \\
        0 & 0 & 0 & \tilde{z}^{2} & 0 \\
        0 & 0 & 0 & 0 & 0 \\
        0 & 0 & 0 & 0 & \lambda \end{array} \right), $$
$|z|=|\tilde{z}|=1$, $r,\tilde{r} \in \Re>0$, then $T$ has form (\ref{ttt}),
the submatrix $T_{1}$ having form (\ref{apr&}) and $t_{11}=t_{33}$. Since
$|t_{33}|=1$ (condition (\ref{unitar4})), we can assume that $t_{11}=t_{33}=1$.
Replace $T_{6}$ by $T_{1}^{*-1}$ and apply (\ref{eqvvv3}); we have
$\tilde{z}^{2}=z^{2}$. Now substitute $T_{1}^{*-1}$ for $T_{6}$ and
$-T_{1}T_{5}^{*}$ for $T_{2}$ in (\ref{eqvvv2}). We obtain
\begin{eqnarray}
t_{35} & = & \tilde{z} t_{12}  \label{apr&&} \\
r-z\overline{t_{12}}/\overline{t_{22}} & = & \tilde{r}t_{22}-z^{2}t_{22}\overline{t_{35}}  \label{apr&&&} \\
z/\overline{t_{22}} & = & \tilde{z} t_{22}.  \label{apr&&&&}
\end{eqnarray}
From (\ref{apr&&&&}) it follows that $|t_{22}|=1$, $\tilde{z}=z$. Hence,
$1/\overline{t_{22}}=t_{22}$, $t_{35}=zt_{12}$ and
$r=\tilde{r}t_{22}$. Therefore, $r=\tilde{r}|t_{22}|$, i.e., $\tilde{r}=r$.
Thus, the numbers $z,r$ are $H$-unitary invariants of canonical
form (\ref{lemma8.2}). That $z$ is an $H$-unitary invariant of
(\ref{lemma8.1}) can be checked in the similar way.

There remains to prove that matrices (\ref{lemma8.1}) and (\ref{lemma8.2})
are indecomposable. The proof is by reductio ad absurdum. Suppose some
nondegenerate subspace $V$ is invariant for $N$ and $N^{[*]}$ ($N$ has form
(\ref{lemma8.1})). As $min\{dim \:V,dim\:V^{[\perp]}\}\leq 2$, we can assume
that
$dim\: V \leq 2$. Since there exists a vector $w_{1} \neq 0 \in S_{0}:\:
w_{1} \in V$, there exists also a vector $w_{2}=av_{3}+bv_{4}+cv_{5}+v \in V$
($v \in S_{0}$, $|b|+|c| \neq 0$). As the vectors
$(N-\lambda I)w_{2}=av_{1}+b(v_{2}+zv_{3})$ and
$(N^{[*]}-\overline{\lambda}I)w_{2}=a\overline{z}v_{1}+bv_{3}+cv_{1}$
must be linearly dependent, we obtain $b=0$. But in this case
the subspace $V$ will be degenerate because $[(N-\lambda I)w_{2},w_{2}]=0$.
This contradiction proves the indecomposability of (\ref{lemma8.1}).
Now let us check the indecomposability of (\ref{lemma8.2}). Suppose
a nondegenerate subspace $V$ is invariant both for $N$ and $N^{[*]}$.
Then, as before, $\exists w_{1}\neq 0 \in S_{0}: w_{1} \in V$ and
$ \exists w_{2}=av_{3}+bv_{4}+cv_{5}+v \in V$ ($v \in S_{0}$, $|b|+|c| \neq 0$).
Therefore, the vectors $(N-\lambda I)w_{2}-z^{2}(N^{[*]}-\overline{\lambda}I)w_{2}=
brv_{2}-crz^{2}v_{1}$ and $(N-\lambda I)w_{2}=av_{1}+brv_{2}+bz^{2}v_{3}+czv_{2}$
have to be linearly dependent. Hence, $b=0 \Rightarrow c=0$. The contradiction
obtained proves that (\ref{lemma8.2}) is also indecomposable. The proof of
the lemma is completed.
\subsubsection{$n=6$}
The matrix $N - \lambda I$ has the form
$$ N - \lambda I=\left( \begin{array}{ccc}
        0 & N_{1} & N_{2} \\
        0 & 0 & N_{3} \\
	0 & 0 & 0 \end{array}
	\right), \; \mbox{where} \;\;
N_{1}=\left( \begin{array}{cc}
        a & b \\
	c & d \end{array}
	\right). \]
The submatrix $N_{1}$ is not equal to zero because then condition
(\ref{h_normal}) of the $H$-normality of $N$ implies $N_{3}=0$ so that
$v_{3},v_{4} \in S_{0}$, which is impossible. Thus, we must consider two
alternatives: $rg\:N_{1}=2$ and $rg\: N_{1}=1$.

(a) ${rg\:N_{1} =2. \:}$ At first apply the transformation
$T=N_{1} \oplus I \oplus N_{1}^{*-1}$; it takes $N_{1}$ to $I$. Since
$N_{1}$ has become equal to $I$, $N_{3}$, according to (\ref{h_normal}),
has become unitary. Recall that
any unitary matrix is unitarily similar to some diagonal one with
nonzero terms of modulus $1$;
moreover, this representation is unique to within order of
diagonal terms. Thus,
$\exists U \; (UU^{*}=I): \;\; \widetilde{N_{3}}=U^{*}N_{3}U$, where
\begin{equation}
\widetilde{N_{3}}=\left( \begin{array}{cc}
        z_{1} & 0 \\
	0 & z_{2} \end{array}
	\right), \; | z_{1} |=| z_{2} |=1, \; arg\:z_{1} \leq arg\:z_{2}. \label{n333}
\end{equation}
If we subject $N - \lambda I$ to the transformation $T=U \oplus U \oplus U$,
then $N_{3}$ maps to (\ref{n333}) and $N_{1}=I$ does not change.

Note that if $z_{1} \neq z_{2}$, $N$ is decomposable. To check this it is
sufficient to reduce
\begin{equation}
N_{2}=\left( \begin{array}{cc}
e & f \\
g & h \end{array}
\right) \label{apr_n2}
\end{equation}
to the diagonal form by means of
transformation (\ref{BT}) with the submatrix
$$ T_{2}=\left( \begin{array}{cc}
0 & (\overline{g}-\overline{z_{1}}f)/(1-\overline{z_{1}}z_{2}) \\
(\overline{f}-\overline{z_{2}}g)/(1-z_{1}\overline{z_{2}}) & 0 \end{array}
\right) $$
(this transformation does not change $N_{1}$ and $N_{3}$).
Now the nondegenerate subspace
$V=span\{v_{1},v_{3},v_{5} \}$ is invariant for $N$ and $N^{[*]}$,
hence, $N$ is decomposable.

Thus, for $N$ to be indecomposable $N_{3}$ must be equal to $zI$. Show
that in case when $z=-1$ $N$ is also decomposable. Indeed, apply
the transformation
$$ T=\left( \begin{array}{ccc}
        U & -\frac{1}{2}N_{2}U & -\frac{1}{8}N_{2}N_{2}^{*}U \\
        0 & U & \frac{1}{2}N_{2}^{*}U \\
	0 & 0 & U \end{array}
	\right), $$
where $U$ is a unitary matrix reducing $N_{2}+N_{2}^{*}$ to the diagonal
form ($U$ is known to exist). Then $N_{2}$ becomes diagonal; we already
know that in this case $N$ is decomposable.

Thus, $N=zI$, $z \neq -1$. Now we will apply only transformations preserving
the submatrices $N_{1}$ and $N_{3}$. First let us take (\ref{BT}) with
$$ T_{2}=\left( \begin{array}{cc}
0 & 0 \\
\overline{f} & 0 \end{array}
\right) $$
and carry submatrix (\ref{apr_n2}) to the form
$$ N_{2}=\left( \begin{array}{cc}
e' & 0 \\
g' & h' \end{array}
\right). $$
Further, apply transformation (\ref{BT}) with
$$ T_{2}=\left( \begin{array}{cc}
t_{13} & 0 \\
0 & t_{24} \end{array}
\right), $$
where ${\cal R}e \{\overline{t_{13}}+zt_{13}\}={\cal R}e\:e'$,
${\cal R}e \{\overline{t_{24}}+zt_{24}\}={\cal R}e\:h'$
(since $z \neq -1$, these equations are solvable for any $e'$ and $h'$).
After this transformation
$$ N_{2}=\left( \begin{array}{cc}
ir_{1} & 0 \\
g' & ir_{2} \end{array}
\right). $$

One can assume that $g'=r_{3} \in \Re \geq 0$. To this end it is sufficient to
put $\widetilde{v_{2}}=e^{i\:arg\:g'}v_{2}$,
$\widetilde{v_{4}}=e^{i\:arg\:g'}v_{4}$,
$\widetilde{v_{6}}=e^{i\:arg\:g'}v_{6}$.
Now apply the transformation
$$ T=\left( \begin{array}{ccc}
        T_{1} & T_{1}T_{2} & -\frac{1}{2}T_{1}T_{2}T_{2}^{*}  \\
        0 & T_{1} & -T_{1}T_{2}^{*} \\
	0 & 0 & T_{1} \end{array}
	\right), \; \mbox{where} $$
$$ T_{1}=1/\sqrt{2}\left( \begin{array}{cc}
1 & 1 \\
-\overline{(z+1)}/|z+1| & \overline{(z+1)}/|z+1| \end{array}
\right), $$
$$ T_{2}=\frac{1}{2}\left( \begin{array}{cc}
-r_{3}/|z+1| & 0 \\
(ir_{2}-ir_{1})-r_{3}\overline{(z+1)}/|z+1| & r_{3}/|z+1| \end{array}
\right). $$
We get:
$$ N_{2}=\left( \begin{array}{cc}
ir_{1}' & 0 \\
g'' & ir_{1}' \end{array}
\right), \;\; r_{1}'=\frac{1}{2}(r_{1}+r_{2}). $$
As above, we can assume that $g'' \in R \geq 0$. For $N$ to be indecomposable
$g''$ must be nonzero so that $g''>0$. This is the final form of the matrix
$N - \lambda I$:
\begin{equation}
  N - \lambda I=\left( \begin{array}{cccccc}
        0 & 0 & 1 & 0 & ir_{1} & 0 \\
	0 & 0 & 0 & 1 & r_{2} & ir_{1} \\
	0 & 0 & 0 & 0 & z & 0 \\
	0 & 0 & 0 & 0 & 0 & z \\
	0 & 0 & 0 & 0 & 0 & 0 \\
	0 & 0 & 0 & 0 & 0 & 0	\end{array}
\right), \; \begin{array}{c} |z|=1,\;z \neq -1,
           \\ r_{1},r_{2} \in \Re,\;r_{2}>0.
           \end{array} \label{apr@}
\end{equation}

Let us show that $z$, $r_{1}$, $r_{2}$ are $H$-unitary invariants. To this
end suppose that an $H$-unitary matrix $T$ reduces (\ref{apr@}) to the form
$$ \tilde{N} - \lambda I=\left( \begin{array}{ccc}
        0 & I & \widetilde{N_{2}} \\
	0 & 0 & \tilde{z}I \\
	0 & 0 & 0 \end{array}
\right), \;
\widetilde{N_{2}}=\left( \begin{array}{cc}
        i\widetilde{r_{1}} & 0  \\
	\widetilde{r_{2}} & i\widetilde{r_{1}} \end{array}
\right), \;
\begin{array}{c} |\tilde{z}|=1,\;\tilde{z} \neq -1,
\\ \widetilde{r_{1}},\widetilde{r_{2}} \in \Re,\;\widetilde{r_{2}}>0.
\end{array} $$
By Corollary of Proposition~1, $T$ must have block triangular form (\ref{ttt}),
therefore,
systems (\ref{eqvvv1}) - (\ref{eqvvv3}) and (\ref{unitar1}) - (\ref{unitar4})
must hold. From (\ref{eqvvv1}), (\ref{unitar4}), and (\ref{unitar1})
it follows that $T_{1}=T_{4}=T_{6}=T_{6}^{*-1}$. Now from (\ref{eqvvv3})
it follows that $\tilde{z}=z$. Combining (\ref{unitar2}) and (\ref{eqvvv2}),
we get $N_{2}=T_{1}\widetilde{N_{2}}T_{1}^{*}+zT_{2}T_{1}^{*}+T_{1}T_{2}^{*}$.
If we denote
$$ T_{2}'=T_{2}T_{1}^{*}=\left( \begin{array}{cc}
        t_{11}' & t_{12}' \\
	t_{21}' & t_{22}' \end{array}
	\right) $$
and write out the general form for $2 \times 2$ unitary matrix
\begin{equation}
T_{1}=\left( \begin{array}{cc}
	\varrho s_{1} & \sqrt{1-\varrho^{2}} s_{2} \\
	\sqrt{1-\varrho^{2}} s_{3} & - \varrho \overline{s_{1}}s_{2}s_{3} \end{array}
	\right), \;\; \varrho \in [0,1],\; | s_{1} |=| s_{2} |=| s_{3} |=1, \label{2-unitar}
\end{equation}
then we obtain
\begin{eqnarray*}
 ir_{1} & = & i\widetilde{r_{1}}+\varrho\sqrt{1-\varrho^{2}}\overline{s_{1}}s_{2}\widetilde{r_{2}}+zt_{11}'+\overline{t_{11}'} \\
 ir_{1} & = & i\widetilde{r_{1}}-\varrho\sqrt{1-\varrho^{2}}\overline{s_{1}}s_{2}\widetilde{r_{2}}+zt_{22}'+\overline{t_{22}'}.
\end{eqnarray*}
Summing these equalities, we get
$$ 2ir_{1} = 2i\widetilde{r_{1}}+zt_{11}'+\overline{t_{11}'}+zt_{22}'+\overline{t_{22}'}. $$
It is easy to check that if ${\cal R}e\{zt+\overline{t}\}=0$ ($z \neq -1$), then
${\cal I}m\{zt+\overline{t}\}=0$. In our case $t_{11}'+t_{22}'$ plays the role of $t$,
therefore, we have $zt_{11}'+\overline{t_{11}'}+zt_{22}'+\overline{t_{22}'}=0$.
Hence $\widetilde{r_{1}}=r_{1}$.
Let us check that from the obtained equality $\widetilde{r_{1}}=r_{1}$ it follows
that $\widetilde{r_{2}}=r_{2}$. Indeed, $zN_{2}^{*}-N_{2}=T_{1}(z\widetilde{N_{2}}^{*}-\widetilde{N_{2}})T_{1}^{*}$.
$$ zN_{2}^{*}-N_{2}=\left( \begin{array}{cc}
	-ir_{1}(z+1) & zr_{2} \\
	-r_{2} & -ir_{1}(z+1) \end{array}
	\right); $$
the determinant of $zN_{2}^{*}-N_{2}$, which does not change the similarity,
is equal to $-r_{1}^{2}(z+1)^{2}+zr_{2}^{2}$, hence
$r_{2}^{2}=\widetilde{r_{2}}^{2}$.
Since sign of $r_{2}$ coincides with that of $\widetilde{r_{2}}$,
$\widetilde{r_{2}}=r_{2}$. The proof of the $H$-unitary invariance
of the numbers $r_{1}$, $r_{2}$ is completed.

(b) ${rg\:N_{1}=1. \:}$ Let us show that in this case
$N$ is decomposable. In fact,
$$ N_{1}=\left( \begin{array}{cc}
        ka & kb \\
	la & lb \end{array}
	\right), \; |a|+|b| \neq 0, \; |k|+|l| \neq 0. $$
Taking $T=T_{1}\oplus I \oplus T_{1}^{*-1}$, where
$$ T_{1}=\left( \begin{array}{cc}
        t_{11} & k \\
	t_{21} & l \end{array}
	\right), \; \;lt_{11} \neq kt_{21}, $$
we reduce $N_{1}$ to the form
$$ N_{1}=\left( \begin{array}{cc}
	0 & 0 \\
	a & b \end{array}
	\right). $$
Without loss of generality one can assume that $a \neq 0$ and, therefore,
that $a=1$ (this may be achieved by putting $\widetilde{v_{2}}=av_{2}$,
$\widetilde{v_{6}}=v_{6}/\overline{a}$). If $b \neq 0$, apply the
transformation $T_{1} \oplus T_{4} \oplus T_{1}^{*-1}$,
where
$$ T_{1}=\left( \begin{array}{cc}
	1 & 0 \\
	0 & 1/\sqrt{|b|^{2}+1} \end{array}
	\right), \;\;
T_{4}=\left( \begin{array}{cc}
1/\sqrt{|b|^{2}+1} & |b|/\sqrt{|b|^{2}+1} \\
\overline{b}/\sqrt{|b|^{2}+1} & -e^{-i\:arg\:b}/\sqrt{|b|^{2}+1} \end{array}
\right), $$
to the matrix $N - \lambda I$ (we mean that $a=1$). Then we obtain
$$ N_{1}=\left( \begin{array}{cc}
	0 & 0 \\
	1 & 0 \end{array}
	\right). $$
According to (\ref{h_normal}),
$$ N_{3}=\left( \begin{array}{cc}
0 & z_{1}cos \alpha \\
0 & z_{2}sin \alpha \end{array}
\right), \; |z_{1}|=|z_{2}|=1, \; 0 \leq \alpha \leq \pi/2. $$
Since $v_{4} \overline{\in} S_{0}$, $sin \alpha \neq 0$. Therefore,
we can apply the transformation $T$ of form (\ref{BT}), where
$$ T_{2}=\left( \begin{array}{cc}
	\overline{g} & (f-z_{1}\overline{g}cos \alpha)/(z_{2}sin \alpha) \\
	0 & 0 \end{array}
	\right) $$
($N_{2}$ has form (\ref{apr_n2})).
Under the action of $T$ the submatrices $N_{1}$ and $N_{3}$ do not change
but the submatrix $N_{2}$ becomes diagonal. Now the nondegenerate
subspace $V=span \{v_{1},v_{5}\}$ is invariant for $N$ and $N^{[*]}$,
hence, $N$ is decomposable.

\begin{lemma}
If an indecomposable $H$-normal operator $N$ ($N: \: C^{6} \rightarrow C^{6}$)
has the only eigenvalue $\lambda$, $dim \:S_{0}=2$,
then the pair $\{N,H\}$ is unitarily similar to canonical pair
\{(\ref{lemma9.1}),(\ref{lemma9.2})\}:
$$ N=\left( \begin{array}{cccccc}
        \lambda & 0 & 1 & 0 & ir_{1} & 0 \\
	0 & \lambda & 0 & 1 & r_{2} & ir_{1} \\
	0 & 0 & \lambda & 0 & z & 0 \\
	0 & 0 & 0 & \lambda & 0 & z \\
	0 & 0 & 0 & 0 & \lambda & 0 \\
	0 & 0 & 0 & 0 & 0 & \lambda	\end{array}
\right), \; \begin{array}{c} |z|=1,\;z \neq -1,
                      \\ r_{1},r_{2} \in \Re,\;r_{2}>0,
             \end{array}  $$
$$ H=\left( \begin{array}{ccc}
                0 & 0 & I_{2} \\
		0 & I_{2} & 0 \\
		I_{2} & 0 & 0 \end{array}
	\right),  $$
where $z,r_{1},r_{2}$ are  $H$-unitary invariants.
\end{lemma}

{\bf{Proof:}} It is necessary to prove only the indecomposability of
the canonical form because the rest was proved before the lemma. Suppose
that a nondegenerate subspace $V$ satisfies the conditions $NV \subseteq V$,
$N^{[*]}V \subseteq V$. As above, we can assume that $dim \: V \leq 3$
(see the proofs of the previos lemmas). Since
$\exists w_{1} \neq 0 \in S_{0}: w_{1} \in V$,
$\exists w_{2}=av_{5}+bv_{6}+v \in V$ ($v \in (S_{0}+S)$, $|a|+|b| \neq 0$).
The vectors $(N-\lambda I)(N^{[*]}-\overline{\lambda}I)w_{2}=av_{1}+bv_{2}$
and $(N-\lambda I-z(N^{[*]}-\overline{\lambda}I))w_{2}=air_{1}(1+z)v_{1}-br_{2}zv_{1}+bir_{1}(1+z)v_{2}+ar_{2}v_{2}$
must be linearly dependent because otherwise $S_{0} \subset V$ and
$dim\:V \geq 4$. Therefore, $-b^{2}r_{2}z=a^{2}r_{2}$. Since $z \neq -1$,
$a=b=0$. This contradiction proves that $N$ is indecomposable.
The proof of the lemma is completed.
\subsubsection{$n=7$}
The matrix $N - \lambda I$ has the form
$$ N- \lambda I=\left( \begin{array}{ccc}
        0 & N_{1} & N_{2} \\
        0 & 0 & N_{3} \\
        0 & 0 & 0 \end{array}
        \right),  \; \mbox{where} \;\;
N_{1}=\left( \begin{array}{ccc}
        a & b & c \\
        d & e & f \end{array}
       \right). \]
As in case when $n=6$, one can check that $N_{1} \neq 0$, therefore,
we must consider the cases $rg\:N_{1}=1$ and $rg\:N_{1}=2$. Show that the
former alternative is also impossible. Indeed, if $rg\:N_{1}=1$, then
$$ N_{1}=\left( \begin{array}{ccc}
        ka & kb & kc \\
        la & lb & lc \end{array}
        \right), \; |a|+|b|+|c| \neq 0, \; |k|+|l| \neq 0. \]
Applying the transformation $T=T_{1} \oplus I \oplus T_{1}^{*-1}$,
where
$$T_{1}= \left(\begin{array}{cc}
              t_{11} & k  \\
              t_{21} & l
              \end{array}
          \right), \;\; lt_{11} \neq kt_{21}, $$
we reduce $N_{1}$ to the form
$$ N_{1}=\left( \begin{array}{ccc}
        0  &  0 &  0  \\
        a  &  b &  c  \end{array}
        \right). $$
Then from condition (\ref{h_normal}) of the $H$-normality of $N$
it follows that
$$ N_{3}=\left( \begin{array}{cc}
        0  &  s \\
        0  &  u \\
        0  &  w \end{array}
        \right).  $$
Since there exists a nontrivial solution $\{\alpha_{i}\}_{1}^{3}$ of the system
\begin{eqnarray*}
  a \alpha_{1}+b \alpha_{2}+ c \alpha_{3}  & = & 0 \\
  \overline{s} \alpha_{1}+ \overline{u}\alpha_{2} + \overline{w}\alpha_{3} & = & 0,
\end{eqnarray*}
the nonzero vector
$v=\alpha_{1}v_{3}+\alpha_{2}v_{4}+\alpha_{3}v_{5}$ belongs to
$S_{0}$, which contradicts the condition $S_{0} \cap S=\{0\}$.

Thus, $rg \: N_{1}=2$. Then without loss of generality it can be assumed that
  $$det \left(
   \begin{array}{cc}
    a & b \\
    d & e \end{array}
    \right) \neq 0.  $$
Take the block diagonal transformation $T_{1} \oplus I \oplus T_{1}^{*-1}$, 
where
$$ T_{1}=\left(
   \begin{array}{cc}
    a & b \\
    d & e \end{array}
    \right). $$
It reduces $N_{1}$ to the form
$$ N_{1}=\left( \begin{array}{ccc}
        1  &  0 &  c'  \\
        0  &  1 &  f'  \end{array}
        \right). $$
Further, apply the transformation $T_{1} \oplus T_{2} \oplus T_{1}^{*-1}$,
where
$$ T_{1}=\left( \begin{array}{cc}
        1 & 0 \\
        0 & \sqrt{1+|f'|^{2}} \end{array}
        \right), \;
T_{2}=\left( \begin{array}{ccc}
        1 & 0 & 0 \\
        0 & 1/ \sqrt{1+|f'|^{2}} & -f'/ \sqrt{1+|f'|^{2}} \\
        0 & \overline{f'}/ \sqrt{1+|f'|^{2}} & 1/ \sqrt{1+|f'|^{2}} \end{array}
        \right). $$
Then we get
$$ N_{1}=\left( \begin{array}{ccc}
        1 & b'' & c'' \\
        0 & 1   & 0 \end{array}
        \right). $$
Now take $T=T_{1} \oplus T_{2} \oplus T_{1}^{*-1}$, where
$$ T_{1}=\left( \begin{array}{cc}
        \sqrt{1+|c''|^{2}} & b'' \\
        0 & 1 \end{array}
        \right), \;
T_{2}=\left( \begin{array}{ccc}
        1/ \sqrt{1+|c''|^{2}} & 0 & -c''/ \sqrt{1+|c''|^{2}} \\
        0 & 1 & 0 \\
        \overline{c''}/ \sqrt{1+|c''|^{2}} & 0 & 1/ \sqrt{1+|c''|^{2}} \end{array}
        \right), $$
and get the final form of the submatrix $N_{1}$:
$$ N_{1}=\left( \begin{array}{ccc}
        1 & 0 & 0  \\
        0 & 1 & 0  \end{array}
        \right). $$
Now consider the submatrix
$$ N_{3}=\left( \begin{array}{cc}
        r & s \\
        t & u \\
        v & w \end{array}  \right). $$
If $v$ and $w$ are both equal to zero, then $v_{5} \in S_{0}$. Therefore,
we can assume that $|v|^{2}+|w|^{2} \neq 0$ and can apply the transformation
$T=T_{1} \oplus T_{1} \oplus I \oplus T_{1}$, where
$$ T_{1}=\left( \begin{array}{cc}
        w/\sqrt{|v|^{2}+|w|^{2}} & \overline{v}/\sqrt{|v|^{2}+|w|^{2}}  \\
       -v/\sqrt{|v|^{2}+|w|^{2}} & \overline{w}/\sqrt{|v|^{2}+|w|^{2}} \end{array}
        \right).  $$
Then
$$ N_{3}=\left( \begin{array}{cc}
        r' & s'  \\
        t' & u' \\
        0 & w' \end{array}
        \right), \; w'=\sqrt{|v|^{2}+|w|^{2}}>0. $$
If $s' \neq 0$, replace $s'$ by $|s'|$ by putting
$\widetilde{v_{1}}=e^{i\:arg\:s'}v_{1}$,
$\widetilde{v_{3}}=e^{i\:arg\:s'}v_{3}$,
$\widetilde{v_{6}}=e^{i\:arg\:s'}v_{6}$. If $s'=0$, then apply the
transformation $\widetilde{v_{1}}=e^{-i\:arg\:t'}v_{1}$,
$\widetilde{v_{3}}=e^{-i\:arg\:t'}v_{3}$,
$\widetilde{v_{6}}=e^{-i\:arg\:t'}v_{6}$ and replace $t'$ by $|t'|$.
Now we can assume that $s' \in \Re \geq 0$ and if $s'=0$, then
$t' \in \Re \geq 0$.

Now let us apply condition (\ref{h_normal}) of the $H$-normality of $N$.
We obtain:
$$N_{3}= \left( \begin{array}{cc}
-z_{1} \overline{z_{2}} cos \alpha & sin \alpha cos \beta  \\
z_{1} sin \alpha  & z_{2} cos \alpha cos \beta  \\
0 & sin \beta \end{array}
\right), $$
$|z_{1}|=|z_{2}|=1$, $0\leq \alpha, \beta \leq \pi /2$, $\beta \neq 0$,
$z_{1}=1$ if $sin \alpha cos \beta=0$, $z_{2}=1$ if $\alpha=\pi/2$.
Let us show that in case when $\alpha=0$ $N$ is decomposable. Indeed, under the
action of (\ref{BT}), where
$$T_{2}=\left( \begin{array}{ccc}
0 & \overline{p} & (h - \overline{p}z_{2}cos \alpha cos \beta)/sin \beta  \\
0 & 0 & 0  \end{array}
\right), $$
the submatrix
$$N_{2}=\left( \begin{array}{cc}
g & h  \\
p & q  \end{array}
\right) $$
becomes diagonal. The nondegenerate subspace $V=span \{ v_{1},v_{3},v_{6}\}$
is now invariant for $N$ and $N^{[*]}$, hence, $N$ is decomposable.

Thus, $\alpha \neq 0$. Applying transformation (\ref{BT}) with
$$T_{2}=\left( \begin{array}{ccc}
0 & t_{14} & t_{15}  \\
0 & t_{24} & t_{25}  \end{array}
\right), $$
where
\begin{eqnarray*}
  t_{14} & = & g/(z_{1}sin \alpha) \\
  t_{15} & = & (h-t_{14}z_{2}cos \alpha cos \beta)/sin \beta \\
  t_{24} & = & (p-\overline{t_{14}})/(z_{1}sin \alpha)   \\
  t_{25} & = & (q-\overline{t_{24}}-t_{24}z_{2}cos \alpha cos \beta)/\sin \beta,
\end{eqnarray*}
we reduce $N_{2}$ to zero without changing $N_{1}$ and $N_{3}$. This is the
final form of the matrix $N - \lambda I$:
$$ N- \lambda I=\left( \begin{array}{ccccccc}
        0 & 0 & 1 & 0 & 0 & 0 & 0 \\
        0 & 0 & 0 & 1 & 0 & 0 & 0 \\
        0 & 0 & 0 & 0 & 0 & -z_{1} \overline{z_{2}}cos \alpha & sin \alpha cos \beta \\
        0 & 0 & 0 & 0 & 0 & z_{1} sin \alpha & z_{2}cos \alpha cos \beta  \\
        0 & 0 & 0 & 0 & 0 & 0 & sin \beta \\
        0 & 0 & 0 & 0 & 0 & 0 & 0 \\
        0 & 0 & 0 & 0 & 0 & 0 & 0 \end{array}
        \right), $$
\begin{center}
$ |z_{1}|=|z_{2}|=1, \; 0< \alpha, \beta \leq \pi/2$, $z_{1}=1$ if $\beta=\pi/2$,
$z_{2}=1$ if $\alpha=\pi/2$.
\end{center}

Show that $z_{1}$, $z_{2}$, $\alpha$, $\beta$ are $H$-unitary invariants.
Suppose an $H$-unitary matrix $T$ reduces $N - \lambda I$ to the form
$$ \tilde{N}- \lambda I=\left( \begin{array}{ccc}
        0 & N_{1} & 0 \\
        0 & 0 & \widetilde{N_{3}} \\
        0 & 0 & 0 \end{array}
        \right), \;\; \mbox{where} $$
$$N_{1}=\left( \begin{array}{ccc}
        1 & 0 & 0 \\
        0 & 1 & 0 \end{array}
        \right), \;\;\;
\widetilde{N_{3}}=\left( \begin{array}{cc}
-\widetilde{z_{1}} \overline{\widetilde{z_{2}}}cos \tilde{\alpha} & sin \tilde{\alpha} cos \tilde{\beta} \\
 \widetilde{z_{1}} sin \tilde{\alpha} & \widetilde{z_{2}}cos \tilde{\alpha} cos \tilde{\beta}  \\
  0 & sin \tilde{\beta} \end{array}
        \right), $$
$ |\widetilde{z_{1}}|=|\widetilde{z_{2}}|=1$,
$0< \tilde{\alpha}, \tilde{\beta} \leq \pi/2$, $\widetilde{z_{1}}=1$ if
$\tilde{\beta}=\pi/2$, $\widetilde{z_{2}}=1$ if $\tilde{\alpha}=\pi/2$.
Therefore, $T$ has block triangular form (\ref{ttt}) and conditions
(\ref{eqvvv1}) - (\ref{unitar4}) hold. Combining (\ref{eqvvv1}),
(\ref{unitar4}), and (\ref{unitar1}), we get: $T_{4}=T_{1} \oplus t_{55}$
($|t_{55}|=1$), $T_{1}=T_{6}=T_{6}^{*-1}$. Now from (\ref{eqvvv3}) it follows
that $T_{4}=t_{11} \oplus t_{22}$ ($|t_{11}|=|t_{22}|=1$),
\begin{eqnarray*}
t_{22} sin\alpha cos\beta & = & t_{11} sin\tilde{\alpha} cos\tilde{\beta} \\
t_{11}z_{1}sin \alpha & = & t_{22}\widetilde{z_{1}} sin \tilde{\alpha} \\
t_{22}sin \beta & = & t_{55} sin \tilde{\beta},
\end{eqnarray*}
hence $t_{11}=t_{22}=t_{55}$, hence $N_{3}=\widetilde{N_{3}}$, i.e.,
$\tilde{\alpha}=\alpha$, $\tilde{\beta}=\beta$,
$\widetilde{z_{1}}=z_{1}$, $\widetilde{z_{2}}=z_{2}$.
Thus, $\alpha$, $\beta$, $z_{1}$, $z_{2}$ are $H$-unitary invariants.

\begin{lemma}
If an indecomposable $H$-normal operator $N$ ($N: \: C^{7} \rightarrow C^{7}$)
has the only eigenvalue $\lambda$, $dim \:S_{0}=2$,
then the pair $\{N,H\}$ is unitarily similar to canonical pair
\{(\ref{lemma10.1}),(\ref{lemma10.2})\}:
$$ N=\left( \begin{array}{ccccccc}
        \lambda & 0 & 1 & 0 & 0 & 0 & 0 \\
        0 & \lambda & 0 & 1 & 0 & 0 & 0 \\
        0 & 0 & \lambda & 0 & 0 & -z_{1} \overline{z_{2}}cos \alpha & sin \alpha cos \beta \\
        0 & 0 & 0 & \lambda & 0 & z_{1} sin \alpha & z_{2}cos \alpha cos \beta  \\
        0 & 0 & 0 & 0 & \lambda & 0 & sin \beta \\
        0 & 0 & 0 & 0 & 0 & \lambda & 0 \\
        0 & 0 & 0 & 0 & 0 & 0 & \lambda \end{array}
        \right), $$
$$
|z_{1}|=|z_{2}|=1, \; 0<\alpha, \beta\leq \pi/2, \;
 z_{1}=1 \; if \; \beta=\pi/2, \; z_{2}=1\; if \; \alpha=\pi/2.
$$
$$
H=\left( \begin{array}{ccc}
        0 & 0 & I_{2} \\
        0 & I_{3} & 0 \\
        I_{2} & 0 & 0 \end{array}
        \right), $$
where $z_{1},\:z_{2},\:r,\: \alpha,\: \beta$ are $H$-unitary invariants.
\end{lemma}

{\bf{Proof:}} We have to prove only the indecomposability of the canonical
form because the rest was proved above. The proof, as is customary,
is by inductio ad absurdum. Suppose a nondegenerate subspace $V$ is
invariant for $N$ and $N^{[*]}$; then we can assume (see the proofs of the
previous lemmas) that $dim\:V \leq 3$ and $\exists w_{2}=av_{6}+bv_{7}+v \in V$
($v \in (S_{0}+S)$, $|a|+|b| \neq 0$). Then some nontrivial linear combination
of the vectors $(N^{[*]}-\overline{\lambda}I)w_{2}=av_{3}+bv_{4}+v'$
($v' \in S_{0}$) and
$(N-\lambda I)w_{2}=a(-z_{1}\overline{z_{2}}cos \alpha v_{3}+z_{1}sin \alpha v_{4})+b(sin \alpha cos \beta v_{3}+z_{2}cos \alpha cos \beta v_{4} +sin \beta v_{5})+v''$
($v'' \in S_{0}$) must belong to $S_{0}$. This implies $b=0 \Rightarrow a=0$.
The contradiction obtained proves that $N$ is indecomposable. The proof is
completed.
\subsubsection{$n=8$}
In this case
$$N- \lambda I=\left( \begin{array}{ccc}
        0 & N_{1} & N_{2} \\
        0 & 0 & N_{3} \\
        0 & 0 & 0 \end{array}
        \right), \; \mbox{where} \; \;
N_{1}=\left( \begin{array}{cccc}
        a & b & c & d \\
        e & f & g & h \end{array}
        \right). \]
As in case when $n=7$, one can check that for the condition
$S \cap S_{0}=\{0\}$ to hold the rank of $N_{1}$ must be equal to $2$.
Without loss of generality it can be assumed that
$$ det\left( \begin{array}{cc}
		a  &  b \\
		e  &  f \end{array}
		\right) \neq 0. \]
As before (in case when $n=7$), taking the block diagonal transformation
$T=T_{1} \oplus I \oplus T_{1}^{*-1}$, where
$$ T_{1}=\left( \begin{array}{cc}
		a  &  b  \\
		e  &  f  \end{array}
		\right), \]
we reduce $N_{1}$ to the form
$$N_{1}=\left( \begin{array}{cccc}
        1 & 0 & c' & d' \\
        0 & 1 & g' & h' \end{array}
        \right). $$

The results for the previous case $n=7$ let reduce the submatrix $N_{1}$
to the form $(I \;\; 0)$. Indeed, there exists a transformation
$$ T=T_{1} \oplus T_{2} \oplus T_{1}^{*-1}, \;{\mbox{where}}\;\;
T_{2}=T_{2}^{*-1}=\left( \begin{array}{cccc}
		t_{33} & t_{34} & t_{35} & 0  \\
		t_{43} & t_{44} & t_{45} & 0 \\
		t_{53} & t_{54} & t_{55} & 0 \\
		0      &      0 &      0 & 1  \end{array}
		\right), $$
that reduces the submatrix $N_{1}$ to the form
$$ N_{1}=\left( \begin{array}{cccc}
		1  &  0 &  0 &  d' \\
		0  &  1 &  0 &  h' \end{array}
		\right)  $$
and there exists a transformation
$$ T=T_{1} \oplus T_{2} \oplus T_{1}^{*-1}, \;{\mbox{where}}\;\;
T_{2}=T_{2}^{*-1}=\left( \begin{array}{cccc}
		t_{33} & t_{34} & 0 & t_{36}  \\
		t_{43} & t_{44} & 0 & t_{46} \\
		0      &  0     & 1 & 0      \\
		t_{63} & t_{64} & 0 & t_{66}  \end{array}
		\right), $$
that reduces the obtained submatrix $N_{1}$ to the desired form
\begin{equation}
 N_{1}=\left( \begin{array}{cccc}
		1  &  0 &  0 &  0 \\
		0  &  1 &  0 &  0 \end{array}
		\right). \label{apr_n1}
 \end{equation}

Now consider the submatrix $N_{3}$ and its submatrices $N_{3}'$ and $N_{3}''$:
$$ N_{3}=\left( \begin{array}{c}
		N_{3}' \\
		N_{3}'' \end{array}
		\right), \;\;
 N_{3}'=\left( \begin{array}{cc}
		p & q \\
		r & s \end{array}
		\right), \;
 N_{3}''=\left( \begin{array}{cc}
		t & u \\
		v & w \end{array}
		\right).  $$
Note that $N_{3}''$ must be nondegenerate because otherwise the system
\begin{eqnarray*}
  \overline{t} \alpha_{1}+\overline{v} \alpha_{2} & = & 0  \\
  \overline{u} \alpha_{1}+ \overline{w}\alpha_{2} & = & 0
\end{eqnarray*}
has a nontrivial solution $\{\alpha_{i}\}_{1}^{2}$, hence, the nonzero
vector $v=\alpha_{1}v_{5}+\alpha_{2}v_{6}$ belongs to $S_{0}$.

Thus, $N_{3}''$ is nondegenerate. Recall that any nondegenerate matrix is a
product of some selfadjoint positive definite matrix and some unitary one.
Consequently, $N_{3}''=RU$, where $R$ is selfadjoint positive definite and
$U$ is unitary. Let $U_{1}$ be a unitary matrix reducing $R$ to the
real positive diagonal form. Taking
$T=U^{*}U_{1} \oplus U^{*} U_{1} \oplus U_{1} \oplus U^{*} U_{1}$,
we carry $N_{3}''$  into the form
$$ N_{3}''=\left( \begin{array}{cc}
r_{1} & 0 \\
0 & r_{2} \end{array}
\right), \; r_{1},r_{2} \in \Re, \; 0<r_{1}\leq r_{2}  $$
without changing the submatrix $N_{1}$. Now we have
$$ N_{3}=\left( \begin{array}{c}
		N_{3}' \\
		N_{3}'' \end{array}
		\right)=
\left( \begin{array}{cc}
                p' & q' \\
		r' & s' \\
                r_{1} & 0 \\
		0 & r_{2} \end{array}
		\right). $$
Further, apply  transformation (\ref{BT}) with
$$ T_{2}=\left( \begin{array}{cccc}
0  &  \overline{m} & (k-r'\overline{m})/r_{1} & (l-s'\overline{m})/r_{2} \\
0  &  0 &  0 &  n/r_{2} \end{array} \right)  $$
and reduce  the submatrix
$$ N_{2}=\left( \begin{array}{cc}
                k & l \\
		m & n \end{array}
		\right) $$
to zero. Finally apply condition (\ref{h_normal}) of the $H$-normality of
$N$. We get: $r_{2} \leq 1$. Show that if $r_{1}=r_{2}$, then
$N$ is decomposable. In fact, if $r_{1}=r_{2}=1$, then from (\ref{h_normal})
it follows that $N_{3}'=0$, hence, the nondegenerate subspace
$V=span\{v_{1},v_{3},v_{5},v_{7}\}$ is invariant for $N$ and $N^{[*]}$,
hence, $N$ is decomposable. If $r_{1}=r_{2} < 1$, then the matrix
$N_{3}'/\sqrt{1-r_{1}^{2}}$ is unitary, therefore, there exists a
 unitary matrix $U$ that reduces $N_{3}'$ to the diagonal form. Then the
transformation $T=U\oplus U\oplus U\oplus U$ does not change the submatrices
$N_{1}=(I \;\; 0)$, $N_{2}=0$, $N_{3}''=r_{1}I$ and reduces $N_{3}'$ to the
diagonal form. Now it is seen that $N$ is decomposable
($V=span \{v_{1},v_{3},v_{5},v_{7}\}$ is nondegenerate, $NV \subseteq V$,
$N^{[*]}V \subseteq V$). Thus, in either case $N$ is decomposable.

There remains to consider the case when $r_{1} < r_{2}$. If $q'\neq 0$,
let us replace $q'$ by $|q'|$ by means of the transformation
$\widetilde{v_{1}}=e^{i\:arg\:q'}v_{1}$,
$\widetilde{v_{3}}=e^{i\:arg\:q'}v_{3}$,
$\widetilde{v_{5}}=e^{i\:arg\:q'}v_{5}$,
$\widetilde{v_{7}}=e^{i\:arg\:q'}v_{7}$.
If $q'=0$, let us put
$\widetilde{v_{1}}=e^{-i\:arg\:r'}v_{1}$,
$\widetilde{v_{3}}=e^{-i\:arg\:r'}v_{3}$,
$\widetilde{v_{5}}=e^{-i\:arg\:r'}v_{5}$,
$\widetilde{v_{7}}=e^{-i\:arg\:r'}v_{7}$.
Then $r'$ will be replaced by $|r'|$. Thus, one can assume that $q' \in \Re \geq 0$
and if $q'=0$, then $r'\in \Re \geq 0$. Applying (\ref{h_normal}) and
renaming the terms of $N_{3}$, we get
\begin{equation}
N_{3}= \left( \begin{array}{cc}
-z_{1}\overline{z_{2}}sin \alpha\/ cos \beta & cos \alpha\/ cos \gamma \\
z_{1}cos \alpha\/ cos \beta & z_{2}sin \alpha \/ cos \gamma  \\
sin \beta & 0  \\
0 & sin \gamma \end{array}
\right),  \label{apr_n3}
\end{equation}
$|z_{1}|=|z_{2}|=1$, $0 < \beta < \gamma \leq \pi/2$, $0 \leq \alpha \leq \pi/2$,
$z_{1}=1$ if $cos \alpha cos \gamma=0$, $z_{2}=1$ if $\alpha=0$.
We already know that if $N_{3}'$ is diagonal, $N$ is decomposable.
Therefore, $\alpha \neq \pi/2$. As a result, we have:
\begin{equation}
 N- \lambda I=\left( \begin{array}{ccc}
        0 & N_{1} & 0 \\
        0 & 0 & N_{3} \\
        0 & 0 & 0 \end{array} \right), \;\;
 N_{1}=\left( \begin{array}{cccc}
        1 & 0 & 0 & 0 \\
        0 & 1 & 0 & 0 \end{array} \right), \label{apr_n0}
\end{equation}
$N_{3}$ has form (\ref{apr_n3}),
\begin{equation}
 \begin{array}{c}
 |z_{1}|=|z_{2}|=1, \; 0 < \beta < \gamma \leq \pi/2, \; 0 \leq \alpha < \pi/2, \\
 z_{1}=1 \; if \; \gamma=\pi/2, \;z_{2}=1 \; if \; \alpha=0.
\end{array} \label{apr_ogr}
\end{equation}

Check the $H$-unitary invariance of the numbers $\alpha$, $\beta$,
$\gamma$, $z_{1}$, and $z_{2}$. To this end suppose that an $H$-unitary
matrix $T$ reduces $N - \lambda I$ to the form $\tilde{N} - \lambda I$,
where $N - \lambda I$ has form (\ref{apr_n0}), (\ref{apr_n3}), (\ref{apr_ogr}),
$$ \tilde{N} -\lambda I=\left( \begin{array}{ccc}
                0 & N_{1} & 0 \\
                0 & 0 & \widetilde{N_{3}} \\
		0 & 0 & 0 \end{array}
		\right), $$
$N_{1}$ has form (\ref{apr_n1}), $N_{3}$ has form (\ref{apr_n3}),
$$\widetilde{N_{3}}= \left( \begin{array}{cc}
-\widetilde{z_{1}}\overline{\widetilde{z_{2}}}sin \tilde{\alpha}\/ cos \tilde{\beta} & cos \tilde{\alpha}\/ cos \tilde{\gamma} \\
\widetilde{z_{1}}cos \tilde{\alpha}\/ cos \tilde{\beta} & \widetilde{z_{2}}sin \tilde{\alpha} \/ cos \tilde{\gamma}  \\
sin \tilde{\beta} & 0  \\
0 & sin \tilde{\gamma} \end{array}
\right), $$
$$ \begin{array}{c}
|z_{1}|=|z_{2}|=1, \;0 < \tilde{\beta} < \tilde{\gamma} \leq \pi/2, \; 0 \leq \tilde{\alpha} < \pi/2, \\
\widetilde{z_{1}}=1 \; if \; \tilde{\gamma}=\pi/2, \;\widetilde{z_{2}}=1 \; if \; \tilde{\alpha}=0.
\end{array} $$
Then $T$ has form (\ref{ttt}) and conditions (\ref{eqvvv1}) - (\ref{unitar4})
hold. From (\ref{eqvvv1}), (\ref{unitar4}), and (\ref{unitar1}) it follows
that $T_{4}=T_{1}\oplus T_{4}'$, $T_{4}'T_{4}'^{*}=I$,
$T_{1}=T_{6}=T_{6}^{*-1}$. From (\ref{eqvvv3}) it follows that
$N_{3}''T_{1}=T_{4}'\widetilde{N_{3}''}$. Taking into account general form
(\ref{2-unitar}) of a $2 \times 2$ unitary matrix, we can check that
this equality implies $T_{4}'=T_{1}=t_{11} \oplus t_{22}$
($|t_{11}|=|t_{22}|=1$), $\tilde{\beta}=\beta$, $\tilde{\gamma}=\gamma$.
Applying (\ref{eqvvv3}) again, we get
\begin{eqnarray*}
t_{22}cos \alpha cos \gamma & = & t_{11} cos \tilde{\alpha} cos \tilde{\gamma} \\
t_{11}z_{1}cos \alpha cos \beta & = & t_{22} \widetilde{z_{1}} cos \tilde{\alpha} cos \tilde{\beta},
\end{eqnarray*}
hence $t_{11}=t_{22}$, hence $\widetilde{N_{3}}=N_{3}$, i.e.,
$\tilde{\alpha}=\alpha$, $\widetilde{z_{1}}=z_{1}$, $\widetilde{z_{2}}=z_{2}$.

\begin{lemma}
If an indecomposable $H$-normal operator $N$ ($N: \: C^{8} \rightarrow C^{8}$)
has the only eigenvalue $\lambda$, $dim \:S_{0}=2$,
then the pair $\{N,H\}$ is unitarily similar to canonical pair
\{(\ref{lemma11.1}),(\ref{lemma11.2})\}:
$$ N- \lambda I=\left( \begin{array}{cccccccc}
        0 & 0 & 1 & 0 & 0 & 0 & 0 & 0 \\
        0 & 0 & 0 & 1 & 0 & 0 & 0 & 0 \\
        0 & 0 & 0 & 0 & 0 & 0 & -z_{1}\overline{z_{2}}sin \alpha \/cos \beta & cos \alpha \/ cos \gamma \\
        0 & 0 & 0 & 0 & 0 & 0 & z_{1}cos \alpha\/ cos \beta & z_{2}sin \alpha \/ cos \gamma \\
        0 & 0 & 0 & 0 & 0 & 0 & \sin \beta & 0 \\
        0 & 0 & 0 & 0 & 0 & 0 & 0 & \sin \gamma \\
        0 & 0 & 0 & 0 & 0 & 0 & 0 & 0 \\
        0 & 0 & 0 & 0 & 0 & 0 & 0 & 0 \end{array}
        \right), $$
$$ \begin{array}{c}
|z_{1}|=|z_{2}|=1, \; 0 \leq \alpha< \pi/2, \; 0< \beta< \gamma \leq \pi/2, \\
 z_{1}=1 \;\;if\; \gamma=\pi/2, \;\; z_{2}=0 \;\;if\; \alpha=0. \end{array}
$$
$$ H=\left( \begin{array}{ccc}
        0 & 0 & I_{2} \\
        0 & I_{4} & 0 \\
        I_{2} & 0 & 0\end{array}
        \right), $$
where $z_{1},z_{2},\alpha,\beta,\gamma$ are $H$-unitary invariants.
\end{lemma}

{\bf{Proof:}} We must  prove only the indecomposability of the
canonical form. Assume the converse. Then (see the proofs of the previous lemmas)
we can assume that $dim\:V \geq 4$, $w_{2}=av_{7}+bv_{8}+v \in V$
($v \in (S_{0}+S)$, $|a|+|b| \neq 0$).
The vectors $(N-\lambda I)(N^{[*]}-\overline{\lambda}I)w_{2}=av_{1}+bv_{2}$,
$(N^{[*]}-\overline{\lambda} I)^{2}w_{2}=a(-\overline{z_{1}}z_{2}sin \alpha cos \beta v_{1}+cos \alpha cos \gamma v_{2})+
b(\overline{z_{1}}cos \alpha cos \beta v_{1}+\overline{z_{2}}sin \alpha cos \gamma v_{2})$
and $(N-\lambda I)^{2}w_{2}=a(-z_{1}\overline{z_{2}}sin\alpha cos \beta v_{1}+
z_{1}cos \alpha cos \beta v_{2})+b(cos \alpha cos \gamma v_{1}+z_{2}sin \alpha
cos \gamma v_{2})$ must be collinear because otherwise we get
$S_{0} \subset V$, but since the condition $NS_{1} \subset (S_{1}+S_{0})$
does not hold, we obtain $dim\:V>4$. Thus, let us write the conditions
of the linear dependence (if $a$ or $b$ is equal to zero, the vectors are not
collinear):
\begin{eqnarray*}
-\overline{z_{1}}z_{2}sin \alpha cos \beta +\overline{z_{1}}cos \alpha cos \beta \frac{b}{a}=
cos \alpha cos \gamma \frac{a}{b}+\overline{z_{2}}sin \alpha cos \gamma \\
-z_{1}\overline{z_{2}}sin \alpha cos \beta+cos \alpha cos \gamma \frac{b}{a}=
z_{1}cos \alpha cos \beta \frac{a}{b}+z_{2}sin \alpha cos \gamma.
\end{eqnarray*}
If we replace the last condition by its complex conjugate and subtract it
from the first, we obtain:
$$ \overline{z_{1}}cos \alpha cos \beta \frac{b}{a} -cos \alpha cos \gamma (\frac{\overline{b}}{\overline{a}})=
cos \alpha cos \gamma \frac{a}{b}-\overline{z_{1}}cos \alpha cos \beta (\frac{\overline{a}}{\overline{b}})
$$
or
$$ \overline{z_{1}}cos \alpha cos \beta \frac{|a|^{2}+|b|^{2}}{a\overline{b}}=
cos \alpha cos \gamma \frac{|a|^{2}+|b|^{2}}{\overline{a}b}.
$$
Modulus of the left hand side must be equal to that of the right hand side,
i.e., $cos \alpha cos \beta=cos \alpha cos \gamma$. Since $cos \alpha \neq 0$,
$cos \beta=cos \gamma$, hence, $\beta=\gamma$.
But for our canonical form $\beta < \gamma$. This contradiction proves
the indecomposability of the operator $N$.

We have considered all alternatives for an indecomposable operator $N$
and have obtained canonical forms for each case. Thus, we have proved
Theorem~2.

\section*{{\large\bf Appendix}\hfill\break
Canonical Forms for $2 \times 2$ Matrices under Congruence}
\begin{proposition}
   Any invertible matrix $A$ of order $2 \times 2$ is congruent to one
and only one of the following canonical forms:
\begin{equation}
A=\left( \begin{array}{cc}
          z & \varrho e^{-i\pi/3}z \\
          0 & e^{i\pi/3}z \end{array}
\right), \; |z|=1, \; \varrho \in \Re \geq \sqrt{3}, \; 0 \leq arg\:z<\pi \; if \; \varrho>\sqrt{3},  \label{canonic1}
\end{equation}
\begin{equation}
A=\left( \begin{array}{cc}
          z_{1} & 0 \\
          0 & z_{2} \end{array}
\right), \;\; |z_{1}|=1, \; |z_{2}|=1, \; arg\:z_{1} \leq arg\:z_{2}, \label{canonic2}
\end{equation}
where $z$, $z_{1}$, $z_{2}$, $\varrho$ form a complete a minimal set of invariants.
\end{proposition}
{\bf{Proof:}} Consider the matrix $A'=AA^{*-1}$. If $\widetilde{A}=TAT^{*}$, then
$\widetilde{A'}=TA'T^{-1}$ so that spectral properties of $A'$ do not
change under congruence of $A$. Reduce $A'$ to the Jordan normal form.
Since $|det\:A'|=1$, there exist three such forms:
\begin{equation}
A'=\left( \begin{array}{cc}
          x_{1} & 0 \\
          0 & x_{2} \end{array} \right),
\;\; x_{1} \neq x_{2}, \; |x_{1}x_{2}|=1, \; |x_{1}| \leq 1, \label{jordan1}
\end{equation}
\begin{equation}
A'=xI, \;\;  |x|=1, \label{jordan2}
\end{equation}
\begin{equation}
A'=\left( \begin{array}{cc}
          x & 1 \\
          0 & x \end{array} \right), \;\; |x|=1.  \label{jordan3}
\end{equation}

(a) $A'$ is reduced to form (\ref{jordan1}). Since $A=A'A^{*}$, we have
\begin{equation}
A'=\left( \begin{array}{cc}
          a & b \\
          c & d \end{array} \right)
    =\left( \begin{array}{cc}
          \overline{a}x_{1} & \overline{c}x_{1} \\
          \overline{b}x_{2} & \overline{d}x_{2}  \end{array} \right)=A'A^{*}.
 \label{AJNF}
\end{equation}
It is seen that either $b=c=0$ or $arg \:x_{1}=arg\: x_{2}$.

If $|x_{1}| < 1$, then from (\ref{AJNF}) it follows that $a=d=0$; since $A$
is invertible, $b$ and $c$ are nonzero, therefore, $arg\: x_{1}=arg\:x_{2}$.
Now let us consider  the function
$f(\varrho)=\frac{1}{2}(1-\varrho^{2}-\sqrt{(\varrho^{2}+1)(\varrho^{2}-3)})$
of the real variable $\varrho$. It monotonically decreases on the interval
$(\sqrt{3},+\infty)$, $f(\sqrt{3})=-1$, and
$\lim_{\varrho \rightarrow +\infty}f(\varrho)=-\infty$, therefore, the equation
$f(\varrho)=s$ has a root $\varrho > \sqrt{3}$ for all $s<-1$. Let $\varrho$
be a root of the equation $f(\varrho)=-|x_{2}|$ and let
$e^{i\:arg\: x_{2}}=-e^{i\pi/3}z^{2}$, where $|z|=1$, $0 \leq arg\:z < \pi$. Then
$x_{1}=\frac{1}{2}e^{i\pi/3}z^{2}(1-\varrho^{2}+\sqrt{(\varrho^{2}+1)(\varrho^{2}-3)})$,
$x_{2}=\frac{1}{2}e^{i\pi/3}z^{2}(1-\varrho^{2}-\sqrt{(\varrho^{2}+1)(\varrho^{2}-3)})$,
and from (\ref{AJNF}) it follows that
$$ A=\left( \begin{array}{cc}
0 & b \\
e^{i\pi/3}z^{2}f(\varrho)\overline{b} & 0 \end{array} \right),
\; \; b \neq 0.  $$
Now the transformation
$$ T=\left( \begin{array}{cc}
   1 & \overline{z}(e^{-i\pi/3}f(\varrho)-1)/(\overline{b}(f(\varrho)^{2}-1)) \\
   e^{2i\pi/3}\varrho f(\varrho)/(e^{i\pi/3}f(\varrho)-1) & -e^{i\pi/3}\overline{z}\varrho/(\overline{b}(f(\varrho)^{2}-1))
\end{array} \right) $$
reduces $A$ to form (\ref{canonic1}) with  $\varrho>\sqrt{3}$. The numbers $\varrho$ and $z$ cannot
be changed under congruence because the eigenvalues of $A'$ are invariants
and from the condition $e^{i\pi/3}z^{2}f(\varrho)=e^{i\pi/3}\tilde{z}^{2}f(\tilde{\varrho})$
($|z|=|\tilde{z}|=1$, $0 \leq arg\:z, arg\:\tilde{z} <\pi$,
$\varrho, \tilde{\varrho} \in \Re >\sqrt{3}$) it follows that
$\tilde{z}=z$, $\tilde{\varrho}=\varrho$.

If $|x_{1}|=1$, then from the condition $x_{1} \neq x_{2}$ it follows that
$arg\:x_{1} \neq arg\:x_{2}$, hence $b=c=0$. By taking $T=D_{2}$
one can interchange the terms $a$ and $d$ of the matrix $A$. Hence, we can assume  that
$arg\:a \leq arg\:d$. Applying the transformation
$$ T=\left( \begin{array}{cc}
1/\sqrt{|a|} & 0 \\
0 & 1/\sqrt{|d|} \end{array} \right), $$
we reduce $A$ to form (\ref{canonic2}) with $z_{1}=e^{i\:arg\:a}$,
$z_{2}=e^{i\:arg\:d}$.

To prove the invariance of $z_{1}$ and $z_{2}$ suppose that
$\widetilde{A}=TAT^{*}$, where $A=z_{1}\oplus z_{2}$,
$\widetilde{A}=\widetilde{z_{1}}\oplus \widetilde{z_{2}}$,
$|z_{1}|=|z_{2}|=|\widetilde{z_{1}}|=|\widetilde{z_{2}}|=1$,
$arg\:z_{1} \leq arg\:z_{2}$,
$arg\:\widetilde{z_{1}} \leq arg\:\widetilde{z_{2}}$. Then
\begin{eqnarray}
z_{1}|t_{11}|^{2}+z_{2}|t_{12}|^{2} &=& \widetilde{z_{1}}  \label{app1} \\
z_{1}t_{11}\overline{t_{21}}+z_{2}t_{12}\overline{t_{22}}&=&0 \label{app2} \\
z_{1}\overline{t_{11}}t_{21}+z_{2}\overline{t_{12}}t_{22}&=&0 \label{app3} \\
z_{1}|t_{21}|^{2}+z_{2}|t_{22}|^{2} &=& \widetilde{z_{2}}. \label{app4}
\end{eqnarray}
Since $t_{11}\overline{t_{21}}=-\overline{z_{1}}z_{2}t_{12}\overline{t_{22}}$
(condiiton (\ref{app2})), (\ref{app3}) holds only if
$(z_{2}^{2}-z_{1}^{2})\overline{t_{12}}t_{22}=0$. If $z_{1}^{2} \neq z_{2}^{2}$,
then $t_{12}$ must be zero because if $t_{22}=0$, then $t_{11}=0$ and,
therefore, $\widetilde{z_{1}}=z_{2}$, $\widetilde{z_{2}}=z_{1}$, which
contradicts the condition $arg\:\widetilde{z_{1}} \leq arg\: \widetilde{z_{2}}$.
Thus, $t_{12}=0$, hence, $t_{21}=0$, $\widetilde{z_{1}}=z_{1}$,
$\widetilde{z_{2}}=z_{2}$. If $z_{1}=z_{2}$, then, according to
(\ref{app1}) - (\ref{app4}),
$\widetilde{z_{1}}=z_{1}(|t_{11}|^{2}+|t_{12}|^{2})$,
$\widetilde{z_{2}}=z_{1}(|t_{21}|^{2}+|t_{22}|^{2})$, hence
$\widetilde{z_{1}}=\widetilde{z_{2}}=z_{1}=z_{2}$. If $z_{2}=-z_{1}$ and
$\overline{t_{12}}t_{22} \neq 0$, then $t_{11}\overline{t_{21}} \neq 0$
and $\widetilde{z_{1}}=z_{1}(|t_{11}|^{2}-|t_{12}|^{2})$. Since
$|t_{21}|/|t_{22}|=|t_{12}|/|t_{11}|$,
$\widetilde{z_{2}}=z_{1}(|t_{21}|^{2}-|t_{22}|^{2})=-\widetilde{z_{1}}|t_{22}|^{2}/|t_{11}|^{2}$.
As $arg\:\widetilde{z_{1}} \leq arg\:\widetilde{z_{2}}$, we get
$\widetilde{z_{1}}=z_{1}$, $\widetilde{z_{2}}=z_{2}$. The case when  $z_{2}=-z_{1}$
and $\overline{t_{12}}t_{22}=0$ can be considered as before. Thus, we have
proved the invariance of the numbers $z_{1}$ and $z_{2}$.

(b) $A'$ is reduced to form (\ref{jordan2}). Then $A=xA^{*}$, $|x|=1$, this
property being invariant with respect to congruence. Since $A$ is invertible,
$A=RU$, where $R$ is selfadjoint positive definite matrix and $U$ is unitary.
Let $T$ be a unitary matrix reducing $U$  to the diagonal form $\Lambda$.
After the application of $T$ we have: $A=\tilde{R}\Lambda$, where $\widetilde{R}=TRT^{*}$
is also selfadjoint positive definite. Now let $T$ be a lowertriangular
matrix such that $T\widetilde{R}T^{*}=I$. Then we reduce $A$ to the
uppertriangular form $T^{*-1} \Lambda T^{*}$. Since the term $c$ of $A$
is now equal to zero, from the condition $A=xA^{*}$ it follows that
$b$ is also equal to zero, i.e., $A$ is diagonal. We already know that
a diagonal matrix is congruent to (\ref{canonic2}) (see case (a) above).
Thus, $A$ can be reduced to form (\ref{canonic2}).

(c) $A'$ is reduced to form (\ref{jordan3}). Let $x=-e^{i\pi/3}z^{2}$ ($|z|=1$).
Then the application of the condition $A=A'A^{*}$ yields:
$$ A=\left( \begin{array}{cc}
a & b \\
-e^{i\pi/3}z^{2}\overline{b} & 0  \end{array}
\right), \; \;\; b=\overline{a}+e^{-i\pi/3}\overline{z}^{2}a. $$
For $A$ to be invertible $b$ must be nonzero. Since
$|b|=|a+e^{i\pi/3}z^{2}\overline{a}|=|a\overline{z}+e^{i\pi/3}\overline{a}z|=
|a\overline{z}-e^{-2i\pi/3}\overline{a}z|=
|e^{i\pi/3}a\overline{z}-e^{-i\pi/3}\overline{a}z|=
2|{\cal I}m\{e^{i\pi/3}a\overline{z}\}|$, we see that
${\cal I}m\{e^{i\pi/3}a\overline{z}\} \neq 0$. Let us chose $z$ so that
${\cal I}m\{e^{i\pi/3}a\overline{z}\} >0$. Applying the transformation
$$ T=\frac{\sqrt[4]{3}}{\sqrt{|b|}^{3}} \left( \begin{array}{cc}
|b| & \frac{2}{3}i\overline{z}{\cal I}m\{a\overline{z}\}|b|/\overline{b} \\
e^{i\pi/3}\overline{z}\overline{b} & \overline{z}^{2}(-\frac{2}{3}i {\cal I}m\{a\overline{z}\}+a\overline{z}) \end{array}
\right), $$
we reduce $A$ to form (\ref{canonic2}) with $\varrho=\sqrt{3}$. It is clear
that matrix (\ref{canonic2}) with $\varrho=\sqrt{3}$ is not comgruent  to
that with $\varrho>3$ because in the former case $A'$ has the diagonal Jordan
normal form in contrast to the latter. Therefore, we must prove only the
invariance of $z$. Note that if $\widetilde{A}=TAT^{*}$, where
$$ A=\left( \begin{array}{cc}
          z & \sqrt{3} e^{-i\pi/3}z \\
          0 & e^{i\pi/3}z \end{array}
\right), \;\; \widetilde{A}=\left( \begin{array}{cc}
          \tilde{z} & \sqrt{3} e^{-i\pi/3}\tilde{z} \\
          0 & e^{i\pi/3}\tilde{z} \end{array}
\right), \;\;  |z|=|\tilde{z}|=1, $$
then $\tilde{z}^{2}=z^{2}$ because the eigenvalue $x=-e^{i\pi/3}z^{2}$  of $A'$
does not  change under congruence of $A$. Therefore,
$$ A'=z^{2}\left( \begin{array}{cc}
          1-3e^{i\pi/3} & \sqrt{3} \\
          \sqrt{3} & e^{2i\pi/3} \end{array}
\right)=\widetilde{A'}. $$
For $T$ to satisfy the condition $A'T=TA'$ the matrix $T$ must has the form
$$ T=\left( \begin{array}{cc}
          t_{11} & t_{12} \\
          t_{12} & t_{11}+it_{12} \end{array} \right). $$
Now from the condition $\widetilde{A}=TAT^{*}$ it follows that
\begin{eqnarray}
z|t_{11}|^{2}+\sqrt{3}e^{-i\pi/3}zt_{11}\overline{t_{12}}+e^{i\pi/3}z|t_{12}|^{2} & = &\tilde{z} \label{app5} \\
z\overline{t_{11}}t_{12}+\sqrt{3}e^{-i\pi/3}z|t_{12}|^{2}+e^{i\pi/3}z(t_{11}\overline{t_{12}}+i|t_{12}|^{2}) & = & 0. \label{app6}
\end{eqnarray}
If $t_{12} \neq 0$, from (\ref{app6}) it follows that
$$ e^{-i\pi/6}\frac{\overline{t_{11}}}{\overline{t_{12}}}+\sqrt{3}e^{-i\pi/2}+
e^{i\pi/6}\frac{t_{11}}{t_{12}}+e^{2i\pi/3}=0, $$
which is impossible because the imaginary part of the left hand side is equal
to ${\cal I}m \{\sqrt{3}e^{-i\pi/2}+e^{2i\pi/3}\}=-\sqrt{3}/2$. Therefore,
$t_{12}=0$, hence (condition (\ref{app5})) $\tilde{z}=z$, i.e., $z$ is an
invariant. This concludes the proof of the proposition.

\end{document}